\documentclass[10pt]{article}
\usepackage{graphicx}
\usepackage{amsmath,amsfonts,latexsym,amssymb}
\usepackage{maple2e}
\usepackage[latin1]{inputenc}
\usepackage[T1]{fontenc}
\usepackage{ae,aecompl}
\usepackage{dsfont}
\newtheorem{theorem}{Theorem}[subsection]
\newtheorem{lemma}[theorem]{Lemma}
\newtheorem{proposition}[theorem]{Proposition}
\newtheorem{coro}[theorem]{Corollary}
\newtheorem{definition}[theorem]{Definition}

\newcommand{\pr}{{\mathbb R\mathbb P}^1}

\newcommand{\sln}{{\rm PSL}(n,\mathbb R)}
\newcommand{\sld}{{\rm PSL}(2,\mathbb R)}

\newcommand{\hol}{\operatorname{Hol}}
\newcommand{\bb}{{\bf B}}
\newcommand{\bfb}{{\bf b}}
\newcommand{\proof}{\vskip 0.2 cm\par\noindent{\sc Proof: }}
\newcommand{\rmks}{\vskip 0.5 truecm\noindent{\sc Remarks: }}
\newcommand{\rmk}{{\vskip 0.5 truecm\noindent\sc Remark: }}
\newcommand{\qed}{~{\sc Q.e.d.}\vskip 0.2 cm}

\newcommand{\comp}{\ast}
\newcommand{\goodRepeta}{{\bf R}}
\newcommand{\goodRepmu}{{\bf M}}
\newcommand{\diagon}{{\Delta}}
\newcommand{\ddd}{\underline}
\newcommand{\invol}{{\mathrm J}}
\newcommand{\grf}{\pi_1 (\Sigma)}
\newcommand{\Hol}{\operatorname{Hol}}
\newcommand{\bgrf}{\partial_\infty \pi_1 (\Sigma)}
\newcommand{\scal}[2]{\langle #1\vert #2\rangle}
\newcommand{\bir}[4]{\frac{\langle #1\vert #2\rangle}{\langle #3\vert #4\rangle}}

\newcommand{\mapping}[4]
{
\left\{
\begin{array}{rcl}
#1 &\rightarrow& #2\\
#3 &\mapsto& #4 
\end{array}
\right.
}
\newcommand\HH{\mathbb{H}^2}
\newcommand\semic{\phi}

\newcommand\EE{\mathbb{R}^n}
\newcommand\RPN{\mathbb{P}(\mathbb{R}^{n})}
\newcommand\torus{\mathbb{T}}
\begin{document}
\title{Cross Ratios and Identities for Higher Teichmüller-Thurston Theory}
\author{François LABOURIE \thanks{Univ. Paris-Sud, Laboratoire de Mathématiques, Orsay F-91405 Cedex; CNRS, Orsay cedex, F-91405, France, ANR Repsurf : 06-BLAN-0311}\and Gregory McSHANE\thanks{Laboratoire Emile Picard, CNRS UMR 5580, Université Paul Sabatier, ANR Repsurf : 06-BLAN-0311}
}

\maketitle

\begin{abstract} We generalise the McShane-Mirzakhani identities from hyperbolic geometry to  arbitrary cross ratios.  We define and study Hitchin representations of open surface groups to ${\rm PSL}(n,\mathbb R)$. We associate to these representations cross ratios and then give explicit  expressions for our generalised identities  in terms of (a  suitable choice of) Fock--Goncharov coordinates.
\end{abstract}
\section{Introduction}

In \cite{McShane:1998}, the second author establishes an identity  for lengths of simple closed geodesics on punctured hyperbolic surfaces.
The results holds  for surfaces with multiple  cusps, however, in order
to simplify the exposition, we consider  in this section only  the case where  $\Sigma$ is  a complete connected orientable finite area  hyperbolic surface with  a single cusp.
The length  $\ell(C)$ of a homotopy class   $C$ of  closed curves is  
defined to be the infimum of the set of  lengths of  curves freely homotopic to  $C$ with respect to the hyperbolic metric; 
this naturally extends to families $\{C_i\}_{i}$ by setting $\ell(\{C_i\}_{i}) = \sum_{i} \ell(C_{i})$. We denote by $\mathcal P$ the set of  embedded pants (with marked boundary) in $\Sigma$  up to isotopy, such that the  first   boundary component of the pair of pants is the cusp.

With this notation, McShane's identity is
\begin{eqnarray}
1=\sum_{P\in\mathcal P}
\frac
{
1}
{
e^{\frac{\ell(\partial P)}{2}}+1}.\label{McShane:1998}
\end{eqnarray}

In \cite{Mirzakhani:2007a}, using McShane's method, 
M. Mirzakhani  extends  this identity to  hyperbolic surfaces  with totally geodesic boundary. 
Let $\Sigma$ be a complete  finite area connected hyperbolic surface with  a single totally geodesic  boundary component $\partial \Sigma$. Mirzakhani's identity is

\begin{eqnarray}
\ell(\partial \Sigma)=\sum_{P\in\mathcal P}\log\left(
\frac
{
e^{\frac{\ell(\partial P)}{2}}+e^{\ell(\partial \Sigma)}}
{
e^{\frac{\ell(\partial P)}{2}}+1}\right),\label{Mirzakhani:2007a}
\end{eqnarray}
where $\mathcal P$ is the set of  embedded pants (with marked boundary) up to homotopy,
 such that the  first boundary component of the pair of pants is $\partial \Sigma$.

The purpose of the paper is twofold. Firstly, 
we show that  the identity above has a natural formulation in terms of (generalised) cross ratios. 
Then, using this formulation, 
we study  identities arising from the cross ratios constructed by the first author in \cite{Labourie:2005} 
for representations  from fundamental groups of surfaces to ${\rm PSL}(n,\mathbb R)$. 

{\em In order to avoid lengthy hypotheses, we shall use the following conventions throughout this article: unless otherwise stated,  all surfaces are assumed to be compact, connected and oriented, all closed surfaces have genus at least  2, all surfaces with non empty boundary are of finite type and their double have genus at least  2.}

We now give a  brief overview of  the main ideas.
\vskip 0.2 truecm
Our first aim in this article is to show that the identities above have a more general interpretation. Let  $\Sigma$ be  a closed surface (of genus at least 2 accordinng to our convention).  The boundary at infinity $\bgrf$ of the fundamental group $\grf$ is a one dimensional compact connected Hölder manifold -- hence Hölder homeomorphic to the circle ${\mathbb T}$ -- equipped with an action of $\grf$ by Hölder homeomorphisms.

A {\em cross ratio} on $\bgrf$ is a  Hölder  function $\bb$  defined on  $$
\bgrf^{4*}=\{(x,y,z,t)\in\bgrf^4\ \vert\ x\not=t \hbox{ and } y\not= z\},$$ invariant under the diagonal action of $\grf$, and 
which satisfies some algebraic  rules. 
Roughly speaking, these rules  encode two conditions which constitute a normalisation, 
and  two multiplicative cocycle identities which hold for some of the variables:
\begin{eqnarray*}
{\hbox{\sc Normalisation:}}\ \ \bb(x,y,z,t)&=&0\ \ \Leftrightarrow x=y \,,\hbox{ or } z=t,\\ 
{\hbox{\sc Normalisation:}}\ \ \bb(x,y,z,t)&=&1 \ \ \Leftrightarrow x=z \,,\hbox{ or } y=t,\\
{\hbox{\sc Cocycle identity:}}\ \ \bb(x,y,z,t)&=&{\bb(x,y,z,w)}{\bb(x,w,z,t)},\\
{\hbox{\sc Cocycle identity:}}\ \ \bb(x,y,z,t)&=&\bb(x,y,w,t)\bb(w,y,z,t).
\end{eqnarray*}
These relations imply an essential symetry
\begin{eqnarray*}
{\hbox{\sc Symmetry:}}\ \ \bb(x,y,z,t)&=&\bb(z,t,x,y).
\end{eqnarray*}
The  {\em period} of  a nontrivial element $\gamma$ of $\grf$ with respect to $\bf B$ is the following real number
\begin{eqnarray*}
\log\vert \bb(\gamma^-,\gamma y,\gamma^+,y)\vert:=\ell_\bb(\gamma),
\end{eqnarray*}
where  $\gamma^{+}$ (respectively $\gamma^{-}$) is  the attracting (respectively  repelling) fixed point of 
$\gamma$ on $\partial_\infty\pi_1 (\Sigma)$ and $y$ is any  element of
 $\partial_\infty \pi_1 (\Sigma) \setminus \{\gamma^{+}, \gamma^{-}\}$;
 as the notation suggests, one checks immediately that  the period is independent of $y$ (see Paragraph  \ref{def:period}).   
In Section \ref{opCR} we extend  the above definitions to cover open surfaces of finite type.

These  definitions are closely related to those given by Otal in \cite{Otal:1990,Otal:1992},
 those discussed from various perspectives  by Ledrappier in \cite{Ledrappier:1995} 
and those of  Bourdon in \cite{Bourdon:1996} in the context of ${\rm CAT}(-1)$-spaces.

The archetype of a cross ratio arises in   hyperbolic geometry in the following way.
A complete hyperbolic metric on $\Sigma$ gives rise to an identification of $\bgrf$ with the real projective line. 
Thus, the projective cross ratio on the projective line gives rise to a cross ratio on $\bgrf$ -- called a {\em hyperbolic cross ratio}.
 The period of 
$\gamma$ is just the hyperbolic length of the closed geodesic associated to  $\gamma$. 
 We recall that the projective cross ratio ${\bf b}$ on the projective line
  -- which we shall refer to  in the sequel as the {\em classical cross ratio} -- is given in projective coordinates by
$$
{\bf b}(x,y,z,t)=\frac{(x-y)(z-t)}{(x-t)(z-y)},
$$
so that ${\bf b}(x,0,1,\infty)=x$.

More generally, as was observed in \cite{Labourie:2005} and as  described in the following  paragraphs, 
cross ratios are associated to certain representations from $\grf$ to ${\rm PSL}(n,\mathbb R)$ for a closed surface $\Sigma$.

Let $P $ be the -- oriented -- sphere minus three points. We choose once and for all a {\em marking} of $P $, that is three elements $\alpha_0$,$\beta_0$, and $\gamma_0$ of $\pi_1(P )$ represented by loops freely homotopic to the oriented boundary loops and such that $\alpha_0\gamma_0\beta_0=1$. An {\em isotopy  class of  pair of pants} in $\Sigma$ is an isotopy  class of embeddings  of $P $ in $\Sigma$ such that the images  $\alpha$, $\beta$ and $\gamma$ of $\alpha_0$, $\beta_0$ and $\gamma_0$ are non trivial. The triple $(\alpha,\beta,\gamma)$ is well defined up to conjugacy by elements of $\grf$ and satisfies  $\alpha\gamma\beta=1$. Henceforth, we refer to  $\alpha$ as the {\em first} boundary component of the  pair of pants.

Given a cross ratio $\bf B$ on $\bgrf$,  we define the value of   the {\em pant gap function}  $G_{\bb}$  associated to $\bf B$
of an  isotopy class of  pair of pants  $P$ represented by $(\alpha,\beta,\gamma)$ to be
$$
G_\bb(P)=\log\vert\bb(\alpha^+,\gamma^-,\alpha^-,\beta^+) \vert.
$$
Note that the first boundary component $\alpha$ plays a special role.
We shall prove 
\begin{theorem}
Let $\Sigma$ be compact  surface with exactly one boundary component representing   a non trivial element $\alpha$ in $\grf$. 
Let $\bb$ be a cross ratio on $\bgrf$. Let $\mathcal P$ be the space of isotopy  classes of pairs of pants  in $\Sigma$ whose first boundary component is $\alpha$, then
$$
\ell_\bb(\alpha)=\sum_{P\in\mathcal P}G_\bb(P).
$$
\end{theorem}
Moreover, after suitably extending  the notion of cross ratio  (see Subsection \ref{opCR}),
the theorem generalises to more general open surfaces of finite type. 
It also generalises ``at a cusp" so that  we  recover Formula (\ref{McShane:1998}). 
The complete results -- Theorem \ref{fondbound} and Theorem \ref{cuspcase} -- are  proved in Section \ref{fondid}.

Examining more closely 
the case of hyperbolic cross ratios, we observe that the pant gap function for $P$ can be computed in terms 
of the lengths of just the boundary components of the pair of pants $P$.
Indeed recall  that every hyperbolic pair of pants with totally geodesic boundary is  determined up to isometry
 by the length of its three boundary components and hence every geometric function is a function of these  three parameters. 
In Section \ref{hypcase}, using Thurston's {\em shear coordinates},  described by Bonahon in \cite{Bonahon:1996}, and  elementary manipulations involving the classical  cross ratio
  -- as opposed to hyperbolic trigonometry in the original proofs -- 
 we recover  Mirzakhani-McShane formulae (\ref{McShane:1998}) and (\ref{Mirzakhani:2007a})  for the pant gap function.

We now discuss briefly  additional ideas that allow our  approach  to extend  to higher dimensions.
For a closed surface $\Sigma$, the first author gives   in \cite{Labourie:2005} an interpretation of
{\em  Hitchin representations} --  which fill a connected component of the space of representations
 from $\grf$ to ${\rm PSL}(n,\mathbb R)$ --
 as the space of rank $n$ cross ratios  on  $\bgrf$  (see precise Definition in \cite{Labourie:2006}).   
The second aim of this article -- accomplished  in Section \ref{sec:hitchpos} -- is to extend these results to open surfaces. 

Let us be more precise.
An element in ${\rm PSL}(n,\mathbb R)$ is {\em purely loxodromic}
  if it is real split and has simple eigenvalues, or in other words all its eigenvalues are real with mutilplicity 1. Let $\Sigma$ be a compact surface possibly with  boundary. 
   A representation  from $\grf$ to ${\rm PSL}(n,\mathbb R)$ is {\em Fuchsian}
  if it factors as
a discrete faithful representation without parabolics into ${\rm PSL}(2,\mathbb R)$ 
composed with 
 the irreducible representation from ${\rm PSL}(2,\mathbb R)$ to ${\rm PSL}(n,\mathbb R)$.
A representation  from $\grf$ to ${\rm PSL}(n,\mathbb R)$ is  {\em Hitchin} if
 the boundary components have purely loxodromic images under  the representation,
and if the representation can be deformed into a Fuchsian representation so that the images
 of the boundary components stay purely loxodromic.

In Section \ref{sec:hitchpos}, we obtain  the following generalisation of results in \cite{Labourie:2005}
\begin{theorem}\label{hcintro}
Let $S$ be a compact connected orientable surface with boundary whose double has genus at least  2.   Let $\rho$ be a Hitchin representation of $\pi_1(\Sigma)$. Then the image under $\rho$ of every non trivial element of the fundamental group is purely loxodromic and, moreover, there exists a cross ratio $\bb$ on $\partial_\infty\pi_1(S)$ whose periods satisfy
$$
\forall \gamma\in \pi_1(S),\ \ \ell _\bb(\gamma)=\log\left(\left\vert\frac{\lambda_{\max}(\rho(\gamma))}{\lambda_{\min}(\rho(\gamma))}\right\vert\right),
$$
where $\lambda_{\max}(\rho(\gamma))$ and $\lambda_{\min}(\rho(\gamma))$ are the  eigenvalues  of respectively maximum and minimum absolute values of the element $\rho(\gamma)$.

Furthermore, let $\Sigma$ be an incompressible surface embedded in $S$.
Then $\rho$ restricted to $\pi_1(\Sigma)$ is a Hitchin  representation.
\end{theorem}

 Theorem \ref{hitchinbord} is a more detailed  statement of this result. 
Its proof involves a doubling construction  described in Theorem \ref{goodhomo} and
 which also yields 
 \begin{coro}
 Let $\Sigma$ be a  compact surface with boundary. Let $\rho$ be a Hitchin representation of $\pi_1(\Sigma)$. Then there exists  a {\em closed} surface $S$ containing $\Sigma$ and a Hitchin representation $\widehat\rho$ of $\pi_1(S)$ such that $\rho$  is the restriction of  $\widehat\rho$.
 \end{coro} 
  
  In a series of articles \cite{Labourie:2006,Labourie:2005a,Labourie:2005}, 
the first author shows that Hitchin representations for closed surfaces are discrete and faithful,
  that every non trivial element is purely loxodromic 
 and that the mapping class group acts properly on the moduli space of Hitchin representations.  
A  consequence of  the above corollary is that these  concepts and results
for  closed surfaces carry over  for surfaces with boundary. 
 
 \vskip 0.2 truecm  
For a cross ratio associated to a representation $\rho$  from $\grf$ to $\sln$, 
the  gap function of a pair of pants associated to  the triple $(\alpha,\beta,\gamma)$ of elements of
 $\grf$  depends only  on the triple $(\rho(\alpha),\rho(\beta),\rho(\gamma))$ of elements of $\sln$ 
(see Proposition \ref{gapsubsurface}). 
 However,  contrary to hyperbolic cross ratios,
 for $n\geq 3$ a simple computation of the  dimensions shows that the pant gap function 
is no longer determined  by just  the eigenvalues of the monodromies of the  three boundary components of the pants: it also  depends on  extra {\em internal} parameters.
In the third part of this article,  we  give a description of  the gap functions  for Hitchin representations which nevertheless 
parallels  our description of  gap functions in  hyperbolic geometry. 
We first recall the construction by Fock and Goncharov in \cite{Fock:2006a}  of {\em Fock--Goncharov coordinates} 
on the {\em Fock--Goncharov moduli space} which are far reaching generalisations of {\em Thurston's shear coordinates}
 on the {\em enhanced Teichmüller space} for rational laminations, 
again  discussed by Bonahon in \cite{Bonahon:1996}. 
 We use  these coordinates to obtain  coordinates on the space of 
Hitchin representations which we show are {\em positive} in the sense of Fock--Goncharov -- see Definition \ref{def:posrep} and 
Theorem \ref{hitchinbord}. 
  Actually, since the Fock--Goncharov moduli space is a ``covering" of the space of Hitchin representations, 
  we obtain $(n!)^3$ different sets of coordinates for the moduli of a pair of pants. 
 In Theorem \ref{theo:goodgap} we finally show that, for a suitable choice of coordinates, the pant gap function has a nice expression.
 On the other hand, using  computer algebra software
 and the explicit description of the holonomies given by V. Fock and A. Goncharov in \cite{Fock:2007}, 
 in Section \ref{Xmap} we show  that,  even when  $n=3$, 
 the pant gap function has a very complicated expression for other choices of coordinates -- see for instance Formula (\ref{comp:1}).

\vskip 0.2 truecm

In conclusion, recall that using her identities, M. Mirzakhani gives in \cite{Mirzakhani:2007a} a recursive formula for the volumes of the moduli space of hyperbolic structures 
identified with the quotient of Teichmüller space by the mapping class group. 
From the work of the first author in \cite{Labourie:2005}, it follows that the mapping class group acts properly on the moduli space of 
Hitchin representations.  It is quite possible that the formula obtained in Theorem \ref{theo:goodgap}
 combined with the use of Fock--Goncharov coordinates can help to compute geometric quantities associated to the corresponding quotient. 
 However the volume is not the right quantity to compute since for $n$ at least 3, one can show it is infinite. 
 
 We also hope that some of our work could be generalised to ${\rm PSL}(n,\mathbb C)$ as McShane's identities were generalised to ${\rm PSL}(2,\mathbb C)$ by B Bowditch \cite{Bowditch:1996},  H. Akiyoshi, H. Miyachi, M. Sakuma \cite{Akiyoshi:2006} and Ser Peow Tan, Yan Loi Wong, Ying Zhang \cite{Tan:2004}.
 
 We conclude by saying that it is a striking fact that
 so many  of the familiar ideas from the world of hyperbolic geometry translate naturally to the world of Hitchin representations.
So much so that  one is tempted to call the latter {\em  a higher (rank) Teichmüller-Thurston theory}.

We are very grateful to  the referees for many extremely  useful comments which helped a lot  in improving the first version of this manuscript.

\tableofcontents

\section{Cross ratios and the boundary at infinity}\label{cross ratiodef}

We first recall the definition of cross ratio  on the boundary at infinity of a closed surface according to \cite{Labourie:2005} and then extend this notion to open surfaces. This extension requires defining a boundary at infinity for surfaces with cusps and boundary.
\subsection{Boundary at infinity for open  surfaces of finite type}\label{opCR}

Let $\Sigma$ be a connected orientable open surface homeomorphic to 
$$\Sigma_0\setminus\{b_1, \ldots, b_n,c_1,\ldots c_p\},$$ 
where $\Sigma_0$ is closed and $b_i, c_i$ are distinct points in $\Sigma_0$.
We have deliberately labelled the points in two different ways. The points $b_i$ will be referred as {\em boundary components} 
 and the points $c_i$ as {\em cusps}. 
 We denote this data by $\Sigma(g,n,p)$ meaning that $\Sigma_0$ has genus $g$,
  and that there are $n$ boundary components and $p$ cusps. We assume that $2g-2+p+n>0$.
  \begin{definition}{\sc[admissible metric]}
  A finite volume hyperbolic metric on  $\Sigma$ is {\em admissible}
  if the completion of $\Sigma$ is a surface with $n$ totally geodesic boundary components
 (corresponding to the neighbourhoods of the  points $b_{i}$) and $p$ cusps (corresponding to neighbourhoods of the the points $c_{i}$).  \end{definition}
 
 Admissible metrics exist since $2g-2+p+n>0$. The deck transformations  associated to  an admissible metric $\sigma$  on an oriented surface yield  an embedding $\rho_\sigma$ from $\grf$ into ${\rm PSL}(2,\mathbb R)$ well defined up to conjugacy in ${\rm PSL}(2,\mathbb R)$.

 We recall the following fact
 
 \begin{proposition}{\sc[Limit sets]}\label{limset} If $\sigma_1$ and $\sigma_2$ are two admissible metrics and  let    $\rho_1$ and $\rho_2$ be the embeddings from $\grf$ associated  as above. Let  $\Lambda_1$  and $\Lambda_2$ be the the limit sets of $\rho_1(\grf)$ and $\rho_2(\grf)$ respectively. Then there exists a unique $\grf$-equivariant Hölder homeomorphism 
$\semic$ from  $\Lambda_1$ to $\Lambda_2$.
\end{proposition}

This allows us to give the following

\begin{definition}{\sc[Boundary at infinity for open surfaces]} The {\em boundary at infinity}   $\bgrf_p$ of the open surface $\Sigma$ with $p$ cusps  is 
the boundary at infinity  of the universal cover 
$\tilde{\Sigma}$  of $\Sigma$ equipped with some admissible metric $g$. 
\end{definition}

\rmks

 \begin{itemize}
 \item Let $\rho_g$ be the monodromy representation of an admissible hyperbolic metric on $\Sigma$. The universal cover $\tilde{\Sigma}$  is then  isometric to the interior of the convex hull of the limit set of the group $\rho_g(\grf)$ in the hyperbolic plane. It follows that  the limit set of the group $\rho_g(\grf)$ is identified to the boundary at infinity of $\tilde{\Sigma}$, hence to   $\bgrf_p$ by definition.
\item By Proposition \ref{limset} above,
this definition is independent of the choice of an admissible metric. Moreover, by the same proposition, $\bgrf_p$ is equipped with a  Hölder structure
(coming from the various choices of embedding it as a limit set) so that it makes sense to speak of Hölder functions on $\bgrf_p$ as well as of sets of zero Hausdorff dimension. Finally, $\grf$ acts on $\bgrf_p$ by Hölder homeomorphisms.
\item If there are no cusps and $\Sigma$ is open, $\bgrf_0$ is precisely  the boundary at infinity of the free group $\grf$:
the Cantor set of ends of any Cayley graph of $\grf$.
\item The set $\bgrf_p$ has  a {\em cyclic ordering on points}  depending on the orientation of the surface $\Sigma$. 
By construction of $\bgrf_p$,  an admissible metric  gives rise  to an inclusion
$i $ of $\bgrf_p$ into the boundary at infinity of the hyperbolic plane which is homeomorphic to the circle $\mathbb T$, 
and so it  inherits a  cyclic ordering from the circle. The orientation of $\Sigma$ being fixed,  
this inclusion  $i$ is not canonical but only defined up to a orientation preserving homeomorphism of the circle, 
therefore the cyclic order is a {\em topological} invariant. If we choose the reverse orientation on $\Sigma$, the cyclic ordering is reversed.
This  cyclic ordering  (see Section \ref{simppairs}) 
describes which conjugacy classes of $\grf$ represent simple curves on the surface  $\Sigma$.

\item Even though we shall not use this construction in the sequel, in order to complete this presentation, we sketch a way to obtain the topological sapce $\bgrf_p$ from $\bgrf_0$. 
Choose an admissible metric on $\Sigma$ with $n+p$ boundary components and no cusps. 
Each geodesic boundary component lifts to a union of disjoint  geodesics in ${\mathbb H}$. 
 We consider the equivalence relation $\mathcal R$ on $\bgrf_0$ which identifies the two end points of
  each of the geodesic that comes from points we wish to declare as cusps. 
  Then the quotient  space  $\bgrf/{\mathcal R}$ is homeomorphic to $\bgrf_p$. At this stage however, the Hölder structure is missing.
\end{itemize}

 Nontrivial elements of $\grf$ are of two types according to the number of their fixed points on $\bgrf_{p}$.
\begin{definition}{\sc[Parabolic and hyperbolic elements]}
 The  nontrivial element  $\gamma$ of $\grf$ is {\em hyperbolic}
 if a loop representing $\gamma$ is not freely homotopic to a curve  in a neighbourhood of a cusp. 
 Such a $\gamma$ has precisely two fixed points on $\bgrf_p$: one  attractive $\gamma^+$ and one repulsive $\gamma^-$. 
 
 The element $\gamma$ is {\em parabolic},  if it is nontrivial and if the loop representing $\gamma$ is freely homotopic  to a curve  in a neighbourhood of a cusp. Then
$\gamma$ has precisely one fixed point on $\bgrf_p$ called the  {\em cusp} of $\gamma$. 
 \end{definition}
 
In what follows,  we  adopt the convention that both  $\gamma^+$ and $\gamma^-$ denotes the cusp of a parabolic element of $\gamma$. Every non trivial element of $\grf$ is either parabolic or hyperbolic.

\subsection{Cross ratio for metric spaces}
\subsubsection{Cross ratios}
Let $S$ be a metric space equipped with the action of a group $\Gamma$ by Hölder homeomorphisms.  Let  
$$S^{4*}=\{(x,y,z,t)\in S^4\mid x\not=t , \hbox{ and } y\not= z\}.$$
We equip $S^{4*}$ with the diagonal action of $\Gamma$.

Our main example will be  the boundary at infinity  $\bgrf_p$ of  a surface with $p$ cusps  equipped with the natural  action of $\grf$ by  Hölder homeomorphisms.

\begin{definition}{\sc[Cross ratio]} A {\em cross ratio} on $S$ is a $\Gamma$-invariant Hölder function $\bb$ on $S^{4*}$ with real values  which satisfies the following rules
\label{birrules}\begin{eqnarray}
\bb(x,y,z,t)&=&0\ \ \Leftrightarrow x=y \,\,\hbox{ or } z=t,\\
\bb(x,y,z,t)&=&1 \ \Leftrightarrow x=z\,\, \hbox{ or } y=t\label{bir12},\\
\bb(x,y,z,t)&=&\bb(x,y,w,t)\bb(w,y,z,t)\label{bir11bis},\\
\bb(x,y,z,t)&=&{\bb(x,y,z,w)}{\bb(x,w,z,t)}\label{bir11}.
\end{eqnarray}
The  {\em first cocycle identity  (\ref{bir11bis})} is  a multiplicative  cocycle identity on the first and third arguments, and   the {\em  second cocycle identity (\ref{bir11})} is  a  multiplicative cocycle identity on the second and fourth arguments.
\end{definition}
These rules imply the following symmetries
\begin{eqnarray}
\bb(x,y,z,t)&=&\bb(z,t,x,y), \label{bir100}\\
\bb(x,y,z,t)&=&\bb(z,y,x,t)^{-1},\label{bir101}\\ 
\bb(x,y,z,t)&=&\bb(x,t,y,z)^{-1}\label{bir102}. 
\end{eqnarray}
The classical cross ratio ${\bf b}$  on $\pr$ defined in an affine chart by
$$
{\bf b}(x,y,z,t)=\frac{(x-y)(z-t)}{(x-t)(z-y)},
$$
is an example of a cross ratio with respect to the action of ${\rm PSL}(2,\mathbb R)$.

We will recall in Section \ref{cr-rep} results of \cite{Labourie:2005} which associate cross ratios to representations in ${\rm PSL}(n,\mathbb R)$ as well as other related constructions. In Section \ref{def:hyp}, we explain how a hyperbolic metric on $\Sigma$ gives rise to a cross ratio.

\begin{definition}{\sc[Period]}\label{def:period}
Let $\bb$ be a cross ratio  on $\bgrf_p$ and $\gamma$ be a hyperbolic element in $\grf$.   The {\em period} of $\gamma$ with respect to $\bb$ is \begin{eqnarray}
\ell_\bb(\gamma):=\log\vert \bb(\gamma^-,\gamma y,\gamma^+,y)\vert,\label{period}
\end{eqnarray}
where $\gamma^{+}$ and $\gamma^{-}$ are  respectively  the attracting and  repelling fixed points of $\gamma$ on $\partial_\infty\grf$ and $y$ is any element of $\partial_\infty \grf$ different from $\gamma^+$ and $\gamma^-$.
\end{definition}
Relation (\ref{bir11}) and the invariance under the action of $\gamma$ imply that $\ell_\bb(\gamma)$ does not depend on $y$. Moreover, by Equation (\ref{bir100}), $\ell_\bb(\gamma)=\ell_\bb(\gamma^{-1})$.

\noindent\rmk

The definition given above {\bf does not coincide} with the  definition given for instance in \cite{Hamenstadt:1997,Ledrappier:1995,Otal:1992}
  (even after taking an exponential): indeed, 
we do not require $\bb(x,y,z,t)=\bb(y,x,t,z)$.
However we may observe that if $\bb(x,y,z,t)$ is a cross ratio with our definition, 
so is $\bb^*(x,y,z,t)=\bb(y,x,t,z)$, and finally $\log\vert\bb\bb^*\vert$ is a cross ratio according to the definitions quoted above.

 \subsubsection{Cross ratios and cyclic orderings}
  We observed in the previous paragraph that $\partial_\infty\grf_p$ has a natural cyclic ordering coming from the orientation of the surface. We impose an extra condition on the cross ratio defined on cyclically ordered sets.
  \begin{definition}{\sc[Cross ratio for cyclically ordered sets]}
  Let $S$ be a metric space equipped with a cyclic ordering and an action of a group $\Gamma$ preserving this ordering.
An {\em ordered cross ratio for $S$} is  a cross ratio for the action of  $\Gamma$ on $S$ such that furthermore if  the lexicographical order on the quadruple  $(t,x,y,z)$  coincides with the induced ordering, then both following inequalities are satisfied
  \begin{eqnarray} 
      \bb(x,z,t,y)> 1,\label{ineq1}\\
      \bb(x,y,z,t) < 0. \label{ineq1111}
   \end{eqnarray}
   \end{definition}

\rmks \begin{itemize}
\item This  extra condition is a mild requirement satisfied by all the cross ratios that  are constructed in this article. It does not depend on the choice of the orientation on the surface. 
\item If $\Sigma$ is a closed surface and $\bb$ is a cross ratio on $\bgrf$, then, by continuity, either $\bb$ or $\bb^{-1}$ is ordered. \item In the sequel, when we speak of a cross ratio on $\bgrf_p$, we shall always assume that the cross ratio is ordered and hence implicitely assume Inequalities (\ref{ineq1111}) and (\ref{ineq1}).
\item Finally, if $\gamma$ is a hyperbolic element of $\grf$ and $\bb$ be a (ordered) cross ratio, then $\bb(\gamma^{-},\gamma(x),\gamma^{+},x)>1$. Hence, the absolute values in the definition of the period could be removed and the periods are positive. 
\end{itemize}

 \subsubsection{Cross ratio  and subsurfaces}\label{subsurface}
 
Any  embedding of  $\Sigma(g,n,p)$ in $\check\Sigma(\check g,\check n,\check p)$ injective at the level of fundamental groups,  sending cusps to cusps  and no boundary component to cusp, can be realised by an isometric embedding for some choice of admissible metrics. Hence, such an embedding  induces an Hölder embedding of the corresponding  boundary at infinity. In particular, a cross ratio on $\partial_\infty \pi_1(\check\Sigma)_{\check p}$ induces a cross ratio on $\bgrf_p$. This will be mostly used in the sequel when $\check n=\check p=0$.
  
\subsubsection{Hyperbolic cross ratios}\label{def:hyp}

Assume that $\Sigma$ is equipped with an admissible hyperbolic metric. In particular, deck transformations generate a representation $\rho$  from $\grf$ to $\sld$. Then, $\bgrf_p$ is identified equivariantly with  the limit set of 
$\rho(\grf)$ which is a subset of the boundary at infinity of $\HH$.  In other words, an admissible metric gives rise to a  $\rho$-equivariant map $\iota$ from $\bgrf_p$ to $\partial_\infty\HH\simeq\pr$. Let
${\bf b}$ be  the classical cross ratio. Then,  
$$
\bb_{\rho}(x,y,z,t)={\bf b}(\iota(x),\iota(y),\iota(z),\iota(t)),
$$ is a cross ratio on $\bgrf_p$. 

\begin{definition}{\sc[hyperbolic cross ratio]} A cross ratio obtained through the previous construction is a {\em hyperbolic cross ratio}. 
\end{definition}

For hyperbolic cross ratios, periods coincide with lengths of closed geodesics. 

\section{Gaps and gap functions}\label{fond}

In this section, we discuss the geometry of the space of pairs of distinct points of  the boundary at infinity. We first introduce a terminology to describe intersection of geodesics in terms of their end points. We then present the Birman--Series set and state their remarkable theorem \cite{Birman:1985}. Finally, we describe {\em gaps} and recall that they are related to embedded pairs of pants.

\subsection{The boundary at infinity and embedding of geodesics}
If $\Sigma$ is equipped with an admissible metric, a geodesic in its universal cover gives rise to a pair of points in the boundary at infinity.
We introduce   a terminology based purely on properties of configurations of  points in  $\bgrf_p$ relative to its cyclic ordering, then translate this terminology into properties of the  associated geodesics (compare \cite{McShane:1998,Mirzakhani:2007a})).

\begin{definition}{\sc[Subsets of $\bgrf_p$]}\begin{itemize}
\item 
Let $a$ and $b$ be  distinct points in $\bgrf_p$. An {\em open interval} $]a,b[$ in $\bgrf_p$ is the set of those elements $t$ in $\bgrf_p$ different from $a$ and $b$ and such that  $(a,t,b)$ is positively oriented. The {\em closed interval} $[a,b]$ is  $]a,b[\cup\{a,b\}$.
\item Let $X,Y$ be subsets of $\bgrf_{p}$. We say that  $X$ {\em does not separate } $Y$ if there exists disjoint closed intervals $I$ and $J$ such that $I$ contains $X$ and $J$ contains $Y$. 
\item The subset $X$  {\em injects} if  for all elements $\gamma$ in $\grf$, 
$X$ does not separate $\gamma(X)$.

\item 
The pair $X=\{x,y\}$  is {\em simple} if it injects and \label{simppairs}
is fixed by a nontrivial element  of $\grf$.\end{itemize}
\end{definition}

Observe that intervals may not be connected.

If $x_1$, $y_1$, $x_2$ and $y_2$ are pairwise distinct points of $\bgrf_p$
and  $\bb$  is a cross ratio -- which is now by convention ordered -- on $\bgrf_{p}$, observe that
$\{x_1,x_2\}$ does not separate $\{y_1,y_2\}$
if and only if
$$\bb(x_{1},y_{1},x_{2},y_2)> 0.$$
Indeed, by Properties (\ref{bir101}) and  (\ref{bir102})  this statement is independent under renumbering of $\{x_1,x_2\}$ and $\{y_1,y_2\}$. But, after renumbering,  either $(x_1,x_2,y_1,y_2)$ is oriented in which case $\{x_1,x_2\}$ does not separate $\{y_1,y_2\}$ and $\bb(x_{1},y_{1},x_{2},y_2)> 0$, or $(x_1,y_1,x_2,y_2)$ is oriented  in which case $\{x_1,x_2\}$ separates $\{y_1,y_2\}$ and $\bb(x_{1},y_{1},x_{2},y_2)< 0$.

If $X$ is a simple pair, there exists an element $\alpha$ in $\grf$, well defined up to multiplicity in $\mathbb{Z}^{*}$, so that 
$X=\{\alpha^+,\alpha^-\}$.

\begin{definition}{\sc[Peripheral elements]}  A non trivial  element $\alpha$ in $\grf$ is 
{\em peripheral}  if it is not parabolic and if the  pair  $\{\alpha^+,\alpha^-\}$ does not separate $\bgrf_{p}$.
By extension, we  also say that $\alpha^+$ and $\alpha^-$ are {\em peripheral} points.
\end{definition}

Note that if $\alpha$ is peripheral then $\{\alpha^+,\alpha^-\}$ is  simple.

The following dictionary relates the previous definitions to intersection properties of geodesics.
Let $X=\{x_1,x_2\}$ be a pair. If we equip $\Sigma$ with an admissible metric and denote by $\widetilde{\gamma}_{X}$  the geodesic joining $x_1$ to $x_2$ in the universal cover of $\Sigma$ and by ${\gamma_{X}}$ its projection on $\Sigma$, then we have the dictionary
\begin{itemize}
\item Let $Y$ be a pair. Then, "{\em $X$ does not separate $Y$}" is equivalent to ``{\em $\widetilde{\gamma}_X$ and $\widetilde{\gamma}_Y$ do not intersect}".
\item ``{\em $X$ is simple}" is equivalent to  ``{\em $\gamma_X$ is a  simple closed geodesic}".
\item ``{\em $X$ injects}" is equivalent to ``{\em $\gamma_X$ has no transverse self intersection}", in particular the closure of $\gamma_X$ is a geodesic lamination.
\end{itemize}

\subsubsection{Pairs of pants}

Let $P $ be a  sphere minus three disjoint open discs, together with a base point and a choice of orientation.  
We choose once and for all a {\em marking} of $P $, that is three elements $\alpha_0$,$\beta_0$, and $\gamma_0$ of $\pi_1(P )$
represented by three  loops each of which is freely homotopic to one of the  boundary loops (with the orientation induced from the sphere)
and such that $\alpha_0\gamma_0\beta_0=1$. 

We now define what we mean by an  embedded  pairs of pants
or more properly an  isotopy  class of embedding of a pair of pants.

\begin{definition}{\sc [Embedded  pairs of pants]} 
An {\em isotopy  class of  pairs of pants} in an orientable  surface $\Sigma$, possibly with boundary,
 is an isotopy  class of orientation preserving incompressible  embeddings  
of $P $ in $\Sigma$, i.e.  the induced morphism on the fundamental groups is an injection. 
We write $\alpha, \beta,\gamma$  respectively
for the images in  $\pi_{1}(\Sigma)$ of  $\alpha_0,\beta_0,\gamma_0$;
 the triple $(\alpha,\beta,\gamma)$ is well defined up to conjugacy by  $\grf$ and satisfies 
  $\alpha\gamma\beta=1$. 
For any such a triple we shall call $\alpha$  the {\em first} boundary component of the  pair of pants
and, for brevity, we shall say that  the triple $(\alpha,\beta,\gamma)$ {\em represents} the pair of pants.
   \end{definition}

\rmk \label{technical}
Note that, with our definition above,
 if  $(\alpha, \beta,\gamma)$ 
 represents an embedding of $P$ then 
$(\alpha, \alpha\beta\alpha^{-1},\alpha\gamma\alpha^{-1}) $
still  represents the same isotopy class of embedding of $P$. 
Though this is a departure  from the usual definition for surfaces with boundary
(which we give below and refer to as pointed pairs of pants) 
it will make the statement of our identities as well as their proofs slightly more ``compact".
Informally,  an embedded pair  of pants is the orbit under the Dehn twist round the boundary
loop $\alpha$ of an embedded pointed pair  of pants.
This distinction is important particularly in the proof of Theorem \ref{fondbound}.

\begin{definition}{\sc [Embedded  pointed pair of pants]}
Let   $\Sigma$ be a surface with boundary.
Fix a basepoint $p$ for the pair of pants $P$ on the boundary loop representing  $\alpha_{0}$
and   a fixed point $p'$ for $\Sigma$ on a  distinguished boundary component $\alpha$.
We define an isotopy  class of  pointed pairs of pants in  $\Sigma$ 
to be an isotopy  class of orientation preserving embeddings  as above,
 but with the additional hypothesis that, for \textit{each} such  embedding and each such isotopy,  
 the image of $p$ is $p'$.
With this definition if   $(\alpha, \beta,\gamma)$ 
 represents an embedding of $P$ then 
$(\alpha, \alpha\beta\alpha^{-1},\alpha\gamma\alpha^{-1}) $
represents a distinct ( from  $(\alpha, \beta,\gamma)$) isotopy class of embedding of $P$. 
  \end{definition}

\begin{figure}[htbp] 
   \centering
   \includegraphics[width=2in]{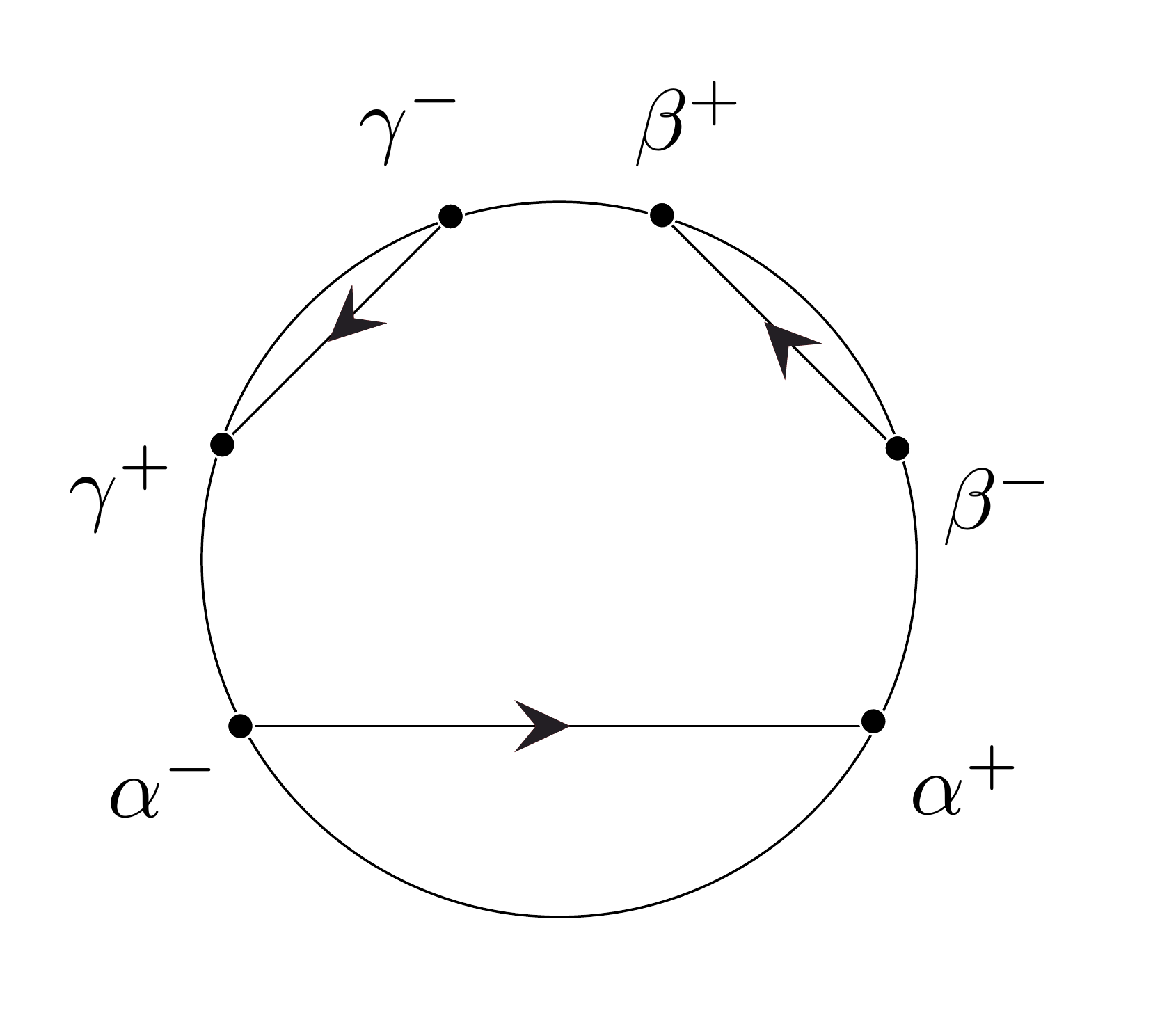}  
   \caption{Associated sextuplet to a pair of pants}
   \label{fig:sextuple}
\end{figure}

Observe then that the sextuplet $(\alpha^-,\alpha^+,\beta^-,\beta^+,\gamma^-,\gamma^+)$ is positively oriented as in Figure \ref{fig:sextuple}.

\subsection{Birman--Series sets}

\subsubsection{The Birman--Series Theorem}

Let 
$$
\bgrf_p^{2*}=\{(x,y)\in\bgrf_p^2\mid x\not= y\}.
$$
\begin{definition}
The {\em Birman--Series set} is the set of  couples in $\bgrf_p^{2*}$ whose associated pairs inject. \end{definition}
The Birman--Series set is a $\grf$-invariant set and by the previous paragraph a couple
$X=(x,y)$ belongs to the Birman--Series set if and only if the projection to the surface $\Sigma$ of the geodesic joining $x$ to $y$ for an admissible metric 
is a complete simple geodesic.
Joan  Birman and Caroline Series  studied  in \cite{Birman:1985} the set of all complete simple geodesic on a 
hyperbolic surface  and  showed that the points of this set are somewhat sparse.  

\begin{theorem}\label{Birman:1985}{\sc[J.Birman--C.Series]}
The union of the simple complete geodesics  of a compact
hyperbolic surface  is closed and has Hausdorff dimension $1$. In particular the Birman--Series set has zero Hausdorff dimension.
\end{theorem}

\subsubsection{Relative Birman--Series sets}
Let $\alpha$ be a  peripheral or parabolic element of $\grf$. 

\begin{definition}{\sc[Relative Birman--Series set]} The  {\em first relative Birman--Series set with respect to $\alpha$} is the subset
 $K_\alpha$ of $\bgrf_p\setminus \{\alpha^+,\alpha^{-}\}$ of those elements $t$ such that $\{\alpha^+,t\}$ injects.

The  {\em second relative Birman--Series set with respect to $\alpha$} is the subset $K_\alpha^*$ of $K_\alpha $ of those elements which do not belong to the $\grf$ orbit of $\alpha^+$.
\end{definition}
In our dictionary,
  $K_{\alpha}$ corresponds to the set of (lifts) of geodesics 
that are simple and spiral to the boundary component $\alpha \subset \Sigma$ in a given direction if $\alpha$ is hyperbolic, 
or go to the cusp defined by $\alpha$  if $\alpha$ is parabolic. These relative sets are subsets of the Birman--Series set, and hence   have zero Hausdorff dimension. 

\begin{proposition}\label{K props}
The sets $K_{\alpha}$ and $K_\alpha^*$ are closed (for the relative topology) subsets  of $]\alpha^+,\alpha^{-}[$.
If $x$ is a fixed point of a peripheral or parabolic element of $\grf$ in $K_\alpha$ then $x$ 
is isolated in $K_{\alpha}$. 
\end{proposition}

\proof  Note first, that the  Birman-Series set is closed so $K_\alpha \cup \{\alpha^+,\alpha^{-}\}$
is immediately seen to be closed too. 

The fact that $K_\alpha^*$ is a closed subset will follow from the second statement, 
since every element in  the orbit of $\alpha^+$ is peripheral or parabolic and hence isolated. 
One shows this second statement as follows. Let $x= \beta^+$ be the attractive  fixed point of a peripheral hyperbolic element $\beta$. We assume that $(\alpha^+,\beta^+,\beta^-,\alpha^-)$ is positively oriented, the other case being treated similarly.  
If $y$ belongs to  $]\beta(\alpha^{+}),\beta^{+}[$, then  $\{\alpha^+,y\}$ does not inject since
$$(\alpha^-,\alpha^+,\beta(\alpha^{+}),y,\beta(y), \beta^{+}),$$ is positively oriented and,  
hence $\{\alpha^+,y\}$ separates $\{\beta(\alpha^+),\beta(y)\}$. Hence the intersection of  $K_\alpha$  with the 
 interval  $]\beta(\alpha^{+}),\beta-[$ contains only $\beta^+$. 
 When $\beta$ is a parabolic, we restrict ouselves to the case $(\alpha^+,\beta^+,\alpha^{-})$ is oriented and  we show using the same argument that $]\beta(\alpha^+),\beta^{-1}(\alpha^-)[$ only contains $\beta^+$. 
\qed

\subsection{Gaps and pair of pants}

Let $\alpha$ be a peripheral or parabolic element.
\begin{definition}{\sc[Gap]}
 A {\em gap} with respect to $\alpha$ is a  couple  $(x,y)$ of elements of 
  $K^*_\alpha$ so that the interval $]x,y[$ does not intersect $K^*_\alpha$.
 \end{definition}

Note that $\bgrf_p\setminus K^*_\alpha$ is a disjoint union of intervals associated to gaps. In the sequel, we shall abusively refer to the interval $]\alpha,\beta[$ when $(\alpha,\beta)$ is a gap also as a gap.

\begin{figure}[htbp] 
   \centering
   \includegraphics[width=3in]{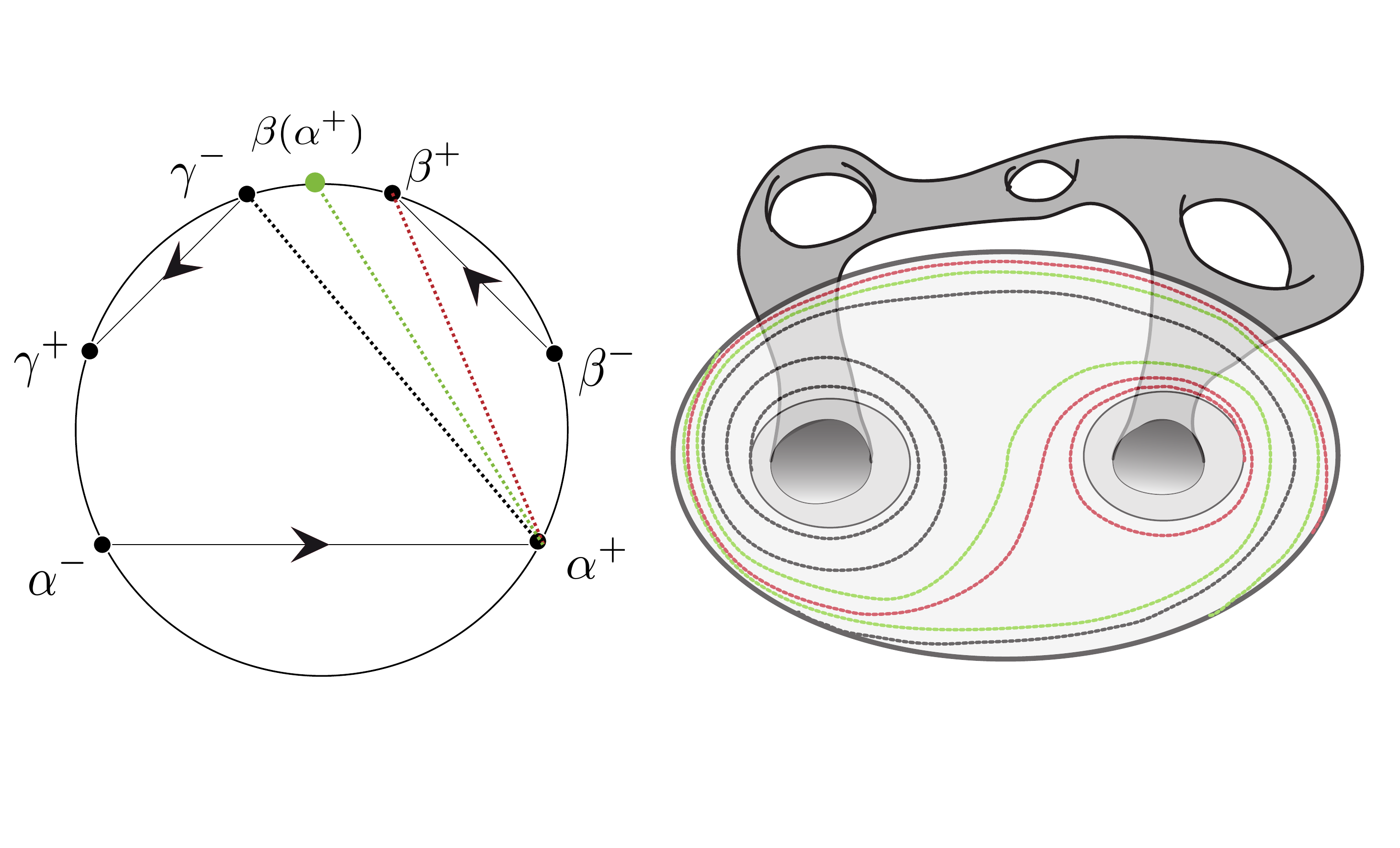}  
   \caption{Gap associated to  a pair of pants}
   \label{fig:2fig}
   \end{figure}

There are exactly two types of gap, those arising from pair of pants, those arising from peripheral elements. The following result is a rephrasing of the corresponding statements in \cite{McShane:1998} and \cite{Mirzakhani:2007a}.
\begin{proposition}\label{pants}{\sc[Gaps and pants]}
If $(\alpha,\beta,\gamma)$ represents a pair of pants then:
\begin{enumerate}
\item  The pair  $(\beta^+,\gamma^-)$ is a gap.
\item  If  furthermore  $\beta$ is peripheral, then  $(\beta^-,\beta^+)$ is a gap.
\end{enumerate}
Conversely, let $(x,y)$ be a gap.
Then, there exist  primitive  elements $\beta$ and $\gamma$  in  $\grf$  such that 
$(x,y) = (\beta^+,\gamma^-)$. Moreover:
\begin{enumerate}
\item If $\beta\not=\gamma$  then, after possibly taking inverses whenever $\gamma$ or $\beta$ are parabolic,  $(\alpha,\beta,\gamma)$ represents a possibly degenerate pair of pants.  Such  group elements  $\beta$ and $\gamma$ are then unique. 
\item If $\beta=\gamma$, then $\beta$ is peripheral and there exists a unique embedded (pointed) pair of pants  $P$ so that 
$\alpha$ and $\beta$ are boundary loops of $P$.
\end{enumerate}

\end{proposition}
 Figure \ref{fig:2fig} represents a gap arising form a pair of pants represented by $(\alpha,\beta,\gamma)$: the left figure represents the gap $(\beta^+,\gamma^-)$ in the universal cover, with the point $\beta(\alpha^+)$ which is the unique point of $K_\alpha^*$ that belongs to the gap and the lift of the geodesic laminations in dotted lines, the right figure  represents the pant and the geodesic laminations which are the projections in the compact surface of the same objects in the previous picture.

\proof  We choose an admissible metric on $\Sigma$. Let $\pi$ be the orthonormal projection from $\bgrf_p$ to the boundary component represented by $\alpha$ if $\alpha$ is peripheral, or to an horosphere centred at the corresponding cusp if $\alpha$ is parabolic. Then  one checks that $x$ belongs to $K_\alpha$ if and only if the orthogonal geodesic -- to either the boundary component or the horosphere --starting from $\pi(x)$ has no self intersection.

Hence, the result follows from Theorem 4.6 and the discussion afterwards  in \cite{Mirzakhani:2007a} when $\alpha$ is peripheral and the corresponding discussion in \cite{McShane:1998} when $\alpha$ is parabolic.\qed

\section{Identities}\label{fondid}

We prove in this section our main results Theorems \ref{fondbound} and  \ref{cuspcase} which state an identity at a boundary component and at a cusp respectively. To treat the cusp case, we shall have to impose an extra regularity assumption on the cross ratio. These identities involve {\em gap functions} which we describe in preliminary paragraphs.

\subsection{The identity at a boundary component}
\subsubsection{Gap functions}
Let  $\bb$ be an ordered  cross ratio on $\bgrf_p$. 
Let $\mathcal P_\alpha$ denote the set of (embeddings of) pairs of pants (possibly degenerate) up to
isotopy which have $\alpha$ as  first boundary component.
Let $\mathcal S_\alpha$ denote the set of pairs of pants up to
isotopy which have $\alpha$ as first boundary component as well as another boundary component of $\Sigma$.

Observe that since we allow isotopies that do not fix $\alpha$ pointwise, the two triples $(\alpha,\beta_0,\gamma_0)$ and  $(\alpha,\beta_1,\gamma_1)$ of elements of $\bgrf^3$
 represent the same class in $\mathcal P_\alpha$ if and only if there exists an integer $n$ such that
\begin{eqnarray*}
\beta_0=\alpha^{n}\beta_1\alpha^{-n},\\ 
\gamma_0=\alpha^{n}\gamma_1\alpha^{-n}. 
\end{eqnarray*}
 
\begin{definition}{\sc[Gap functions]}\begin{itemize}
\item Let $P$ be  a pair of pants  in  $\mathcal P_\alpha$  represented by the triple $(\alpha,\beta,\gamma)$. 
The value of the {\em pant gap function} $G_\bb$ at $P$  is 
$$
G_\bb(P)=\log(\bb(\alpha^+,\gamma^-,\alpha^-,\beta^+)).
$$
\item
Let $P$ be a pair of pants in $\mathcal S_\alpha$   represented by the triple $(\alpha,\gamma_0,\gamma_1)$ so that  $\beta=\gamma_i$  is another peripheral element for $i=0$ or $i=1$. 
The value of the  {\em boundary gap function} $G^r_\bb$  at $P$ is
$$
G^r_\bb(P)=\log (\bb(\alpha^+,\beta^+,\alpha^-,\beta^-)).
$$
\end{itemize}
\end{definition}
Note in the second case that if $\beta$ is parabolic then $\beta^+= \beta^-$ and $G^r_\bb(P)= 0$. Finally,  the gap functions are well defined and have positive values: indeed, according to Figure \ref{fig:sextuple}, $(\alpha^+,\beta^-,\beta^+,\alpha^-)$ and  $(\alpha^+,\beta^+,\gamma^-,\alpha^-)$  are  positively oriented hence by  Inequality (\ref{ineq1}), $\bb(\alpha^+,\gamma^-,\alpha^-,\beta^+)>1$ and $\bb(\alpha^+,\beta^+,\alpha^-,\beta^-)>1$.

\subsubsection{The main result for a boundary component}

We can now state  the  identity. We use the terminology of the beginning of  Section \ref{cross ratiodef}.

\begin{theorem}\label{fondbound}
Let $\Sigma$ be a compact connected oriented surface of finite type with
boundary components and $p$ cusps. Let $\alpha$ be a primitive peripheral element in $\grf$ which represents an oriented boundary component. Let $\bb$ be an ordered  cross ratio on $\bgrf_{p}$.
Then
$$
\ell_\bb(\alpha)=\sum_{P\in\mathcal P_\alpha}G_\bb(P) +\sum_{P\in
\mathcal S_\alpha}G^r_\bb(P).
$$
\end{theorem}

\rmks \begin{itemize}
\item Our identities in the Fuchsian case are equivalent to those of Mirzakhani 
although the way we count contributions from pants is slightly different.
Mirzakhani separates the set of embeddings of pants $\mathcal P_\alpha$ into 
those containing two boundary components  of $\Sigma$ ($\mathcal S_\alpha$ in our notation) and 
those that have a single boundary component ($\mathcal{ P}_\alpha \setminus \mathcal S_\alpha$ in our notation).
On the other hand, our first series counts a contribution from every embedding whether it
has one or two boundary components and the second is a ``correction term" which 
corresponds to the contribution  due to geodesics that ``escape the surface" 
via  the  second boundary component of an embedding in $\mathcal S_\alpha$. This is only a matter of convention and taste.
\item The reader may wonder why the  pair of pants whose boundary components are cusps do not appear $\mathcal S_\alpha$. The answer is that adding them to $\mathcal S_\alpha$ has no effect: these pair of pants  have zero relative gaps functions. 
\end{itemize}

The strategy of the proof, as in \cite{McShane:1998}, 
 is to compute the length of a circle  (the quotient of $]\alpha^{+},\alpha^{-}[$ by $\langle \alpha \rangle$)
 in terms of the lengths of the complementary regions,  
 which are gaps,  of  the set  $K_{\alpha}^*$.

The identity holds for more general functions than cross ratios. Indeed, we shall not use  the first  cocycle identity on the first and third arguments, and only use the second cocycle identity on the second and fourth arguments. 

\proof Let $\bb$ be a cross ratio and $\alpha$ be a primitive peripheral element of $\grf$ as in the statement.
We write $]\alpha^+,\alpha^-[$ for $\bgrf_p\setminus\{\alpha^-,\alpha^+\}$
and we fix a reference point $\zeta\in ]\alpha^+,\alpha^-[$.
 Observe that for any point  $y$ in  $]\alpha^+,\alpha^-[$,
 $\bb(\alpha^+,y,\alpha^-,\zeta)$ is positive by Inequality (\ref{ineq1}) and Equality (\ref{bir102}).
 
Let  $B$ be the map  from $]\alpha^+,\alpha^-[$ to $\mathbb R$ defined by
$B(y)=\log(\bb(\alpha^+,y,\alpha^-,\zeta))$;
Note that, since $\bb$ is Hölder,  $B$ is Hölder.
Note further that  $B(K_\alpha^*)$ has zero Hausdorff dimension, 
since by Theorem  \ref{Birman:1985}, $K_\alpha^*$ has zero Hausdorff dimension.
Moreover, $B$ is injective 
since  $\bb(\alpha^+,y,\alpha^-,\zeta) = \bb(\alpha^+,y',\alpha^-,\zeta))$
implies $\bb(\alpha^+,y,\alpha^-,y') = 1$ by (\ref{bir11bis})
and so $y=y'$ by  (\ref{bir12}). 
Observe moreover that $B$ preserves the ordering  by Inequality (\ref{ineq1}) -- hence is a homeomorphism -- so  we identify $]\alpha^+,\alpha^-[$ with its image under $B$.

Let $\mu$ be the Lebesgue measure on $\mathbb R$. Recall that a set of Hausdorff zero dimension has zero Lebesgue measure. 
Hence  $
\mu(K_\alpha^*)=0.
$

By the second cocycle identity
$$
B(\alpha(z))=B(z)-\ell_\bb(\alpha). 
$$
Let  $\mathbb T=\mathbb R/\ell_\bb(\alpha)\mathbb Z$, and $\pi$ be the projection from $\mathbb R$ to $\mathbb T$.
The set $K_\alpha^*$ is also invariant by $\alpha$ and we continue to  denote by $K_\alpha^*$ its projection on $\mathbb T$.
 Let $\widehat I_P=[\beta^+,\gamma^-]\subset \mathbb R$ for each pair of pants $P=(\alpha,\beta,\gamma)$ 
and  $\widehat J_P=[\beta^-,\beta^+]$
for each pair of pants  $P$ represented by $(\alpha,\beta,\gamma)$
where $\beta$ is peripheral  (see Proposition \ref{pants} and Figure \ref{fig:sextuple}). 

Observe that by Condition (\ref{bir12}) and since  $\alpha\widehat J_P$ and $\widehat J_P$ are disjoint,  $\pi$ is injective from $\widehat J_P$ to $J_P=\pi \widehat J_P$,  similarly from $\widehat I_P$ to $I_P=\pi \widehat I_P$. 
Moreover, the set of intervals $I_P$ is in bijection with  ${\mathcal P}_\alpha$ 
and, since we prefer $\beta$ over $\gamma$, there is a two to one map from  ${\mathcal S}_\alpha$ to the set of intervals $J_P$.  
By Proposition \ref{pants} 
$$
\mathbb T\setminus K_\alpha^*=(\sqcup_{P\in{\mathcal P}_\alpha}I_P)  \sqcup (\sqcup_{P\in \mathcal S_\alpha }J_P).
$$
By construction and the second cocycle identity (\ref{bir11}.) 
\begin{itemize}
\item 
 $\mu(I_P)=G_{b}(P)$, if $P$ belongs to ${\mathcal P}_\alpha$,
 \item
  $\mu(J_P)=G_{b}^r(P)$, if $P$ belongs to $S_\alpha$.
\end{itemize}
The identity  then  follows from 
\begin{eqnarray*}
\ell_\bb(\alpha)&=&\mu(\mathbb T)=\mu(\mathbb T\setminus K_\alpha^*)\\
&=&\sum_{P\in{\mathcal P}_\alpha}\mu(I_P)+\sum_{P\in{\mathcal S}_\alpha}\mu(J_P)
=\sum_{P\in{\mathcal P}_\alpha}G_\bb(P)+\sum_{P\in{\mathcal S}_\alpha}G_b^r(P).
\end{eqnarray*}\qed

\subsection{The identity at a cusp}\label{defcusp}

In order to state and prove the identity at a cusp, we make in this section the following further  assumption on the cross ratio $\bb$. 
\begin{definition}{\sc{[Regularity Hypothesis]}}\label{reghyp}

There is an orientation preserving Hölder embedding of $\bgrf_{p}$ in the circle  $\mathbb T$, such that the action of $\grf$ turns out to be $C^1$  and that moreover
for every pair of distinct points $s,t$ in $\partial_{\infty}\pi_{1}(\Sigma)_{p}$
 the function $$
(x,y)\mapsto \bb(x,s,y,t),
$$
is the restriction of a  $C^1$ function with non zero Hölder derivatives along the diagonal in $\mathbb T^2$. 
\end{definition}

\subsubsection{Auxiliary functions at a cusp}\label{walpha}

Let $\alpha$ be a correctly oriented  primitive parabolic element of $\grf$ with fixed point $\alpha^+$. 
By \textit{correctly oriented} we mean that  
the orientation of a loop in the surface representing $\alpha$ coincides  with the  orientation around the boundary component. 
We define for any  points $s$, $t$, $s_0$ in $\bgrf_p$
$$
W_\alpha(s,t)=\left.\frac{\partial_y \log \bb(\alpha^+,s,y,t)}{\partial_y \log \bb(\alpha^+,s_0,y,\alpha(s_0))}\right\vert_{y=\alpha^+}.
$$ By the hypothesis above, $W_\alpha$ is well defined and Hölder.

\begin{proposition}\label{propW}
The function $W_\alpha(s,t)$ does not depend on the choice of $s_0$. Moreover,
\begin{eqnarray}
W_\alpha(\alpha(s),\alpha(t))&=&W(s,t),\label{winvar}\\
W_\alpha(s,\alpha(s))&=&1,\label{wnorm}\\
W_\alpha(s,u)&=&W_\alpha(s,t)+W_\alpha(t,u).\label{addcoc}
\end{eqnarray}
Finally, if $(\alpha^+,s,t)$ is positively oriented, then $W_\alpha(s,t)$ is negative.
\end{proposition}
\proof
We first remark that the differential of $\alpha$ at $\alpha^+$ is the identity. Indeed, maybe after exchanging $\alpha$ and $\alpha^{-1}$,  the triple $(y,\alpha(y),\alpha^+)$ is always oriented, hence the right derivative of $\alpha$ at $\alpha^+$ is no more  than 1, and the left derivative is no less than 1. Therefore 
$$
\left.\partial_y \log \bb(\alpha^+,s,y,t)\right\vert_{y=\alpha^+}=\left.\partial_y \log \bb(\alpha^+,s,\alpha(y),t)\right\vert_{y=\alpha^+}.$$
This, together with the invariance of $\bb$, shows Equality (\ref{winvar}).

We define 
$$
R_x(s,t,s_0,t_0)=\frac{\partial_y \log \bb(x,s,y,t)}{\partial_y \log \bb(x,s_0,y,t_0)}\bigg\vert_{y=x}.$$
Observe that  $R_x(s_0,t_0,s_0,t_0)=1$ and that by the second cocycle identity,
\begin{eqnarray*}
R_x(s,u,s_0,t_0)&=&R_x(s,t,s_0,t_0)+R_x(t,u,s_0,t_0).
\end{eqnarray*}
We claim  that $R_{\alpha^+}(s,t,s_0,\alpha (s_0))$ does not depend on the choice of $s_0$.
Indeed, by the second cocyle identity and the invariance under $\alpha$
\begin{eqnarray*}
& &\frac{1}{R_{\alpha^+}(s,t,s_0,\alpha (s_0))}-\frac{1}{R_{\alpha^+}(s,t,t_0,\alpha (t_0))}\\
&=&\frac{1}{R_{\alpha^+}(s,t,s_0,t_0)}-\frac{1}{R_{\alpha^+}(s,t,\alpha (s_0),\alpha (t_0))}\\
&=&\frac{\partial_y \log \bb(\alpha^+,s_0,y,t_0)}{\partial_y \log \bb(\alpha^+,s,y,t)}\bigg\vert_{y=\alpha^+}-\frac{\partial_y \log \bb(\alpha^+,\alpha(s_0),y,\alpha(t_0))}{\partial_y \log \bb(\alpha^+,s,y,t)}\bigg\vert_{y=\alpha^+}\\&=&0.
\end{eqnarray*}
Since
$
W_\alpha(s,t)=R_{\alpha^+}(s,t,s_0,\alpha (s_0)),
$
this concludes the proof of the first part of the proposition. 
The last statement follows from the previous observation.\qed

\subsubsection{Cusp gap function}\label{cuspgapfunction}
Let again $\mathcal P_\alpha$ denotes the set of pair of pants on $\Sigma$ which have the cusp of $\alpha$ as a boundary component
and  $\mathcal S_\alpha$ the set of pair of pants  which have the cusp of $\alpha$ as a boundary component
as well as some other boundary $\beta$.  
As before, we consider pairs of pants up to isotopies that leave  $\alpha$ invariant.  

\begin{definition} The {\em cusp gap function} of the pair of pants $P$ represented by $(\alpha,\beta,\gamma)$ in $\grf^3$ is
$$
W_{\bb}(P)=W_\alpha(\gamma^-,\beta^+).
$$
For  $P$ in $\mathcal S_\alpha$, the {\em boundary cusp  gap function} is
$$
W_{\bb}^r(P)=W_\alpha(\beta^+,\beta^-).
$$
\end{definition}
Since $W_\alpha$ is invariant under the action of $\alpha$, it follows that $W_\bb(P)$ and  $W_\bb^r(P)$ only depend on the isotopy class of $P$.

\subsubsection{The identity}
Our identity for cusps is
\begin{theorem}\label{cuspcase}
\begin{equation}
1=\sum_{P\in\mathcal P_\alpha}W_{\bb}(P) +   \sum_{P\in\mathcal S_\alpha}W_{\bb}^r(P)
\end{equation}
\end{theorem}
Observe again that pair of pants with another cusp as boundary component have zero relative gap function.

\proof The proof is almost exactly the same as that of Theorem \ref{fondbound}.
Let   $\alpha$ be a correctly oriented  primitive parabolic element of $\grf$.
Since $\alpha$ is parabolic $\alpha^+= \alpha^-$.
 One then defines a map  of $\bgrf_p\setminus\{\alpha^+\}$ in $\mathbb R$  by  $$
B(y)=W_{\alpha}(s,y).
$$
where $s \in \bgrf_p \setminus \{\alpha^+\}$
is an arbitrary  point. Moreover, by the last assertion in Proposition \ref{propW}, if $(\alpha^+,t,y)$ is oriented, then
$$
B(y)-B(t)=W_{\alpha}(t,y)<0.
$$
Thus $B$ is monotone and injective.

Noting that $B(\alpha(y)) - B(y)= W_\alpha(y,\alpha(y)) = 1$,
 we replace $\ell_{\bb}(\alpha)$ by 1 
in the proof of Theorem \ref{fondbound} to obtain the equation.
 \qed

\section{Gaps and coordinates in hyperbolic geometry}\label{hypcase}

In this section, we restrict ourselves to hyperbolic cross ratios and express the gap functions of a pair of pants as a function  of the lengths of the boundary components, and more generally of the {\em shear coordinates} of the pair of pants.  We recover the previous results by G. McShane and M. Mirzakhani  \cite{Mirzakhani:2007a,McShane:1998} and do not claim any originality about these results. Nevertheless, the method is new and involves only the formal properties of hyperbolic  cross ratios  instead of hyperbolic tri\-go\-no\-me\-try as in the original proofs. Moreover, these computations will be helpful later. This approach also emphasises the importance of the notion of cross ratio.

\subsection{Gap function for a pair of pants} 
The main result of this section is the following

\begin{theorem}
Let $\bfb$ be a hyperbolic cross ratio.
Let $P$ be a pair of pants represented by  $(\alpha,\beta,\gamma)$.  
For the boundary gap functions, we assume that $\beta$ is also a  peripheral element.
 Let $\ell(\alpha)$,
 $\ell(\beta)$ and  $\ell(\gamma)$ be the lengths of the corresponding boundary components. 
 If $\alpha$ is a correctly oriented primitive peripheral element then the gap and boundary gap functions are given by 
\begin{eqnarray}
G_\bfb(P)&=&\log\left(\frac{e^{\frac{\ell(\beta)+\ell(\gamma)}{2}}+e^{\frac{\ell(\alpha)}{2}}}{e^{\frac{\ell(\beta)+\ell(\gamma)}{2}}+e^{-\frac{\ell(\alpha)}{2}}}\right),\\
G^r_\bfb(P)
&=&\log\left(\frac{\cosh \left(\frac{\ell(\gamma)}{2}\right)+\cosh\left(\frac{\ell(\beta)-\ell(\alpha)}{2}\right)}{\cosh \left(\frac{\ell(\gamma)}{2}\right)+\cosh\left(\frac{\ell(\beta)+\ell(\alpha)}{2}\right)}\right).
\end{eqnarray}
Moreover, assume $\alpha$ is a correctly oriented primitive peripheral element in $P$. Then, the cusp gap functions are  given by
\begin{eqnarray*}
W_ {\bfb}(P)&=&\frac{1}{1+e^{-\frac{\ell(\beta)+\ell(\gamma)}{2}}},\\
W_\bfb^r(P)&=&\frac{\sinh\left(\frac{\ell(\beta)}{2}\right)}
{\cosh\left(\frac{\ell(\gamma)}{2}\right)+\cosh\left(\frac{\ell(\beta)}{2}\right)}.
\end{eqnarray*}
\end{theorem} 

We split the theorem in Propositions \ref{gapsl2} and \ref{gapsl2pinch}. We prove them 
 using Thurston's shear coordinates \cite{Bonahon:1996,Bonahon:2001,Thurston:1984} that we describe in the next section. We recover this way Formulae (\ref{McShane:1998}) and (\ref{Mirzakhani:2007a}) given in the introduction.

The proofs will make a heavy use of the extra identity satisfied by a hyperbolic cross ratio $\bfb$, namely
\begin{eqnarray}
1-\bfb(f,v,e,u)=\bfb(u,v,e,f)\label{cr0},
\end{eqnarray}
or equivalently
\begin{eqnarray}
\bfb(x,y,z,t)=1-\frac{1}{\bfb(y,z,x,t)}=\frac{1}{1-\bfb(z,x,y,t)}.\label{projcra}
\end{eqnarray}
These rules implies the following symmetry
\begin{eqnarray}
\bfb(x,y,z,t)=\bfb(y,x,t,z).\label{projsym}
\end{eqnarray}

\subsection{Length functions and shear coordinates}\label{shear} 
 Let $P$ be the sphere minus three points   and let $(\alpha,\beta,\gamma)$ be a marking given by  three  elements of $\pi_1(P)$  such that 
$
\alpha\gamma\beta=1
$.
We assume that $P$ is equipped with a finite area hyperbolic metric whose completion has $p$ cusps and $3-p$ totally geodesic boundary components.
Let   $\alpha_0, \beta_0, \gamma_0$ be  fixed points of $\alpha$, $\beta$, $\gamma$ on $\partial_\infty\pi_1(P)_{p}$
 which we consider as a subset of $\mathbb{RP}^1$. We shall denote  $\alpha_1$ the other fixed point of $\alpha$ if $\alpha$ is hyperbolic. If $\alpha$ is parabolic, we set as usual $\alpha_0=\alpha_1$. Observe that $\alpha^{-1}(\beta_0)=\gamma(\beta_0)$. Let finally $\bfb$ be the hyperbolic cross ratio.

 \begin{definition}{\sc[Shear coordinates]} 
 The {\em shear coordinates} of the hyperbolic pair of pants $P$ are the positive numbers $A$, $B$ and $C$ defined by
\begin{eqnarray*}
B&=&-\bfb(\alpha_0,\beta_0,\gamma_0,\alpha^{-1}(\beta_0))=-\bfb(\alpha_0,\beta_0,\gamma_0,\gamma(\beta_0)),\\
C&=&-\bfb(\beta_0,\gamma_0,\alpha_0,\beta^{-1}(\gamma_0))=-\bfb(\beta_0,\gamma_0,\alpha_0,\alpha(\gamma_0)),\\
A&=&-\bfb(\gamma_0,\alpha_0,\beta_0,\gamma^{-1}(\alpha_0))=-\bfb(\gamma_0,\alpha_0,\beta_0,\beta(\alpha_0)).
\end{eqnarray*}
\end{definition}
We observe that these shear coordinates depend on the choice of the fixed points, hence there are 8 choices of shear coordinates.
We prove now
\begin{proposition}{\sc[lengths and shears]}
Let $\ell_0(\alpha)=\log \bfb(\alpha_0,\alpha(z),\alpha_1,z)$. Then
\begin{eqnarray*}
& &e^{\ell_0(\alpha)}=BC, \ \ e^{\ell_0(\beta)}=AC,\  e^{\ell_0(\gamma)}=AB,\\
& &A=e^{\frac{-\ell_0(\alpha)+\ell_0(\beta)+\ell_0(\gamma)}{2}},\ \
B=e^{\frac{-\ell_0(\beta)+\ell_0(\alpha)+\ell_0(\gamma)}{2}},\ \
C=e^{\frac{-\ell_0(\gamma)+\ell_0(\beta)+\ell_0(\alpha)}{2}}.
\end{eqnarray*}
More generally if $\alpha_0$ and $\alpha_1$ are fixed by $\alpha$ and distinct, then for any points $z$, $t$ and $s$ in $\partial_\infty\pi_1(P)_p$ not fixed by $\alpha$, with $s\not=t$, we have \begin{eqnarray}
\bfb(\alpha_0,\alpha(z),\alpha_1,z)=\bfb(\alpha_0,s,t,\alpha^{-1}(s))\bfb(s,t,\alpha_0,\alpha(t)).\label{l=bc}
\end{eqnarray}
\end{proposition}
\proof It suffices to prove the first equality, the others follow from cyclic permutation. We introduce the following functions of a point $x$ in  $\partial_\infty\pi_1(P)_p$
\begin{eqnarray*}
B\left(x\right)&=&-\bfb\left(x,\beta_0,\gamma_0,\alpha^{-1}\left(\beta_0\right)\right),\\
C\left(x\right)&=&-\bfb\left(\beta_0,\gamma_0,x,\alpha\left(\gamma_0\right)\right).\\
\end{eqnarray*}
By Relation (\ref{projcra}) and (\ref{cr0}), we have
$$
\left(1+\frac{1}{B\left(x\right)}\right)\left(1+C\left(x\right)\right)=\bfb\left(\alpha\left(\gamma_0\right),\gamma_0,x,\beta_0\right)\bfb\left(\gamma_0,x,\beta_0,\alpha^{-1}\left(\beta_0\right)\right).
$$
Using the symmetry (\ref{projsym}) as well as  the  second cocycle identity, we get
$$
\left(1+\frac{1}{B\left(x\right)}\right)\left(1+C\left(x\right)\right)=\bfb\left(\gamma_0, \alpha\left(\gamma_0\right) ,\beta_0,\alpha^{-1}\left(\beta_0\right)\right).
$$
Thus the following quantity
$$
\left(1+\frac{1}{B\left(x\right)}\right)\left(1+C\left(x\right)\right)
$$
does not depend on $x$.
We observe that, again using the second identity, 
\begin{eqnarray*}
B\left(\alpha_1\right)&=&e^{-\ell_0\left(\alpha\right)}B\left(\alpha_0\right)=e^{-\ell_0\left(\alpha\right)}B,\\
C\left(\alpha_1\right)&=&e^{-\ell_0\left(\alpha\right)}C\left(\alpha_0\right)=e^{-\ell_0\left(\alpha\right)}C.
\end{eqnarray*}
Hence, we have
$$
\left(1+\frac{1}{B}\right)\left(1+C\right)=
\left(1+\frac{e^{\ell_0\left(\alpha\right)}}{B}\right)\left(1+\frac{C}{e^{\ell_0\left(\alpha\right)}}\right).
$$
Finally, we remark that the equation
$$
\left(1+\frac{1}{B}\right)\left(1+C\right)=\left(1+\frac{u}{B}\right)\left(1+\frac{C}{u}\right),
$$
 is quadratic in $u$ and its two obvious solutions are $u=1$ and $u=BC$. Equation (\ref{l=bc}) is just another way to restate this last equality.
 \qed

\subsection{Gap functions in terms of  shear coordinates}
We now compute the gap function in terms of the shear coordinates and lengths.

\begin{proposition}\label{gapsl2}
We have the following expression of the pant  gap functions
\begin{eqnarray*}
G_\bfb\left(P\right)&=&\log \left(\frac{1+B^{-1}}{1+C} \right)
=\log \left(\frac{e^{\frac{\ell\left(\beta\right)+\ell\left(\gamma\right)}{2}}+e^{\frac{\ell\left(\alpha\right)}{2}}}{e^{\frac{\ell\left(\beta\right)+\ell\left(\gamma\right)}{2}}+e^{-\frac{\ell\left(\alpha\right)}{2}}} \right),\\
G^r_\bfb\left(P\right)
&=&\log \left(\frac{\cosh \left(\frac{\ell\left(\gamma\right)}{2}\right)+\cosh\left(\frac{\ell\left(\beta\right)-\ell\left(\alpha\right)}{2}\right)}{\cosh \left(\frac{\ell\left(\gamma\right)}{2}\right)+\cosh\left(\frac{\ell\left(\beta\right)+\ell\left(\alpha\right)}{2}\right)} \right).
\end{eqnarray*}
\end{proposition} 
\proof
It follows from the first cocycle identity and  Relation (\ref{projcra}) that
\begin{eqnarray*}
\bfb\left(v,f,u,e\right)=\frac{\bfb\left(w,f,u,e\right)}{\bfb\left(w,f,v,e\right)}=\frac{1-\bfb\left(e,f,u,w\right)}{1-\bfb\left(e,f,v,w\right)}.
\end{eqnarray*}
Using the first cocycle identity (\ref{bir11bis}), we get
\begin{eqnarray*}
\bfb\left(v,f,u,e\right)
&=&\frac{1-\bfb\left(e,f,v,w\right)\bfb\left(v,f,u,w\right)}{1-\bfb\left(e,f,v,w\right)}\label{cr3}.
\end{eqnarray*}
Hence, we have
\begin{eqnarray}
\bfb\left(\alpha_0,\gamma_0,\alpha_1,\beta_0\right)&=&\frac{1-\bfb\left(\beta_0,\gamma_0,\alpha_0,\alpha\left(\gamma_0\right)\right)\bfb\left(\alpha_0,\gamma_0,\alpha_1,\alpha\left(\gamma_0\right)\right)}{1-\bfb\left(\beta_0,\gamma_0,\alpha_0,\alpha\left(\gamma_0\right)\right)}\cr
&=&\frac{1+Ce^{-\ell_{0}\left(\alpha\right)}}{1+C}=\frac{1+B^{-1}}{1+C}\cr 
&=&\frac{1+e^{\frac{-\ell_0\left(\gamma\right)+\ell_0\left(\beta\right)-\ell_0\left(\alpha\right)}{2}}}{1+e^{\frac{-\ell_0\left(\gamma\right)+\ell_0\left(\beta\right)+\ell_0\left(\alpha\right)}{2}}}\cr
&=&\frac{e^{\frac{-\ell_0\left(\beta\right)+\ell_0\left(\gamma\right)}{2}}+e^{-\frac{\ell_0\left(\alpha\right)}{2}}}{e^{\frac{-\ell_0\left(\beta\right)+\ell_0\left(\gamma\right)}{2}}+e^{\frac{\ell_0\left(\alpha\right)}{2}}}.\label{eq:hypgap1}
\end{eqnarray}
We obtain our first result by taking $\alpha_0=\alpha^+$, $\beta_0=\beta^+$ and $\gamma_0=\gamma^-$. Indeed, in this case we have
$$
\ell_0\left(\alpha\right)=-\ell\left(\alpha\right),\ \ \ell_0\left(\beta\right)=-\ell\left(\beta\right),\ \  \ell_0\left(\gamma\right)=\ell\left(\gamma\right).
$$
When changing $\beta_0$ to $\beta_1$, we get
\begin{eqnarray}
\bfb\left(\alpha_0,\gamma_0,\alpha_1,\beta_1\right)&=
&\frac{e^{\frac{\ell_0\left(\beta\right)+\ell_0\left(\gamma\right)}{2}}+e^{\frac{\ell_0\left(\alpha\right)}{2}}} 
{e^{\frac{\ell_0\left(\beta\right)+\ell_0\left(\gamma\right)}{2}}+e^{ - \frac{\ell_0\left(\alpha\right)}{2}}}\label{eq:hypgap2}. 
\end{eqnarray}
Hence combining Equations (\ref{eq:hypgap1}) and (\ref{eq:hypgap2}),  and using the second cocycle identity, we get
\begin{eqnarray*}
\bfb\left(\alpha_0,\beta_1,\alpha_1,\beta_0\right)&=&\frac{\left(e^{-\frac{\ell_0\left(\beta\right)+\ell_0\left(\gamma\right)}{2}}+e^{-\frac{\ell_0\left(\alpha\right)}{2}}\right)}{\left(e^{\frac{-\ell_0\left(\beta\right)+\ell_0\left(\gamma\right)}{2}}+e^{\frac{\ell_0\left(\alpha\right)}{2}}\right)}\frac{\left(e^{\frac{\ell_0\left(\beta\right)+\ell_0\left(\gamma\right)}{2}}+e^{\frac{\ell_0\left(\alpha\right)}{2}}\right)}{\left(e^{\frac{\ell_0\left(\beta\right)+\ell_0\left(\gamma\right)}{2}}+e^{-\frac{\ell_0\left(\alpha\right)}{2}}\right)}\cr
&=&\frac{e^{\ell_0\left(\gamma\right)}+1+e^{\frac{-\ell_0\left(\beta\right)+\ell_0\left(\gamma\right)+\ell_0\left(\alpha\right)}{2}}+e^{\frac{+\ell_0\left(\beta\right)+\ell_0\left(\gamma\right)-\ell_0\left(\alpha\right)}{2}}}{e^{\ell_0\left(\gamma\right)}+1+e^{\frac{-\ell_0\left(\beta\right)+\ell_0\left(\gamma\right)-\ell_0\left(\alpha\right)}{2}}+e^{\frac{\ell_0\left(\beta\right)+\ell_0\left(\gamma\right)-\ell_0\left(\alpha\right)}{2}}}\cr
&=&\frac{e^{\frac{-\ell_0\left(\gamma\right)}{2}}+e^{\frac{-\ell_0\left(\beta\right)+\ell_0\left(\alpha\right)}{2}}+e^{\frac{\ell_0\left(\beta\right)-\ell_0\left(\alpha\right)}{2}}+e^{\frac{\ell_0\left(\gamma\right)}{2}}}{e^{\frac{\ell_0\left(\gamma\right)}{2}}+e^{-\frac{\ell_0\left(\beta\right)-\ell_0\left(\alpha\right)}{2}}+e^{\frac{\ell_0\left(\beta\right)+\ell_0\left(\alpha\right)}{2}}+e^{-\frac{\ell_0\left(\gamma\right)}{2}}}\cr
&=&\frac{\cosh \left(\frac{\ell\left(\gamma\right)}{2}\right)+\cosh\left(\frac{\ell\left(\beta\right)-\ell\left(\alpha\right)}{2}\right)}{\cosh \left(\frac{\ell\left(\gamma\right)}{2}\right)+\cosh\left(\frac{\ell\left(\beta\right)+\ell\left(\alpha\right)}{2}\right)}.
\end{eqnarray*}
The proposition follows. \qed
\subsection{Cusp gap functions}
We use in this paragraph the notations of Paragraph \ref{cuspgapfunction}.
\begin{proposition}\label{gapsl2pinch}If the hyperbolic pair of pants $P$ has a cusp at $\alpha$,
we have
\begin{eqnarray*}
W_{\bfb}(P)=\frac{1}{1+e^{\frac{\ell(\beta)+\ell(\gamma)}{2}}}, & &
W_{\bfb}^r(P)=\frac{\sinh\left(\frac{\ell(\beta)}{2}\right)}
{\cosh\left(\frac{\ell(\gamma)}{2}\right)+\cosh\left(\frac{\ell(\beta)}{2}\right)}.
\end{eqnarray*}
\end{proposition}
\proof
We introduce as before the {\em shear coordinates} of $P$ associated to some choice of fixed points:
\begin{eqnarray*}
B&=&-\bfb(\alpha^+,\beta^+,\gamma^-,\alpha^{-1}(\beta^+))=-\bfb(\alpha^+,\beta^+,\gamma^-,\gamma(\beta^+)),\\
C&=&-\bfb(\beta^+,\gamma^-,\alpha^+,\beta^{-1}(\gamma^-))=-\bfb(\beta^+,\gamma^-,\alpha^+,\alpha(\gamma^-)),\\
A&=&-\bfb(\gamma^-,\alpha^+,\beta^+,\gamma^{-1}(\alpha^+))=-\bfb(\gamma^-,\alpha^+,\beta^+,\beta(\alpha^+)).
\end{eqnarray*}
We now have
\begin{eqnarray*}
& &1=BC, \ \ e^{-\ell(\beta)}=AC,\  e^{\ell(\gamma)}=AB,\\
& &A=e^{\frac{-\ell(\beta)+\ell(\gamma)}{2}},\ \
B=e^{\frac{\ell(\beta)+\ell(\gamma)}{2}},\ \
C=e^{\frac{-\ell(\gamma)-\ell(\beta)}{2}}.
\end{eqnarray*}
By the definition in \ref{walpha},
$$
W_{\bfb}(P)=W_\alpha(\gamma^-,\beta^+).
$$
Recall that
$$
R_x(s,t,s_0,t_0)=\frac{\partial_y \log \bfb(x,s,y,t)}{\partial_y \log \bfb(x,s_0,y,t_0)}\bigg\vert_{y=x}.
$$
In projective coordinates, we obtain by a direct computation of derivatives that
$$
R_x(s,t,s_0,t_0)=\frac{\frac{1}{x-t}-\frac{1}{x-s}}{\frac{1}{x-t_0}-\frac{1}{x-s_0}}.$$
Hence
$$
R_x(s,t,s_0,t_0)
=\frac{(x-t_0)(x-s_0)(t-s)}{(x-t)(x-s)(t_0-s_0)}=\bfb(t_0,x,t,s)\bfb(s_0,x,s,t_0).
$$
Therefore, taking $s_0=\gamma^-$, yields by definition of $W_\alpha$ that 
\begin{eqnarray}
W_{\alpha}(\gamma^-,\beta^+)=R_{\alpha^+}(\gamma^-,\beta^+,\gamma^-,\alpha(\gamma^-))=
\bfb(\alpha(\gamma^-),\alpha^+,\beta^+,\gamma^-).\label{W2-first}
\end{eqnarray}
Hence
\begin{eqnarray}
W_{\bfb}(P)
&=&\frac{1}{\bfb(\alpha(\gamma^-),\gamma^-,\beta^+,\alpha^+)}\ \ \hbox{ by (\ref{bir11bis})}\cr &=&\frac{1}{1-\bfb(\alpha^+,\gamma^-,\beta^+,\alpha(\gamma^-)\,)}\hbox{ by (\ref{projcra})}\cr
&=&\frac{1}{1+C^{-1}}\ \ \hbox{ by (\ref{bir11bis})}\cr 
&=&\frac{1}{1+e^{\frac{\ell(\beta)+\ell(\gamma)}{2}}}.\label{cuspgap1} 
\end{eqnarray}
This proves the first equality in the proposition.

Similarly, taking $\beta^-$ instead of $\beta^+$, yields 
\begin{eqnarray}
W_{\alpha}(\gamma^-,\beta^-)=\frac{1}{1+e^{\frac{-\ell(\beta)+\ell(\gamma)}{2}}}.\label{cuspgap2}
\end{eqnarray}
Recall that 
$$
W_{\bfb}^r(P)=W_{\alpha}(\beta^+,\beta^-).
$$
The additive cocycle identity (\ref{addcoc}) yields
$$
W_{\bfb}^r(P)=W_\alpha(\gamma^-,\beta^-)-W_{\alpha}(\gamma^-,\beta^+).
$$
Using Equations (\ref{cuspgap1}) and   (\ref{cuspgap2}), we obtain 
\begin{eqnarray*}
W_{\bfb}^r(P)&=&\frac{1}{1+e^{\frac{-\ell(\beta)+\ell(\gamma)}{2}}}-\frac{1}{1+e^{\frac{\ell(\beta)+\ell(\gamma)}{2}}}\\
&=&\frac{\sinh\left(\frac{\ell(\beta)}{2}\right)}
{\cosh\left(\frac{\ell(\gamma)}{2}\right)+\cosh\left(\frac{\ell(\beta)}{2}\right)}.
\end{eqnarray*}
This last equality concludes the proof.
\qed

\section{Fuchsian and Hitchin representations}

Let $\Sigma$ be a compact oriented connected surface with or without boundary. If $\Sigma$ is closed we assume it has genus at least two; if $\Sigma$ has a boundary, we assume 
that its double along its boundary has genus at least  two.
\begin{definition}{\sc[Fuchsian  and Hitchin homomorphisms]} An {\em $n$-Fuchsian} homomorphism from $\grf$ to $\sln$ is  a
homomorphism $\rho$ which factorises as
$\rho=\iota\circ\rho_{0}$, where $\rho_{0}$ is a convex cocompact discrete faithful 
homomorphism with values in ${\rm PSL}(2,\mathbb R)$ and $\iota$
is an  irreducible homomorphism from  ${\rm PSL}(2,\mathbb R)$ to 
${\rm PSL}(n,\mathbb R)$.  

An {\em $n$-Hitchin homomorphism} from $\grf$ to $\sln$ is a homomorphism which may  be deformed into an $n$-Fuchsian homomorphism in such a way that the image of each  boundary component stays purely loxodromic at each stage of the deformation.
\end{definition}

In the case of a closed  surface, these definitions agree with the definitions in  \cite{Labourie:2006}. Recall that for a closed surface the image of every non trivial element is purely loxodromic, hence the restriction of a Hitchin representation for a closed surface to a subsurface is again Hitchin.

Let $\Sigma$ be a compact connected oriented surface with $k$ boundary components $C_1,\ldots,C_k$. Let  $A_1,\ldots, A_k$ be conjugacy classes of purely loxodromic elements in ${\rm PSL}(n,\mathbb R)$. We denote by 
$${\operatorname{Hom}}_H(\grf,{\rm PSL}(n,\mathbb R); A_1,\ldots, A_k)$$
the space of Hitchin homomorphisms  from $\grf$ to $\sln$ whose holonomy along the boundary component $C_j$ is conjugated to $A_j$. Observe that $\sln$ acts on this space by conjugation.
Let
\begin{eqnarray*}
& &{\operatorname{Rep}}_H(\grf,{\rm PSL}(n,\mathbb R); A_1,\ldots, A_k)\\ & &=\ \ {\operatorname{Hom}}_H(\grf,{\rm PSL}(n,\mathbb R); A_1,\ldots, A_k)/\sln,
\end{eqnarray*}  be {\em the moduli space of Hitchin representations} whose holonomies along the boundary component $C_j$ is conjugated to $A_j$.
\section{Positivity}
\label{sec:FG}
In this section, we  recall  definitions used by Volodia Fock and Sasha Goncharov in Section 9 of  \cite{Fock:2006a}. 

\subsection{Positivity in flag manifolds}\label{sec:posflag}

\begin{definition}{\sc[flags and basis of flags]} A {\em flag} $F$ in $\mathbb R^n$ is a family $(F^1,\ldots , F^{n-1})$ such that $F^k$ is a $k$-dimensional vector space  
and $F^k\subset F^{k+1}$. 
 A {\em basis for the flag  $F$} is 
a basis $\{f_i\}$  of  $\mathbb R^n$
 such that 
$F^k$ is generated by $\{f_1,\ldots,f_k\}$.  We denote by $\mathcal F(\mathbb R^n)$ the space of all  flags.
\end{definition}
\subsubsection{Triple ratios and positive triple of flags}

Let $(F,G,H)$ be a  triple of flags in {\em general position} that is so that  for every triple of positive integers $(m,l,p)$ with $m+l+p=n$ the sum $F^m+G^l+H^p$ is direct. 
Let $\{f_i\}$, $\{g_i\}$, and $\{h_i\}$ be bases respectively for $F$, $G$ and $H$.  We define if $1\leq p \leq n$
$$
\widehat f^{ p}=f_1\wedge\ldots\wedge f_p, 
$$
and by convention $\widehat f^0=1$.

\begin{definition}{\sc[Positivity and $(m,l,p)$-triple ratio]}\label{triratio}The {\em $(m,l,p)$-triple ratio} of $(F,G,H)$ is 
$$
T^{m,l,p}(G,F,H)=\frac{\Omega(\widehat f^{ m+1} \wedge  \widehat  g^{ l}\wedge  \widehat  h^{ p-1})\Omega(\widehat f^{ m-1}
 \wedge  \widehat   g^{ l+1}\wedge  \widehat  h^{ p})\Omega(\widehat f^{ m} \wedge \widehat  g^{ l-1}\wedge \widehat  h^{ p+1})} 
 {\Omega(\widehat f^{ m+1} \wedge \widehat  g^{ l-1}\wedge \widehat  h^{ p})\Omega(\widehat f^{ m} 
 \wedge \widehat  g^{ l+1}\wedge \widehat  h^{ p-1})\Omega(\widehat f^{ m-1} \wedge \widehat  g^{ l}\wedge \widehat  h^{ p+1})},
$$
where $\Omega$ is any volume form on $\mathbb R^n$. 

The triple  $(F,G,H)$ is {\em positive} if the $(m,l,p)$-triple ratios are positive for all $(m,l,p)$.
\end{definition}

Finally, the following result holds by the remark after Lemma 9.1 in \cite{Fock:2006a} in dimension 3 and an  induction in the dimension for the general case

\begin{proposition}
The  collection of functions $T=(T^{m,l,p})_{m+l+p=n}$  defines a homeomorphism from the space of positive triple of flags 
-- up to the action of $\sln$-- to $(\mathbb R^+)^{\frac{(n-1)(n-2)}{2}}$.
\end{proposition}
The reader is at least  encouraged to check that the dimensions are the same. This is a special case and,  as we just said, a step in the proof of Theorem 9.1 of  \cite{Fock:2006a}.
For $n=3$, there is just one triple ratio and we have \begin{proposition}
 The {\em triple ratio} 
of the three  flags   $F_{1}=(L_1,P_1)$, $F_{2}=(L_2,P_2)$ and $F_{3}=(L_3,P_3)$ 
is
 $$
T(F_1,F_2,F_3)=\frac{\langle \widehat L_1\vert  \widehat P_2 \rangle  \langle \widehat L_2\vert  \widehat P_3 \rangle  \langle \widehat L_3\vert  \widehat P_1 \rangle }{\langle \widehat L_1\vert  \widehat P_3 \rangle  \langle \widehat L_3\vert  \widehat P_2 \rangle  \langle \widehat L_2\vert  \widehat P_1 \rangle }.
$$
where $\widehat L_i$  and $\widehat P_i$ are nonzero vectors in 
$L_i$ and    $P_i^\perp$ respectively.
\end{proposition}

\rmk The reader can check from this proposition or the previous definition that the positivity of a triple is invariant under {\em all} permutations. Therefore positivity has nothing to do with orientation.

\subsubsection{Cross ratios and positive quadruple of flags}
We now consider a quadruple of flags $Q=(X,Y,Z,T)$ in $\mathbb R^n$. 
Throughout this paragraph, we suppose that both $Q_1=(X,Y,Z)$ and $Q_2=(T,X,Z)$ are positive triples of flags.
 We associate  $n-1$ numbers to $Q$  which, together with the triple ratios of  $Q_1$ and $Q_2$, 
completely determine the configuration $Q$ up to the action of $\sln$.

As before we choose  a basis $\{x_1,\ldots,x_n\}$ adapted to the flag  $X$  
and likewise for  $Y$, $Z$ and $T$. We say a  quadruple $(X,Y,Z,T)$ is in {\em general position} if the triples $(X,Y,Z)$ and $(T,X,Z)$ are in general position.  
Following \cite{Fock:2006a}, 
\begin{definition}{\sc[Edge functions]}
The  value of the  {\em edge functions} $\delta_i$ for  $i=1, \ldots, n-1$ at a   quadruple in general position $(X,Y,Z,T)$ is 
$$
\delta_i(X,Y,Z,T)=\frac{\Omega(\widehat x^{i}\wedge \widehat z^{n-i-1}\wedge t_1)\Omega(\widehat x^{i-1}\wedge y_1\wedge\widehat z^{n-i})}{\Omega(\widehat x^{i}\wedge y_1 \wedge\widehat z^{n-i-1})\Omega(\widehat x^{i-1}\wedge \widehat z^{n-i}\wedge t_1)},
$$
where $\Omega$ is any volume form on $\mathbb R^n$.
\end{definition}
 For $n=2$, it is  worth 
  noticing that
\begin{eqnarray}\label{deltabir}
\delta_1(X,Y,Z,T)=\frac{\Omega(x_1 \wedge t_1)\Omega( y_1\wedge z_1)}{\Omega( x_1\wedge y_1)\Omega(z_1\wedge t_1)}
=-{\bf b}(X,T,Z,Y).
\end{eqnarray}
We also observe that
\begin{eqnarray}
\delta_1(U,V,W,R)=\frac{1}{\delta_{n-1}(W,V,U,R)}.\label{eq:inv}
\end{eqnarray}
We write  $\Delta$  for  the $(n-1)$-tuple of functions  $(\delta_1,\ldots,\delta_{n-1})$.
Let $\pi$ be the projection from $\EE$ to the  two-dimensional vector space 
 $P=\EE/(X^{i-1}\oplus Z^{n-i-1})$.
  Let ${\bf b}_P$ be the cross ratio in this plane. One deduces from Equation (\ref{deltabir}) that
\begin{eqnarray}
\delta_i(X,Y,Z,T)=-{\bf b}_P(\pi(X^i),\pi(T^{1}),\pi(Z^{n-i})),\pi(Y^{1})).\label{deltacross}
\end{eqnarray} 
\begin{definition}{\sc[Positive quadruple]}
The quadruple of flags $$Q=(X,Y,Z,T),$$  is {\em positive} if it is in general position and  if all the edge functions $\delta_i$ are positive and if both $Q_1=(X,Y,Z)$ and $Q_2=(T,X,Z)$ are positive triples of flags.
\end{definition}

Positivity of quadruples is invariant under cyclic permutations.

An easy induction on dimension yields.
\begin{proposition}\label{FGquad}
The mapping
$$
Q\longrightarrow (T(Q_1),T(Q_2),\Delta(Q)),$$
is a homeomorphism from the space of positive quadruples of flags 
(up to the action of the projective special linear group) into  $(\mathbb R^+)^{(n-1)^2}$. 
\end{proposition}

\subsubsection{Positive maps and representations}
Following \cite {Fock:2006a}, we have 
 \begin{definition}\begin{itemize}\label{def:posrep}
 \item A {\em positive map}
 from a cyclically  ordered set to the space of flags is  a map such that the image of every positively ordered quadruple is a positive quadruple.
\item
A representation from  $\grf$ to $\sln$ is {\em positive}
if there exists a positive continuous $\rho$-equivariant map from $\bgrf_p$ to $\mathcal{ F}(\EE)$ for some $p$. 
\end{itemize}
\end{definition}

One can easily check that the limit curve $\xi$ of a Fuchsian representation is positive.
We shall prove later that all Hitchin representations are positive.
\section{Frenet curves}\label{repsln}

In this section, we define Frenet curves and recall results from \cite{Labourie:2005,Labourie:2006} that link representations of closed surface groups to cross ratios using Frenet curves. We give a simpler expression of the cusp gap function in the case of the cross ratio associated to a Frenet curve.
In the last paragraph, we recall the relation between positivity and Frenet curves.
\subsection{Frenet curves and cross ratios}\label{par:frenet}

\begin{definition}{\sc[Frenet curve and osculating flag]}\label{defrenet}
A   curve $\xi$ defined from $\mathbb T$ to $\RPN$  is a {\em Frenet
  curve} if there exists a curve
$(\xi^{1},\xi^{2},\ldots,\xi^{n-1})$ defined on $\mathbb T$, called the {\em  osculating flag curve},
with values in the flag variety such that for every $x$ in $\mathbb T$, $\xi(x)=\xi^{1}(x)$, and moreover
\begin{itemize}
\item  For every pairwise distinct points $(x_{1},\ldots,x_{l})$ in $\mathbb T$ and positive  integers $(n_1,\ldots,n_l)$ such that  $\sum_{i=1}^{i=l}n_i\leq n$, then the following sum is direct
\begin{eqnarray}
\xi^{n_i}(x_i)+\ldots+\xi^{n_{l}}(x_{l})\label{fre2}.
\end{eqnarray}
\item For every $x$ in $\mathbb T$ and positive  integers $(n_1,\ldots,n_l)$ such that  $p=\sum_{i=1}^{i=l}n_i\leq n$,
  then
\begin{eqnarray}
\lim_{\substack{(y_1,\ldots,y_l)\rightarrow x,\\ y_i
\text{all distinct}}}
(\bigoplus_{i=1}^{i=l}\xi^{n_i}(y_i))=\xi^{p}(x)\label{fre3}.
\end{eqnarray}
\end{itemize}
We call $\xi^{n-1}$ the {\em osculating hyperplane}.
\end{definition}

By Condition (\ref{fre3}), the image of a Frenet curve is  a
$C^{1}$-submanifold and  the tangent line to $\xi^1(x)$ is $\xi^2(x)$.
 Moreover,  a Frenet curve $\xi$ is {\em hyperconvex} in the following sense: for any $n$-tuple of pairwise distinct points
  $(x_1,\ldots,x_{n})$ we have 
$$
\xi(x_1)+\xi(x_2)+\ldots+\xi(x_n)=\mathbb R^n.
$$

\begin{definition}{\sc[Associated weak cross ratio]}  Let $\xi$ be a Frenet curve and $\xi^*$ be  its associated osculating hyperplane curve.  
The  {\em weak cross ratio} associated to this pair of curves is the function on $\mathbb T^{4*}$  defined by 
$$
\bb_{\xi,\xi^*}(x,y,z,t)=\frac{\langle\widehat\xi(x)\vert\widehat\xi^*(y)\rangle\langle\widehat\xi(z)\vert\widehat\xi^*(t)\rangle}
{\langle\widehat\xi(z)\vert\widehat\xi^*(y)\rangle\langle\widehat\xi(x)\vert\widehat\xi^*(t)\rangle},
$$ 
where for every $u$, we choose an arbitrary nonzero vector $\widehat\xi(u)$ and  $\widehat\xi^*(u)$ in  respectively $\xi(u)$ and $\xi^*(u)$.  
\end{definition}

Observe that the weak cross ratio associated to a Frenet curve is not necessarily a cross ratio. However, we shall see that it is so for Frenet curves  that we shall be interested in. 

\rmks
\begin{enumerate}
\item Let $V=\xi(x)\oplus\xi(z)$. Let ${\eta}(m)=\xi^*(m)\cap V$. Let ${\bf b}_V$ be the classical cross ratio on $\mathbb P(V)$, then
\begin{eqnarray}
\bb_{\xi,\xi^*}(x,y,z,t)={\bf b}_V(\xi(x),{\eta}(y),\xi(z),{\eta}(t)).\label{crossinter}
\end{eqnarray}
\item Symmetrically, let $H=\xi^*(y)\cap \xi^*(t)$, let $\pi$ be the projection on $P=\mathbb R^n/H$, let ${\bf b}_P$ be the classical cross ratio on $\mathbb P(P)$, then
\begin{eqnarray}
\bb_{\xi,\xi^*}(x,y,z,t)={\bf b}_P(\pi(\xi(x)),\pi(\xi^*)(y)),\pi(\xi(z)),\pi(\xi^*(t))).\label{crosspi}
\end{eqnarray}
\end{enumerate}

\subsection{Frenet curves and representations}\label{cr-rep}
The following summarises some of the results of \cite{Labourie:2006,Labourie:2005}.

\begin{theorem}\label{mainA}{\sc[F. Labourie]} 
Let   $\rho$ be an $n$-Hitchin representation of the fundamental group of  a closed connected oriented surface of genus at least two. 
Then,  there exists a $\rho$-equivariant Frenet curve from $\partial_{\infty}\grf$ to
$\mathbb P(\mathbb R^{n})$. Moreover, the osculating flag  curve is Hölder. 

Finally, the weak cross ratio $\bb_\rho=\bb_{\xi,\xi^*}$ is a cross ratio -- called  
the {\em cross ratio associated to  the  Hitchin representation $\rho$} -- and the period of $\gamma$ is given by
$$
\ell_\bb(\gamma)=\log\left(\left\vert\frac{\lambda_{\max}(\rho(\gamma))}{\lambda_{\min}(\rho(\gamma))}\right\vert\right),
$$
where $\lambda_{\max}(\rho(\gamma))$ and $\lambda_{\min}(\rho(\gamma))$ are the  eigenvalues  of respectively maximum and minimum absolute values of the element $\rho(\gamma)$.\end{theorem}
The first main result in \cite{Labourie:2005}  characterises  which cross ratios appear in this theorem. 
In Corollary \ref{hitchindouble}, we will show that every Hitchin representation $\rho$ for a surface $\Sigma$  with boundary is the restriction of  a Hitchin representation, called the Hitchin double,  for the double surface. It follows that we can also associate to such a representation a cross ratio: the restriction of the cross ratio associated to the Hitchin double by Theorem \ref{hcintro}.
Finally,  the following observation follows from the construction and from the density of attracting points of elements of the fundamental group.

\begin{proposition}\label{gapsubsurface}
Let $S$ be a closed surface. Let $\Sigma$ be an incompressible connected surface embedded in $S$. Let $\rho_0$ and $\rho_1$ be two Hitchin representations of 
$\pi_1(S)$ whose restrictions to $\grf$ coincide.

Let $\bb_{\rho_0}$ and $\bb_{\rho_1}$ be the associated cross ratios on $\partial_\infty\pi_1(S)$. Then $\bb_{\rho_0}$ and $\bb_{\rho_1}$ coincide on $\bgrf$. 

In particular, the gap functions for a pair of pants defined by the triple $(\alpha,\beta,\gamma)$  of elements of $\grf$ and a cross ratio associated to a Hitchin representation $\rho$ only depend on the triple $(\rho(\alpha),\rho(\beta),\rho(\gamma))$ of elements of $\sln$. 
\end{proposition}

\subsection{Frenet curves and cusps}
 In this paragraph, we compute  the cusp gap function whenever the cross ratio arises from a Frenet curve. 
For the sake of simplicity, we consider a surface $\Sigma$ with only one cusp. Let $\rho$ be  a representation from $\pi_1(\Sigma)$ to ${\rm PSL}(\EE)$. 
\vskip 0.1 truecm
{\em We assume that there exists a $\rho$-equivariant Frenet curve $\xi$ from $\bgrf_1$ to $\mathbb P(\EE)$ with Hölder osculating flag curve.} .

Note that we do not know if this is the case in general except for Fuchsian representation, since the results of the first author do not cover the cusp case.

\vskip 0.1 truecm
Let $\bb_\rho=\bb_{\xi,\xi^*}$ be the associated cross ratio on $\bgrf_1$. Since $\xi$ is Frenet, there exists  a $C^1$ structure on $\bgrf_1$ such that the map
$
(x,y)\mapsto \bb_\rho(x,s,y,t),
$
is $C^1$ with Hölder derivatives along the diagonal. In other words, the  associated cross ratio satisfies the
 regularity hypothesis \ref{reghyp}.

We evaluate  in this paragraph the quantity $W_\alpha(s,t)$ defined in Paragraph \ref{defcusp}.
Let $F=(L,P)$ be a two-flag in a vector space, that is  a pair $(L,P)$ where $L$ is a line included in a two-plane $P$.
Let $S,T,S_0,T_0$ be four hyperplanes. 
We denote by $\bb$ the classical  cross ratio in $\mathbb P(P)$. We also denote hyperplanes by uppercase letters and their intersection with $P$ by the corresponding  gothic lowercase so that  ${\mathfrak h}=H\cap P$.
Then we define for generic hyperplanes,
$$
\widehat W_F(S,T,S_0,T_0)=\bb({\mathfrak t},{\mathfrak s},L,{\mathfrak s}_0)\bb({\mathfrak t}_0,L,{\mathfrak t},{\mathfrak s}_0). 
$$
We now prove
\begin{theorem}\label{cuspcom}
Let  $(\xi^{1},\xi^{2},\ldots, \xi^{n-1})$ be the osculating flag curve  of $\xi$. 
 Let $F(\alpha)=(\xi^1(\alpha^+),\xi^2(\alpha^+))$. Then
$$
W_\alpha(s,t)=\widehat W_{F(\alpha)}\left(\xi^{n-1}(s),\xi^{n-1}(t),  \xi^{n-1}(s_0),\xi^{n-1}(\alpha(s_0))\right). 
$$
Moreover, when $n=2$,
\begin{eqnarray}
W_\alpha(s,t)=\bb_{\rho}(\alpha(s),\alpha^+,t,s).\label{W2}
\end{eqnarray}
\end{theorem}
Let us make the following remarks when $n=2$
\begin{itemize}
\item
One recover directly from Formula (\ref{W2}) that when $n=2$
$$
W_\alpha(s,t)=W_\alpha(s,u)+W_\alpha(u,t).
$$
Indeed, we choose projective coordinates so that $\alpha^+=+\infty$. In this case $\alpha$ is the translation by a constant $\tau$ and
$$
W_\alpha(s,t)=\frac{t-s}{\tau}.
$$
\item  Formula (\ref{W2}) coincides  with Formula (\ref{W2-first}) obtained through a direct computation.
 \end{itemize}
\subsubsection{A preliminary proposition}
We  introduce some notations and definitions. Let  $u_1,u_2$ be lines in $E$  and  $P_1,P_2$ be hyperplanes, we define  
$$
{\bf b}(u_1, P_1, u_2, P_2)=\frac{\langle\widehat u_1\vert\widehat P_1 \rangle\langle\widehat u_2\vert\widehat P_2\rangle}{\langle\widehat u_1\vert\widehat P_2 \rangle\langle\widehat u_2\vert\widehat P_1 \rangle},
$$
where  $\widehat u_i$ and  $\widehat P_i$ are nonzero vectors in $u_i$ and $P_i^\perp$ respectively.

Now let $(L,P)$ be a two-flag, let $\ddd L$ be a line so that $P$ is generated by $L$ and $\ddd L$. 

\begin{proposition}
We have
\begin{eqnarray}
\widehat W_F(S,T,S_0,T_0)&=&\frac{\bfb(L,S_0,\ddd L, T)-\bfb(L,S_0,\ddd L, S)}{\bfb(L,S_0,\ddd L, T_0)-1}.
\end{eqnarray}
\end{proposition}
\proof  We denote by $A$ the right hand term. Let $\bb$ be the cross ratio in $\mathbb P (P)$. 
We again denote hyperplanes by uppercase letters and their intersection with $P$ by the corresponding  gothic lowercase so that  ${\mathfrak h}=H\cap P$.
 Since $L$ and $\ddd L$ generate $P$, for all hyperplanes $U$ and $V$, we have similarly to  Equation (\ref{crossinter})
$$
\bfb(L,U,\ddd L, V)=\bb(L,{\mathfrak u},\ddd L, {\mathfrak v}).
$$
By applying first Relation (\ref{cr0}), then the first cocycle identity, 
we get 
$$
1-\bfb(L,S_0,\ddd L, T_0)=\bb(\mathfrak t_0,\mathfrak s_0,\ddd L,L)=\bb(\mathfrak t_0,\mathfrak s_0,\mathfrak t,L)\bb(\mathfrak t,\mathfrak s_0,\ddd L,L).
$$
Moreover, using first the second cocycle identity then Relation (\ref{cr0}), we get 
\begin{eqnarray*}
\bb(L,\mathfrak s_0,\ddd L,\mathfrak s)-\bb(L,\mathfrak s_0,\ddd L, \mathfrak t)&=& \bb(L,\mathfrak s_0,\ddd L, \mathfrak s)\left(1-\bb(L,\mathfrak s,\ddd L, \mathfrak t)\right)\\&=&\bb(L,\mathfrak s_0,\ddd L, \mathfrak s)\bb(\mathfrak t,\mathfrak s,\ddd L, L).
\end{eqnarray*}
Hence,
$$
A=\left(\frac{\bb(L,\mathfrak s_0,\ddd L,\mathfrak  s)}{\bb(\mathfrak t_0,\mathfrak s_0,t,L)}\right)\left(\frac{\bb(\mathfrak t,\mathfrak s,\ddd L, L)}{\bb(\mathfrak t,\mathfrak s_0,\ddd L,L)}\right).
$$
Applying the second cocycle identity to the second factor, we obtain 
$$
A=\frac{\bb(L,\mathfrak s_0,\ddd L, \mathfrak s)}{\bb(\mathfrak t,\mathfrak s_0,\ddd L,\mathfrak s)\bb(\mathfrak t_0,\mathfrak s_0,\mathfrak t,L)}=\bb(\mathfrak t,\mathfrak s,L,\mathfrak s_0)\bb(\mathfrak t_0,L,\mathfrak t,\mathfrak s_0).
$$
The result follows\qed
\subsubsection{Proof of the theorem} 
With $t_{0} = \alpha(s_{0}),\,x = \alpha^{+}$, 
let 
$S=\xi^{n-1}(s)$, $S_0=\xi^{n-1}(s_0)$, $T=\xi^{n-1}(t)$, $T_0=\xi^{n-1}(t_0)$, $L=\xi(x)$, $L+\ddd L=\xi^2(x)$. 

Let also  $\widehat y$ be a nonzero vector of $\xi(y)$,
$z$ be a nonzero vector of  $\ddd L$, $(\widehat s_0,\widehat t_0, \widehat s, \widehat t)$ be linear forms 
whose kernels are respectively 
$(S_0,T_0,S,T)$.
Then, by Equation (\ref{bir12}),
$$
R_x(s,t,s_0,t_0)=\frac{\partial_y \log \bb_{\rho}(x,s,y,t)}{\partial_y \log \bb_{\rho}(x,s_0,y,t_0)}\bigg\vert_{y=x}
=\lim_{y\rightarrow x}\frac{\log (\bb_{\rho}(x,s,y,t))}{\log (\bb_{\rho}(x,s_0,y,t_0))}.
$$
Using that $\lim_{s\rightarrow 1}\frac{\log s}{1-s}=1$, we get
\begin{eqnarray*}
R_x(s,t,s_0,t_0)
&=&\lim_{y\rightarrow x}\frac{1-\bb_{\rho}(x,s,y,t)}{1-\bb_{\rho}(x,s_0,y,t_0)}\\
&=&\lim_{y\rightarrow x}\frac{1-\bir{\widehat x}{\widehat s}{\widehat x}{\widehat t}\bir{\widehat y}{\widehat t}{\widehat y}{\widehat s}}{1-\bir{\widehat x}{\widehat s_0}{\widehat x}{\widehat t_0}\bir{\widehat y}{\widehat t_0}{\widehat y}{\widehat s_0}}\\&=&
\lim_{y\rightarrow x}{\bir{\widehat x}{\widehat t_0}{\widehat x}{\widehat t}}{\bir{\widehat y}{\widehat s_0}{\widehat y}{\widehat s}}\left(\frac{\scal{\widehat x}{\widehat t}\scal{\widehat y}{\widehat s}-\scal{\widehat x}{\widehat s}\scal{\widehat y}{\widehat t}}{\scal{\widehat x}{\widehat t_0}\scal{\widehat y}{\widehat s_0}-\scal{\widehat x}{\widehat s_0}\scal{\widehat y}{\widehat t_0}}\right)\\&=&
{\bir{\widehat x}{\widehat t_0}{\widehat x}{\widehat t}}{\bir{\widehat x}{\widehat s_0}{\widehat x}{\widehat s}}\left(\frac{\scal{\widehat x}{\widehat t}\scal{z}{\widehat s}-\scal{\widehat x}{\widehat s}\scal{z}{\widehat t}}{\scal{\widehat x}{\widehat t_0}\scal{z}{\widehat s_0}-\scal{\widehat x}{\widehat s_0}\scal{z}{\widehat t_0}}\right)\\&=&{\bir{\widehat x}{\widehat t_0}{\widehat x}{\widehat t}}{\bir{z}{\widehat t}{z}{\widehat t_0}}\left(\frac{\bir{\widehat x}{\widehat t}{z}{\widehat t}\bir{z}{\widehat s}{\widehat x}{\widehat s}-1}{\bir{\widehat x}{\widehat t_0}{\widehat x}{\widehat s_0}\bir{z}{\widehat s_0}{\widehat x}{\widehat s_0}-1}\right).
\end{eqnarray*}
Using the definition of $\bfb$, we get
\begin{eqnarray*}
R_x(s,t,s_0,t_0)&=&\bfb(L,T_0,\ddd L, T)\frac{1-\bfb(L,T,\ddd L, S)}{1-\bfb(L,T_0,\ddd L, S_0)}.
\end{eqnarray*}
Applying the second cocycle identity thrice  yields
\begin{eqnarray*}
R_x(s,t,s_0,t_0)&=&\frac{\bfb(L,S_0,\ddd L, T)(1-\bfb(L,T,\ddd L, S))}{\bfb(L,S_0,\ddd L, T_0)(1-\bfb(L,T_0,\ddd L, S_0))}\\
&=&\frac{\bfb(L,S_0,\ddd L, T)-\bfb(L,S_0,\ddd L, S)}{\bfb(L,S_0,\ddd L, T_0)-1}.
\end{eqnarray*}
The first part of the result follows from the previous proposition. We obtain the result for $n=2$, by taking $s=s_0$.
\qed

\subsection{Frenet curves and positivity}\label{Frenetarepositive}
We state the elementary Lemma \ref{weylpos} which shows that the osculating flag curve to a Frenet curve is positive.

\begin{definition}{\sc[Compatible flags]}\label{def:compa} Two flags  $F_1,F_2$ are {\em transverse} if the sum
 $F_{1}^k + F_{2}^{n-k}$ is direct for all $1\leq k \leq n$. If $F_1$ and $F_2$ are transverse flags, we 
say a flag  $F$ is {\em  compatible} with  $F_1$ with respect to  $F_2$  -- or in short with $(F_1,F_2)$ --
if  there exist an integer $p$, such that 
\begin{eqnarray*}
k\leq p &\implies & F^k= F_1^{k},\\
k> p &\implies & F^k= F_1^{p}\oplus F_2^{k-p}. 
\end{eqnarray*}
\end{definition}

The following  lemma will  follow easily from the definitions

\begin{lemma}\label{weylpos}{\sc[Frenet curve and positivity]}
Let $\xi$:\, $\mathbb T \rightarrow \mathbb P (\EE)$
be  a  Frenet curve. Then its osculating flag curve $\widehat\xi$ is positive. 
More generally \begin{itemize}
\item 
Let $K\subset \mathbb T$ be a closed set.
\item Let $G=\bigcup_{n\in \mathbb N}\{x^+_i,x^-_i\}$  be a  collection of points of $K$.
\item 
 Let $K_0= K\setminus G$.
 \end{itemize}
Assume that, for every $i$,  $K_0$ lies in one of the connected component of $\mathbb T\setminus \{x^+_i,x^-_i\}$. 
For each  $i$, let  $F_i$  be a compatible flag with  $(\widehat\xi(x^+_i),\widehat\xi(x^-_i))$.
 Then the following map is positive : 
$$
\left\{
\begin{array}{rcl}
K_0\cup_{n\in \mathbb N} \{x^+_i\}&\rightarrow&\mathcal F,\\ 
x_i^+&\mapsto& F_i,\\
x\in K_0&\mapsto& \widehat \xi(x).
\end{array}
\right.
$$
\end{lemma}

\rmk 
\begin{itemize}
\item This lemma is an immediate consequence of Proposition \ref{posquadflag} which we prove in an Appendix.
\item Conversely, we may ask whether a positive curve is always the osculating flag curve of a Frenet curve. This is indeed the case if the positive curve is smooth enough (see  Theorem 1.8 of Fock and Goncharov in \cite{Fock:2006a}).  
However, this does not imply directly that a positive representation preserves a Frenet curve since the osculating flag curve 
 is only Hölder unless the representation is Fuchsian.
\item The related discussion in Section 9 of  \cite{Fock:2006a} covers the case  $K=\mathbb T$ and this general case could also be deduced from general discussion about positivity as in G. Lusztig work \cite{Lusztig:1994}.  \end{itemize}

\section{Hitchin representations for open surfaces}\label{sec:hitchpos}

The aim of this section is to extend Theorem \ref{mainA}
 to representations 
of  the fundamental group of surfaces with boundary and to prove the following

\begin{theorem}\label{hitchinbord}
Let $\Sigma$ be a compact connected oriented surface with or without boundary components, but without cusps. Let $\rho$ be an $n$-Hitchin representation of $\grf$. Then there exists a positive $\rho$-equivariant Hölder map from $\partial_\infty\grf$ to the flag manifold, which is furthermore the restriction of the osculating flag curve of a Frenet curve. Moreover, for all nontrivial elements $\gamma$ in $\grf$, $\rho(\gamma)$ is purely loxodromic.

In particular, every Hitchin representation is a positive  representation. 
Furthermore, let $S$ be  an incompressible surface embedded in $\Sigma$.
Then $\rho$ restricted to $\pi_1(S)$ is a positive representation. 

Finally, there exists an ordered cross ratio $\bb$ on $\partial_\infty\grf$, such that the  period of any non trivial element $\gamma$ satisfies
$$
\ell _\bb(\gamma)=\log\left(\left\vert\frac{\lambda_{\max}(\rho(\gamma))}{\lambda_{\min}(\rho(\gamma))}\right\vert\right),
$$
where $\lambda_{\max}(\rho(\gamma))$ and $\lambda_{\min}(\rho(\gamma))$ are the  eigenvalues  of respectively maximum and minimum absolute values of the element $\rho(\gamma)$.
\end{theorem}

Observe that, as a corollary, this result yields better information about positives curves appearing in \cite{Fock:2006a} for positive representations: these positive curves are Hölder and are the restriction of the osculating flag curve of a Frenet curve. 

Theorem \ref{hitchinbord} is a consequence  of  the "doubling" Corollary \ref{hitchindouble} of Theorem \ref{mainA} which shows that every Hitchin representation is the restriction of a Hitchin representation of the fundamental group of a closed surface.
 
 We prove Theorem \ref{hitchinbord} in Paragraph \ref{hitchindoubleproof}. The main part of this section is devoted to the doubling construction and the proof of Theorem \ref{mainA}.

\subsection{Homomorphisms and boundary components}

\subsubsection{Good homomorphisms}

Let $\Sigma$ be a compact surface with boundary. Every boundary component is supposed to be oriented and we abusively identify each boundary component with a loop.

We introduce  technical definitions whose purpose is to define "better" homomorphisms in a given conjugacy class by fixing some boundary data. 

\begin{definition}{\sc[Good homomorphism]}
A {\em good homomorphism} based at a point $v$ of a boundary component $\partial_v$ is a homomorphism ${\goodRepeta}$   from $\pi_1\left(\Sigma,v\right)$ into $\sln$, so that \begin{itemize}
\item  The matrix  ${\goodRepeta}\left(\partial_{v}\right)$ is  diagonal  with decreasing entries.
\item The image of all boundary components are purely loxodromic.
\end{itemize}
\end{definition}

The purpose of this definition is to study Hitchin representations and we say

\begin{definition}{\sc[Good representative]}
A {\em good representative}, based at a point $v$ of a boundary component $\partial_v$, of a Hitchin representation $\rho$   is a good homomorphism ${\goodRepeta}$, based at $v$, from $\pi_1\left(\Sigma,v\right)$ into $\sln$, which belongs to the conjugacy class defined by  $\rho$.
\end{definition}

\subsubsection{Changing boundary components}

Given $\rho$ and $v$, a good representative is uniquely defined up to conjugation by a diagonal matrix.

The next proposition explains how good representatives based at different  boundary components are related.

Let ${\mathcal D}$ be the group of diagonal matrices. We denote by  $\comp$ the composition of paths, so that $a\comp b$ is the path $a$ followed by $b$. 
Let $\mathcal B_v$ be the set of  classes of paths joining $v$ to a boundary component up to the following equivalence: the path $a$ joining $v$ to a boundary component $\partial_w$ is equivalent to a path  $b$ joining $v$ to  $\partial_w$ if there exists a path $d$ along $\partial_w$ so that $a\comp d$ is homotopic to $b$.  

We now prove

\begin{proposition}\label{exisdoteta}
Let ${\goodRepeta}$ be a good homomorphism. Then there exists
a unique map
$$
\ddd{\goodRepeta}:\mathcal B_v\times{\mathcal D} \rightarrow {\rm PSL}\left(n,\mathbb R\right),
$$
such that the following properties hold
\begin{itemize}
\item For every path $c$, there exists an element  ${\rm K}_c$ of $\sln$ so that for any diagonal matrix  $\Delta$, 
$$
\ddd{\goodRepeta}(c,\Delta)= {\rm K}_c\Delta{\rm K}_c^{-1}.$$
\item Let $c$ be an arc joining $v$ to a point $w$ in another boundary component $\partial_w$ of  $S$. Let ${\bf F}$ be any  good homomorphism based at $w$ conjugated to ${\goodRepeta}$, then
\begin{eqnarray}
{\goodRepeta}\left(c\comp \partial_w\comp  c^{-1}\right)=\ddd{\goodRepeta}\left(c,{\bf F}\left(\partial_w\right)\right).\label{ghv}
\end{eqnarray}
\end{itemize}
\end{proposition}

Note that, fixing $c$, the map $\Delta\mapsto \ddd{\goodRepeta}(c,\Delta) $ is a conjugation between two maximal split tori of $\sln$.

\proof Let $c$ be an element of $\mathcal B_v$ that is a path joining $v$ to a point $w$ in a boundary component $\partial_w$ of  $S$.
We observe that 
\begin{eqnarray*}
{\goodRepeta}\left(c\comp \partial_w\comp  c^{-1}\right)&=&{\rm K}\diagon_w {\mathrm K}^{-1},
\end{eqnarray*}
where $\diagon_w$ is a  diagonal matrix with decreasing entries and ${\rm K}$ is an element of ${\rm PSL}\left(n,\mathbb R\right)$ well defined up to multiplication by a diagonal matrix on the right. It follows that for every diagonal matrix $\diagon$
$$
\ddd{\goodRepeta}\left(c,\diagon\right):={\rm K}\diagon {\rm K}^{-1},
$$
depends only on the equivalence class of $c$ and $\diagon$. We observe that this map $\ddd{\goodRepeta}$ satisfies the properties of the proposition and is characterised by them.
\qed

The following  proposition summarises the properties of the map that we have just constructed.
\begin{proposition}\label{propGoodRepeta}
With the notation of the previous proposition, the map $\ddd{\goodRepeta}$  enjoys the following properties.
\begin{enumerate}
\item For all loop $\gamma$ based at $v$, all path $c$ and all  diagonal matrix $\diagon$
 \begin{eqnarray}
\ddd{\goodRepeta}\left(\gamma\comp c,\diagon\right)={\goodRepeta}\left(\gamma\right)\ddd{\goodRepeta}\left(c,\diagon\right){\goodRepeta}\left(\gamma\right)^{-1}.\label{ghg}
\end{eqnarray}
\item Fixing a path  $c$ and a diagonal matrix $\Delta$, the map ${\goodRepeta}\mapsto \ddd{\goodRepeta}(c,\Delta)$ is a continuous map from the space of good homomorphims to $\sln$.
\item If $\diagon$ is a diagonal matrix
\begin{eqnarray}
\ddd{\diagon{\goodRepeta}\diagon^{-1}}&=&\diagon\ddd{\goodRepeta}\diagon^{-1},\label{ghd}\\
\ddd{\goodRepeta}\left(c,\diagon\right)^{-1}&=&\ddd{\goodRepeta}\left(c,\diagon^{-1}\right).\label{invdoteta}
\end{eqnarray}
\end{enumerate}

\end{proposition}

\proof We fix  a path $c$ joining $v$  to a point $w$ in a boundary component $\partial_w$ as well as a good homomorphism $\bf F$ based at $w$ and conjugated to $\goodRepeta$.

One then checks that Equation (\ref{ghg}) holds for $\Delta={\bf F}(\partial_w)$ hence for all diagonal matrices $\Delta$ by the construction of $\ddd{\goodRepeta}$. The continuity follows a similar argument.  The last properties follow easily from the construction.
\qed

\subsection{The doubling construction}
We present the main doubling construction. Let $\Sigma$ be a closed surface with boundary. Let  $\widehat \Sigma$ be its double. Let  $j_0$  and $j_1$ be the two injections of $\Sigma$ in $\widehat \Sigma$, and $j:x\mapsto\bar x$ the involution of $\widehat \Sigma$ fixing all points in the boundary of $\Sigma$ such that $j\circ{j_0}=j_1$. 

\subsubsection{Extensions of good homomorphisms}
We first define

\begin{definition} Let $\invol$ be an involution of $\operatorname{PGL}(n,\mathbb R)$ commuting with all diagonal matrices. Let $\goodRepeta$ be a good homomorphism from $\pi_1(\Sigma)$ into $\sln$. A {\em $\invol$-extension} of $\goodRepeta$ is a homomorphism $\widehat\goodRepeta$ from $\pi_1(\widehat\Sigma)$ into $\sln$ so that 
\begin{enumerate}
\item We have \begin{eqnarray}
\widehat{\goodRepeta}\circ (j_0)_*={\goodRepeta}.
\end{eqnarray}
\item For all element $\gamma$ in $\pi_1(\widehat\Sigma)$, we have \begin{eqnarray}
\widehat{\goodRepeta} \left(\bar \gamma\right)=\invol \widehat{\goodRepeta} \left(\gamma\right)  \invol .\label{Jext2}\end{eqnarray}
\item For all arc $c$ 
 in $\Sigma$ joining $v$ to a boundary component of $\Sigma$,  we have
\begin{eqnarray}
\widehat{\goodRepeta}\left(c\comp {\bar c}^{-1}\right)= \ddd{\goodRepeta}\left(c,\invol \right)\invol,
\label{Jext3}\end{eqnarray}
where $\ddd{\goodRepeta}$ is defined in Proposition \ref{exisdoteta}.
\end{enumerate}

\end{definition}

The following proposition  shows that we have a uniquely defined extension in the case of good homomorphims.

\begin{proposition}\label{goodhomo0}

Let $\Sigma$ be a compact surface with non empty  boundary. Let  $\widehat \Sigma$ be its double.  Let $\invol$ be an involution of $\operatorname{PGL}(n,\mathbb R)$ commuting with all diagonal matrices. Then any good homomorphism admits a unique $\invol$-extension.

Finally the map which associates to a good homomorphism its $\invol$-extension is continuous.
\end{proposition}
\proof The uniqueness follows from the fact that $\pi_1(\widehat \Sigma)$ is generated by the groups  $(j_0)_*\pi_1\left(\Sigma\right)$, and $(j_1)_*\pi_1\left(\Sigma\right)$ as well as the homotopy classes of loops of the form $c\comp{\bar c}^{-1}$ where $c$ describes the set of arcs joining $v$ to another boundary component.

We now prove the existence of the $\invol$-extension. We first recall the description of  the fundamental group of the double in combinatorial terms. Let $\partial_0,\ldots,\partial_m$ be the boundary components of  $\Sigma$ such that $v$ belongs to $\partial_0$. Let $c_1,\ldots,c_m$ be arcs joining $v$ to the boundary components $\partial_1,\ldots,\partial_m$. Let $\partial_i^c$ be the elements of $\pi_1\left(\Sigma,v\right)$ given by
$$
\partial_i^c=c_i\comp \partial_i\comp  c_i^{-1}.
$$
Let  $F_m$ be the free group on $m$ generators $x_1,\ldots,x_m$. Let $$G=\pi_1\left(\Sigma\right)\star\pi_1\left(\Sigma\right)\star F_m,$$  where $A\star B$ denotes the free product of $A$ and $B$. Let $i_0$, and $i_1$ be the injections of $\pi_1\left(\Sigma\right)$ in $G$ given by the first and second factor respectively.
Let $H$ be the group normally generated by $i_0\left(\partial_0\right)^{-1}\cdot i_1\left(\partial_0\right)$ and  the elements
$$
i_0\left(\partial_p^c\right)^{-1}\cdot x_p\cdot  i_1\left(\partial_p^c\right)\cdot  x_p^{-1},
$$
for ${1\leq p\leq m}$. Let
$$
\Gamma=G/H.
$$
We identify $x_p$ with its image in $\Gamma$ and $i_0$ and $i_1$ with injections of $\grf$ in $\Gamma$.
Let $i$ be the  involution of $G$ such that for all $p$ in $\{1,\ldots,m\}$
\begin{eqnarray*}
i\left(x_p\right)&=&x_p^{-1},\\
i\circ i_0&=&i_1.
\end{eqnarray*}
Observe that $i$ preserves $H$ so that $i$ gives rise to an involution also denoted $i$ on $\Gamma$. We now recall that  $\Gamma$ is isomorphic with $\pi_1(\widehat \Sigma)$  and that the isomorphism is given by the map $\varphi$ defined by
\begin{eqnarray*}
\varphi\left(x_p\right)&=&c_p\comp( \bar c_p)^{-1},\\
\varphi\circ i_k&=&(j_k)_*,\\
\varphi\circ i&=& \bar\varphi,
\end{eqnarray*}
for all $p$ in $\{1,\ldots,m\}$ and $k$ in $\{0,1\}$.
We identify once and for all $\pi_1(\widehat \Sigma)$ with $\Gamma$ using $\varphi$.

We now construct $\widehat{\goodRepeta}$. Let ${\goodRepeta}$ be a good homomorphism of $\rho$. Using the notation of Proposition \ref{goodhomo}, we define $\widetilde{\goodRepeta}$ as the morphism from $G$ to $\sln$ uniquely characterised by
\begin{eqnarray*}
\widetilde{\goodRepeta}\circ i_0&=&{\goodRepeta},\\
\widetilde{\goodRepeta}\circ i_1&=&\invol {\goodRepeta} \invol  ,\\
\widetilde{\goodRepeta} \left(x_p\right)&=&\ddd{\goodRepeta}\left(c_p,\invol \right)\invol,
\end{eqnarray*}
for all $p$ in $\{1,\ldots,m\}$.
We now prove that $\widetilde{\goodRepeta}$ vanishes on elements of $H$. First, we have for all $p$ in $\{1,\ldots,m\}$
\begin{eqnarray*}
\widetilde{\goodRepeta}\left(i_0\left(\partial_p^c\right)^{-1}\cdot  x_p\cdot  i_1\left(\partial_p^c\right)\cdot  x_p^{-1}\right)&=&
\widetilde{\goodRepeta}\left(i_0\left(\partial_p^c\right)\right)^{-1}\cdot \widetilde{\goodRepeta}\left( x_p\right)\cdot \widetilde{\goodRepeta} \left( i_1\left(\partial_p^c\right)\right)\cdot  \widetilde{\goodRepeta}\left(x_p^{-1}\right)\\
&=&{\goodRepeta}\left(\partial_p^c\right)^{-1}\ddd{\goodRepeta}\left(c_p,\invol \right){\goodRepeta}\left(\partial_p^c\right) \ddd{\goodRepeta}\left(c_p,\invol \right)^{-1}.
\end{eqnarray*}
Recall now from the definition of $\ddd{\goodRepeta}$  that we have for every $p$, a diagonal matrix $\diagon$ and a matrix ${\rm K}$ such that
$$
{\goodRepeta}\left(\partial^c_p\right)={\rm K}\diagon {\rm K}^{-1},\ \ddd{\goodRepeta}\left(c_p,\invol \right)={\rm K}\invol  {\rm K}^{-1}.
$$
Hence
\begin{eqnarray*}
\widetilde{\goodRepeta}\left(i_0\left(\partial_p^c\right)^{-1}\cdot   x_p\cdot    i_1\left(\partial_p^c\right)\cdot   x_p^{-1}\right)&=&
{\rm K}\diagon^{-1} {\rm K}^{-1}{\rm K}\invol {\rm K}^{-1} {\rm K}\diagon {\rm K}^{-1}{\rm K}\invol {\rm K}^{-1}\\
&=&{\rm K}\diagon^{-1} \invol  \diagon  \invol  {\rm K}^{-1}=\operatorname{Id}.
\end{eqnarray*}
Therefore $\widetilde{\goodRepeta}$ vanishes on the elements  generating  $H$ normally. Hence,  $\widetilde{\goodRepeta}$ gives rise to a homomorphism  $\widehat{\goodRepeta}$ from $\Gamma$ to $\sln$.

The homomorphism $\widehat{\goodRepeta}$ satisfies by construction the first condition to be a $\invol$-extension. 

Observe that by Equation (\ref{invdoteta}) and the construction of the generators of $\Gamma$, we have 
\begin{eqnarray}\widehat{\goodRepeta}\circ i&=&\invol {\goodRepeta} \invol \label{JRJ}.\end{eqnarray} Hence $\widehat{\goodRepeta}$ satisfies the second property (\ref{Jext2}) of the definition of $\invol$-extension.

Finally, we  check Property (\ref{Jext3}). Let  $c$ be a curve in $\Sigma$ joining $v$ to a boundary component $\partial_p$. We observe that
there exists a loop in $\Sigma$ based $\gamma$ at $v$ such that as elements of $\mathcal B_v$
$c=\gamma\comp c_p$,
thus
$$
c\comp {\bar c}^{-1}=\gamma\comp  c_p\comp  \bar c_p^{-1} \comp  \bar \gamma^{-1}.
$$
Hence, using Equation (\ref{JRJ}) and the definition of $\widehat{\goodRepeta}$
$$
\widehat{\goodRepeta} \left( c\comp {\bar c}^{-1} \right)=\widehat{\goodRepeta}\left(\gamma\right)\ddd{\goodRepeta}\left(c_p,\invol \right)\invol  \widehat{\goodRepeta}\left(\bar\gamma^{-1}\right)= {\goodRepeta}\left(\gamma\right)\ddd{\goodRepeta}\left(c_p,\invol \right){\goodRepeta}\left(\gamma^{-1}\right)\invol.
$$
It follows from Equation (\ref{ghg}) that
\begin{eqnarray}
\widehat{\goodRepeta} \left(c\comp {\bar c}^{-1}\right)&=&  \ddd{\goodRepeta}\left(c,\invol \right)\invol \label{cJc}.
\end{eqnarray}
Thus we have completed the proof of the existence and uniqueness of $\invol$-extension.

The continuity statement follows from the construction and the continuity statement in Proposition \ref{propGoodRepeta}.
 \qed

\subsubsection{Doubling Hitchin representations}

Let  again $\Sigma$ be a compact surface with non empty  boundary. Let  $\widehat \Sigma$ be its double.  Let $\invol$ be an involution of $\operatorname{PGL}(n,\mathbb R)$

We define extension for representations.

\begin{definition}
Let $\rho$ be a Hitchin representation from $\grf$ to $\sln$. A $\invol$-extension of $\rho$ is a representation $\widehat\rho$ from $\pi_1(\widehat\Sigma)$ to $\sln$ such that the $\invol$-extension of any good representative of $\rho$ belongs to the conjugacy class of $\widehat\rho$. 
\end{definition}

We now prove the main result of this section

\begin{theorem}\label{goodhomo}{\sc[Doubling]}
Let $\rho$ be a Hitchin representation from $\grf$ to $\sln$. Then there exists a unique $\invol$-extension $\widehat\rho$ of $\rho$.

Finally, the map $\rho\mapsto\widehat\rho$ is continuous.
\end{theorem}

\begin{definition}{\sc[Hitchin double]} Let  $$
\invol_n =\left(
\begin{array}{cccc}
1 & 0 & 0 &\ldots\\
0 &-1 & 0 &\ldots\\
0 & 0 & 1 &\ldots\\
\ldots & \ldots & \ldots &\ldots
\end{array}
\right).
$$
The $\invol_n$-extension of  a Hitchin representation $\rho$ is  the {\em Hitchin double} of $\rho$.
\end{definition}
The discussion extends to any real reductive split group. The following Corollary (proved in Paragraph \ref{hitchindoubleproof}) is important
\begin{coro}\label{hitchindouble}
Let $\Sigma$ be a compact connected oriented surface with  non empty boundary whose double $\widehat\Sigma$ has genus at least two. Let $\rho$ be a Hitchin representation of $\pi_1\left(\Sigma\right)$.  Then the Hitchin double of $\rho$ is a Hitchin representation of a closed surface. In particular, there exists a Hitchin representation $\widehat\rho$ of $\pi_1(\widehat\Sigma)$ whose restriction to $\pi_1\left(\Sigma\right)$ is $\rho$.
\end{coro}

\proof
We prove the first assertion of Theorem \ref{goodhomo}. The uniqueness of $\widehat\rho$ is the consequence of the uniqueness of $\invol$-extensions.

To prove the existence, we show  that the $\invol$-extension  $\widehat{\goodRepeta}$ of a good representative $\goodRepeta$ based at $v$ is independent -- up to conjugation -- of the choices made for ${\goodRepeta}$. We have two degrees of freedom in the construction of ${\goodRepeta}$
\begin{itemize}
\item the choice of a good representative based at $v$,
\item the choice of $v$.
\end{itemize}

The choice of a different  good representative based at $v$ amounts in conjugating ${\goodRepeta}$ by a diagonal matrix $\diagon$. It follows from Equation 
(\ref{ghd}) that $\diagon\widehat{\goodRepeta}\diagon^{-1}$ is the $\invol$-extension of $\diagon{\goodRepeta}\diagon^{-1}$. Thus the conjugacy class of $\widehat{\goodRepeta}$ is invariant under the choice of a good representative when the based point is fixed.

Secondly, we can choose another base point $w$ in a  boundary component $\partial_w$  of $\Sigma$. Let $\gamma$ be an arc from $v$ to $w$. Let ${\rm K}$ be a matrix  such that 
\begin{eqnarray}
{\goodRepeta}\left(\gamma\comp  \partial_w\comp  \gamma^{-1}\right)&=&{\rm K} \diagon_w {\rm K}^{-1},\\
\ddd{\goodRepeta}(\gamma,\Delta)&=&{\rm K}\Delta {\rm K}^{-1}\label{KDK}
\end{eqnarray}
where $\diagon_w$ is a diagonal matrix with decreasing entries.
It follows that the homomorphism  $\goodRepmu$ given by 
$$
\goodRepmu\left(h\right) = {\rm K}^{-1} {\goodRepeta}\left(\gamma\comp  h\comp  \gamma^{-1}\right) {\rm K},
$$
is a good representative based at $w$. 
We are going to prove that the homomorphism ${\bf F}$ defined by 
\begin{eqnarray}
{\bf F} \left(c\right) &=&{\rm K}^{-1}\widehat{\goodRepeta}\left(\gamma\comp   c\comp   \gamma^{-1}\right) {\rm K},\label{defF}
\end{eqnarray}
is the $\invol$-extension of $\goodRepmu$.
First, we note that if $c$ is an arc from $v$ to a component $\partial_u$, then
$$
\goodRepmu \left(\gamma^{-1}\comp  c\comp  \partial_u\comp   c^{-1}\comp  \gamma\right)={\rm K}^{-1} {\goodRepeta}\left(c \comp  \partial_u\comp   c^{-1}\right){\rm K}.
$$
It follows that
\begin{eqnarray}
\ddd\goodRepmu \left( \gamma^{-1}\comp c,  \diagon\right)&=&{\rm K}^{-1}\ddd{\goodRepeta}\left(c,\diagon\right) {\rm K}.\label{Mu}
\end{eqnarray}
We now prove that ${\bf F}$ is the $\invol$-extension of $\goodRepmu$. The first property of $\invol$-extension is obviously satisfied. For the second one, observe that  by definition.
\begin{eqnarray*}
{\bf F} \left({\bar c}\right) &=&{\rm K}^{-1} \widehat{\goodRepeta}\left(\gamma\comp   \bar c\comp   \gamma^{-1}\right)  {\rm K}\\
&=&{\rm K}^{-1} \widehat{\goodRepeta}\left({\gamma\comp   \bar\gamma^{-1}}\right) \widehat{\goodRepeta}\left(\overline{\gamma\comp  c\comp   \gamma^{-1}}\right) \widehat{\goodRepeta}\left(\overline{\gamma\comp   \bar\gamma^{-1}}\right) {\rm K}.
\end{eqnarray*}
It follows by the second and third properties of $\invol$-extensions  (\ref{Jext2}) and  (\ref{Jext3})  that
\begin{eqnarray*}
{\bf F} \left({\bar c}\right)
&=&{\rm K}^{-1}\widehat{\goodRepeta}\left(\gamma\comp   \bar \gamma^{-1}\right) \invol\widehat{\goodRepeta}\left({\gamma^{-1}\comp  c\comp   \gamma}\right)\invol \invol\widehat{\goodRepeta}\left(\gamma\comp   \bar \gamma^{-1}\right) \invol{\rm K}\\
&=&{\rm K}^{-1}\ddd{\goodRepeta}\left(\gamma,\invol \right)\widehat{\goodRepeta}\left(\gamma^{-1}\comp  c\comp  \gamma\right)\ddd{\goodRepeta}\left(\gamma,\invol \right) {\rm K}.
\end{eqnarray*}
Finally using (\ref{KDK}) for $\Delta=\invol$,
\begin{eqnarray*}
{\bf F} \left({\bar c}\right)
=\invol  {\rm K}^{-1}\widehat{\goodRepeta}\left(\gamma^{-1}\comp  c\comp  \gamma\right){\rm K} \invol 
=\invol {\bf F}\left(c\right)\invol.
\end{eqnarray*}
Thus the second property (\ref{Jext2}) of the definition of a $\invol$-extension is satisfied by $\bf F$.
For the last property (\ref{Jext3}) of the definition of $\invol$-extension, let  $c$ be a curve from $w$ to another boundary component, then by definition
\begin{eqnarray*}
{\bf F} \left(c\comp{\bar c}^{-1}\right)&=&{\rm K}^{-1} \widehat{\goodRepeta}\left(\gamma\comp c \comp  {\bar c}^{-1} \comp  \gamma^{-1}\right)  {\rm K}\\
&=&{\rm K}^{-1} \widehat{\goodRepeta}\left(\gamma\comp c\comp \overline{(\gamma\comp c)}^{-1}\right)
 \widehat{\goodRepeta}\left(\bar\gamma\comp \gamma^{-1}\right) {\rm K}.
\end{eqnarray*}
Hence, using (\ref{Jext2}),
\begin{eqnarray*}
{\bf F} \left(c\comp{\bar c}^{-1}\right)
&=&{\rm K}^{-1}  \ddd{\goodRepeta}\left( \gamma\comp c,\invol \right)\invol  
 \widehat{\goodRepeta}\left(\bar\gamma\comp \gamma^{-1}\right) {\rm K}\\
&=&{\rm K}^{-1} \ddd{\goodRepeta}\left(\gamma\comp c,\invol \right)   
 \widehat{\goodRepeta}\left(\gamma\comp \bar\gamma^{-1}\right)\invol {\rm K}\\
&=&{\rm K}^{-1}\ddd{\goodRepeta}\left(\gamma\comp c,\invol \right) \ddd{\goodRepeta}\left(\gamma,\invol \right){\rm K}\\
&=& {\rm K}^{-1}\ddd{\goodRepeta}\left(\gamma\comp c,\invol \right) {\rm K}\invol \end{eqnarray*}
It follows that by (\ref{Mu})
\begin{eqnarray*}
{\bf F} \left(c\comp {\bar c}^{-1}\right)
&=&\ddd{\bf M}\left(c,\invol \right)\invol .
\end{eqnarray*}
It follows that ${\bf F}$ is the $\invol$-extension of $\goodRepmu$. Hence $\widehat\rho$ is well defined as a conjugacy class. We have completed the proof of the existence and uniqueness of $\widehat\rho$.

The continuity statement follows from  the continuity statement in Proposition \ref{goodhomo0}.\qed
\vskip 1 truecm
\subsubsection{Proof of  Theorem \ref{hitchinbord} and Corollary \ref{hitchindouble}}\label{hitchindoubleproof}
 
A Fuchsian representation $\rho$  from $\grf$ to $\operatorname{PSL}\left(2,\mathbb R\right)$ is a monodromy representation of hyperbolic metric on  $\Sigma$ with totally geodesic boundary. The monodromy of the closed hyperbolic surface obtained by isometric gluing 
of $\Sigma$ along its boundary is the $\invol_2$-extension of $\rho$. Indeed $\invol_2$ is the symmetry with respect to a geodesic. Hence, 
the  Hitchin double  $\widehat\rho$ is the monodromy of the closed hyperbolic surface homeomorphic to the topological double 
of $\Sigma$ and,  in particular, it  is also Fuchsian.

We consider the irreducible representation of dimension $n$ given by the action of $\operatorname 
{GL}\left(2,\mathbb R\right)$ on homogeneous polynomials of two variables of degree $n-1$. Then the image of $\invol_2$  is $\invol_n$. It follows that the Hitchin double of an $n$-Fuchsian representation is an $n$-Fuchsian representation 
 and so belongs to the Hitchin component.
 Now recall that any   Hitchin representation $\rho_{1}$ of  $\grf$ 
 is obtained by continuous deformation of some Fuchsian representation $\rho_{0}$;
 that is there is a continuous path $\rho_t, t \in [0,1]$ of representations connecting $\rho_{0}$ to $\rho_{1}$.
 The ``doubling map" $\rho \mapsto \widehat{\rho}$ is continuous so that  $\widehat{\rho}_{t}$ 
 is a path of representations of the fundamental group of the topological double of $\Sigma$.
 This, in fact,  proves Corollary \ref{hitchindouble} as
  $\widehat{\rho}_{t}$ starts in the Hitchin component at $\widehat{\rho}_{0}$
  and so  remains in the Hitchin component for all $t\in [0,1]$;
  in particular $\widehat{\rho}_{1}$  is a Hitchin representation.

 Now, since the Hitchin double $\widehat{\rho}$ is a  Hitchin representation,
   Theorem \ref{mainA} applies and  
 we see that  $\widehat{\rho}$   preserves a Frenet curve, hence a positive curve by Lemma \ref{weylpos}. Hence, 
 the Hitchin double is a positive representation.
Finally, it follows directly from  the definition of positivity that 
 the restriction of  the 
Hitchin double $\widehat{\rho}$ to the 
fundamental group of any connected  incompressible surface  in the double surface is a positive representation.  

The statement involving the cross ratio is a consequence of the same result for closed surfaces -- see Theorem \ref{hcintro}.
This concludes the proof of  Theorem \ref{hitchinbord}. \qed

\section{Gap functions and coordinates for Hitchin rep\-re\-sen\-ta\-tions}

Let $\rho$ be a Hitchin representation of the fundamental group of  a surface $\Sigma$ and $P$ be a pair of pants in $\Sigma$. By Theorem \ref{hitchinbord}, $\partial_\infty\pi_1(\Sigma)$ inherits a cross ratio $\bb_\rho$. 

In order to complete our circle of ideas and fully generalise the results of Section \ref{hypcase}, we wish to describe the gap functions for $P$. 
Recall that by Proposition \ref{gapsubsurface}, these gap functions only depend on $\rho(\pi_1(P))$.  The aim of this section is to compute these functions using coordinates -- generalising shear coordinates -- on the space of Hitchin representations of $\pi_1(P)$ which come from coordinates on the Fock--Goncharov moduli space ${\mathcal M}_{FG}(P )$ -- see below for definitions.

Therefore,  the aim of this section -- parallel to Section \ref{hypcase} --  is threefold
\begin{itemize}
\item  Following \cite{Fock:2006a}, we describe  the {\em Fock--Goncharov moduli space} of the pair of pants and its {\em Fock-Goncharov coordinates} in the next paragraph. This Fock-Goncharov moduli space is mapped onto the space of positive representations by a {\em holonomy map}.
\item In Paragraph \ref{lifthitch}, we describe the preimage  in the Fock--Goncharov moduli space of a Hitchin representation under the holonomy map. This construction induces various {\em Fock--Goncharov coordinates for the moduli space of Hitchin representation} of the fundamental group of a pair of pants Paragraph \ref{coordpant}.
\item  We finally express in Theorem \ref{theo:goodgap} the pant gap function for a good choice of coordinates on the space of Hitchin representations.
\end{itemize}

\subsection{Fock--Goncharov moduli space for the pair of pants}

We start with the canonical ideal triangulation ${\mathcal T}$ of  the 3-punctured sphere $P$ obtained by gluing two triangles $X$ and $Z$. 
This ideal triangulation has two faces $X$ and $Z$, three ideal  vertices $\alpha$, $\beta$ and $\gamma$, and three edges $A$, $B$  and $C$
  where $A$ is the edge opposite  the vertex $\alpha$ {\em etc}.
\begin{definition}{\sc[Fock--Goncharov Moduli space]}\label{def:fgm}
Following  \cite{Fock:2006a}, an element of the {\em Fock--Goncharov moduli space} ${\mathcal M}_{FG}(P )$ is 
a  configuration of six flags 
$$
\mathcal S=(X_\alpha,Z_\gamma, X_\beta,  Z_\alpha, X_\gamma, Z_\beta)\in {\mathcal F}(\mathbb R^n)^6/{\rm PSL}(n,\mathbb R),$$ 
such that 
\begin{enumerate} 
\item The triple  $T_X=(X_\alpha, X_\beta, X_\gamma)$ is positive.
\item  The three  triples $(Z_\alpha, X_\gamma, X_\beta)$, $(X_\alpha, X_\gamma,Z_\beta)$ and  $(X_\alpha, Z_\gamma,X_\beta)$ are positive and  are equivalent under the action of ${\rm PSL}(n,\mathbb R)$ to some positive triple $T_{Z}$.\label{equiv-fgm}
\item  The three quadruples
 \begin{eqnarray*}
Q_A&=&(X_\gamma, Z_\alpha,X_\beta,X_\alpha),\\
Q_B&=&(X_\alpha,Z_\beta,X_\gamma,X_\beta),\\
Q_C&=&(X_\beta,Z_\gamma,X_\alpha,X_\gamma),
\end{eqnarray*}
 are positive.
 \end{enumerate}
 \end{definition}

We associate  to an element of the Fock--Goncharov moduli space a representation of the fundamental group of the pair of pants in $\sln$ as follows. 

We    
 choose a complete finite volume hyperbolic metric on ${P }$ such that the edges or $\mathcal T$ are geodesics. Let $\widetilde{P}$ be the universal cover of ${P}$.  
 Let  $\widehat{\mathcal T}$ be the pull back of the ideal triangulation ${\mathcal T}$ on $\widetilde{P }$. 
 Then let $\widehat{\mathcal V}$ 
 denote the set of ideal vertices   of $\widehat{\mathcal T}$. We consider $\widehat{\mathcal V}$ as a subset of $\partial_\infty{\pi_1({P })_3}$, and observe that  this identification is independent of  the choice of the hyperbolic metric.
 
Note that the fundamental group $\pi_1({P })$ acts on $\widehat{\mathcal V}$.  We describe this action as follows. For every lift  $X^i$  of our original triangle $X$ -- respectively $Z^i$ of $Z$ --  we denote by $(X^i_\alpha,X^i_\beta,X^i_\gamma)$  -- respectively $(Z^i_\alpha,Z^i_\beta,Z^i_\gamma)$ -- the lift of the vertexes $(\alpha,\beta,\gamma)$, and follow a similar convention for edges.

\begin{figure}[here]
   \centering
   \includegraphics[width=3in]{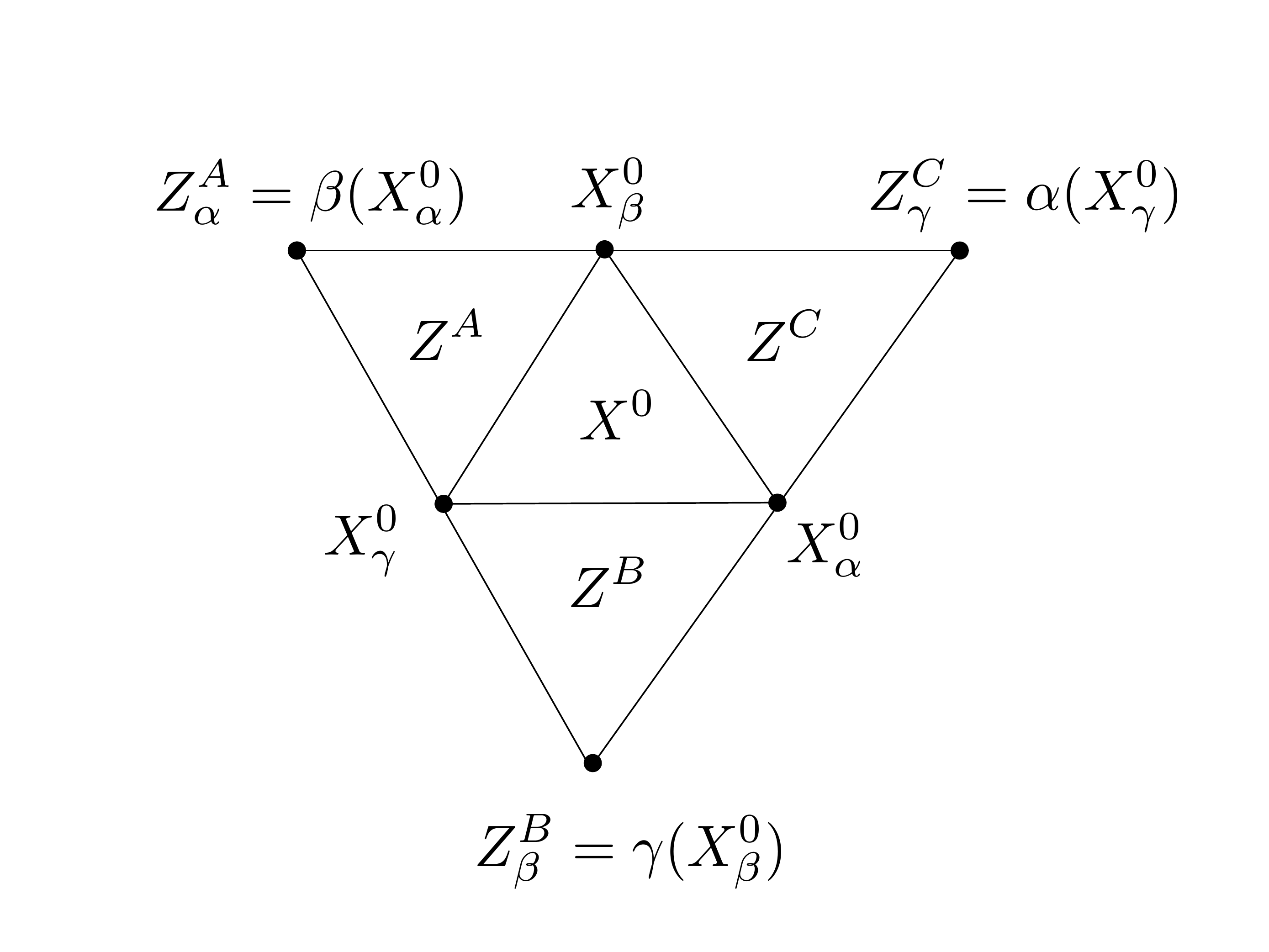} 
   \caption{\em The action of the fundamental group
   \protect\label{fig:actionpi}}
\end{figure}

We now choose a particular lift $X^0$ of $X$. Let then $Z^A$ be the lift of $Z$ adjacent to $X^0_A$ and similarly for the other edges $B$ and $C$.
We then denote  -- again abusively --
  by $\alpha$ the element of $\pi_1(P)$ that fixes the vertex $X^0_\alpha$ and sends $X^0_\gamma$ to  
$Z^C_\gamma$, and symmetrically for $\beta$ and $\gamma$ as in Figure \ref{fig:actionpi}.

 According to Theorem 7.1 of  \cite{Fock:2006a} -- which is another incarnation of the general principle  associating to a geometric structure a holonomy representation -- an element of the Fock--Goncharov moduli space defines a pair $(f,\rho)$ where  $\rho$ is  a homomorphism from $\pi_1( {P })$ to ${\rm PSL}(n,\mathbb R)$ and $f$ is a $\rho$-equivariant map from $\widehat{\mathcal V}$ to $\mathcal F(\mathbb R^n)$
-- well defined up to congugacy by ${\rm PSL}(n,\mathbb R)$. The pair $(f,\rho)$ is characterised by  $\mathcal S$ in the following way
\begin{eqnarray*}
X_\alpha=f(X^0_\alpha),& \ \ &Z_\gamma=\rho(\beta)\left(f(X^0_\alpha)\right),\\
X_\beta=f(X^0_\beta),&\ \ &Z_\alpha=\rho(\gamma)\left(f(X^0_\beta)\right),\\
X_\gamma=f(X^0_\gamma),&\ \ &Z_\beta=\rho(\alpha)\left(f(X^0_\gamma)\right).
\end{eqnarray*}
The existence of $\rho(\alpha)$, $\rho(\beta)$ and $\rho(\gamma)$,  as well as the fact that $\rho$ is indeed a homomorphism,  is guaranteed by Condition \ref{equiv-fgm}  of Definition \ref{def:fgm}.

Then, by construction,  if $(U,V,W)$ and $(Y,U,W)$ are  pairs of adjacent triangles  in $\widehat{\mathcal T}$ with $Z\neq T$,  then
 $\left(f(U),f(V),f(W),f(Y)\right)$ is a positive quadruple.

\begin{definition} The representation $\rho$  -- denoted by  $\hol (\mathcal S)$ --
described above is the {\em holonomy} of the element $\mathcal  S$ of the Fock--Goncharov moduli space. 
\end{definition}

This representation is positive and, moreover, it follows from Proposition 5.8  and Theorem 7.1 in \cite{Fock:2006a} that all positive representations are obtained this way.

\subsubsection{Fock--Goncharov coordinates for a pair of pants} \label{coordpant} 
Following \cite{Fock:2006a} and by Proposition \ref{FGquad}, a parametrisation of ${\mathcal M}_{FG}({P })$ is given   by the following collections of functions of $\mathcal S$  
$$(X^{m,l,p},Z^{m,l,p},\Delta^A_k,\Delta^B_k,\Delta^C_k),$$
 for  $1\leq k\leq n-1$ and  $m+l+p=n$, which we  define as follows using the notations of Section \ref{sec:posflag}
\begin{itemize}
\item  triple ratios describing $T_X$ : $X^{m,l,p}(\mathcal S)=T^{m,l,p}(T_X)$ for $m+l+p=n$,
\item  triple ratios describing $T_Z$ : $Z^{m,l,p}(\mathcal S)=T^{m,l,p}(T_Z)$ for $m+l+p=n$,
\item  edges functions \begin{eqnarray*}
\Delta^A_k(\mathcal S)&=&\delta_k(Q_A),\\
\Delta^B_k(\mathcal S)&=&\delta_k(Q_B),\\
\Delta^C_k(\mathcal S)&=&\delta_k(Q_C),
\end{eqnarray*} for $1\leq k\leq n-1$.
\end{itemize}

\subsection{Hitchin representations for the pair of pants}
\subsubsection{Lifts of Hitchin representations}\label{lifthitch}

 Let ${\operatorname{Hom}}_H({P })$ be the  space of Hitchin homomorphism of $\pi_1(P)$ into $\sln$. The space ${\operatorname{Hom}}_H({P })$ have one or two connected components depending on the parity of $n$, and the group of outer automorphisms of $\sln$ acts transitively by conjugation  on it set of connected component . Let
 $$
 \operatorname{Rep}_H(P)=\operatorname{Hom}_H(P)/\operatorname{PGL}(n,\mathbb R),
 $$ 
 be the  {\em space of Hitchin representations} of $\pi_1(P)$ into $\sln$.  Let ${\mathcal M}_{FG}({P })$ be the Fock--Goncharov moduli space. Let $(u,v,w)$ be a triple of elements of the Weyl group $\mathcal W$ of $\sln$. 
 
Observe first that given a purely loxodromic element ${\rm M}$, we have a well defined action of the Weyl group -- identified  with the symmetric group of $\{1,\ldots,n\}$  -- on the set of invariant flags by ${\rm M}$ characterised as follows:  let $(L_1,L_2,\ldots, L_n)$ be the tuples of eigenlines ordered by decreasing eigenvalues;  let $u$ be an element of the Weyl group; then the image of the attracting flag ${\rm M}^+$ is 
$$
{\rm M}^u=(L_{u(1)},L_{u(1)}\oplus L_{u(2)},\ldots).
$$
\begin{proposition} Let $(u,v,w)$ be a triple of elements of the Weyl group.
There exists a map
$$ 
{\mathcal S}^{u,v,w}:
\mapping
{{\operatorname{Rep}}_H({P })}
{{\mathcal M}_{FG}({P }),}
{\rho}
{{\mathcal S}^{u,v,w}_\rho,}$$
 such that ${\hol}\circ {\mathcal S}^{u,v,w}={\rm Id}$. 
This map is onto a connected component of ${\hol}^{-1}({\operatorname{Rep}}_H({P }))$.
\end{proposition}

\proof
 Recall first  that by Corollary \ref{hitchindouble}, the double $\widehat\rho$ preserves a Frenet curve, which we denote by $\xi$.

Let  now $\alpha$ be a boundary component of ${P }$. Recall that since $\rho(\alpha)$ is purely loxodromic, its dynamics on the flag manifolds has the following property : the action of $\rho(\alpha)$ on the flag manifold has exactly one attractive point that we call the {\em attractive flag} and denote by $A^+$. The attractive flag is $A^+=(L_1,L_1\oplus L_2, \ldots)$ where $L_i$ is the eigenspace for the eigenvalue $\lambda_i$, where $\lambda_1>\lambda_2>\ldots$. Finally it follows from \cite{Labourie:2006} that the attractive flag  $A^+$ is the osculating flag of the curve $\xi$ at the point $\xi(\alpha^+)$.

Symmetrically the action of $\rho(\alpha)$ has exactly one repulsive point  that we call the {\em repulsive flag}  and denote by $A^-$.

Let  $\mathcal W$ be the Weyl group of $\sln$. Then for every element $u$ in $\mathcal W$ in $\sln$, we denote a in the introduction of this section
$$
A^u=(L_{u(1)},L_{u(1)}\oplus L_{u(2)},\ldots),
$$
the image of $A^+$ by $u$. Similarly we have flags $B^v$ and $C^w$ for the boundary components associated to $\beta$ and $\gamma$ and the elements $v$ and $w$ of the Weyl group.

It follows from Section 7 and 9 of \cite{Fock:2006a}, 
 that the configuration given by the following sextuplet is positive 
$$
\mathcal S_\rho^{u,v,w}= (A^u,\rho(\beta)(A^u),
B^v, \rho(\gamma)(B^v),C^w, \rho(\alpha)(C^w)).
$$
Hence $\mathcal S^{u,v,w}_\rho$ determines 
 an element of the Fock--Goncharov moduli space for the pair of pants whose image under  ${\hol}$ is $\rho$.
 
 Conversely it follows from the description in the previous paragraph that for any ${\mathcal S}$ in ${\mathcal M}_{FG}({P })$  such that 
 $$ 
 {\hol}({\mathcal S})=\rho,
 $$
 there exists a unique triple  $(u,v,w)$ in the Weyl group such that 
 $$
 {\mathcal S}=\mathcal S^{u,v,w}_\rho.
 $$
Indeed any flag invariant by $\rho(\alpha)$ is of the form $A^u$ for some $u$ in the Weyl group. 
Thus, the map $(\rho,u,v,w)\mapsto \mathcal S_\rho^{u,v,w}$ from ${\operatorname{Rep}}_H(P)\times \mathcal W^3$ to ${\hol}^{-1}({\operatorname{Rep}}_H({P }))$ is a homeomorphism. Since ${\operatorname{Rep}}_H(P)$ is connected it follows that  $S^{u,v,w}$ is onto a connected component of ${\hol}^{-1}({\operatorname{Rep}}_H({P }))$.\qed

\begin{figure}[here]
   \centering
   \includegraphics[width=3in]{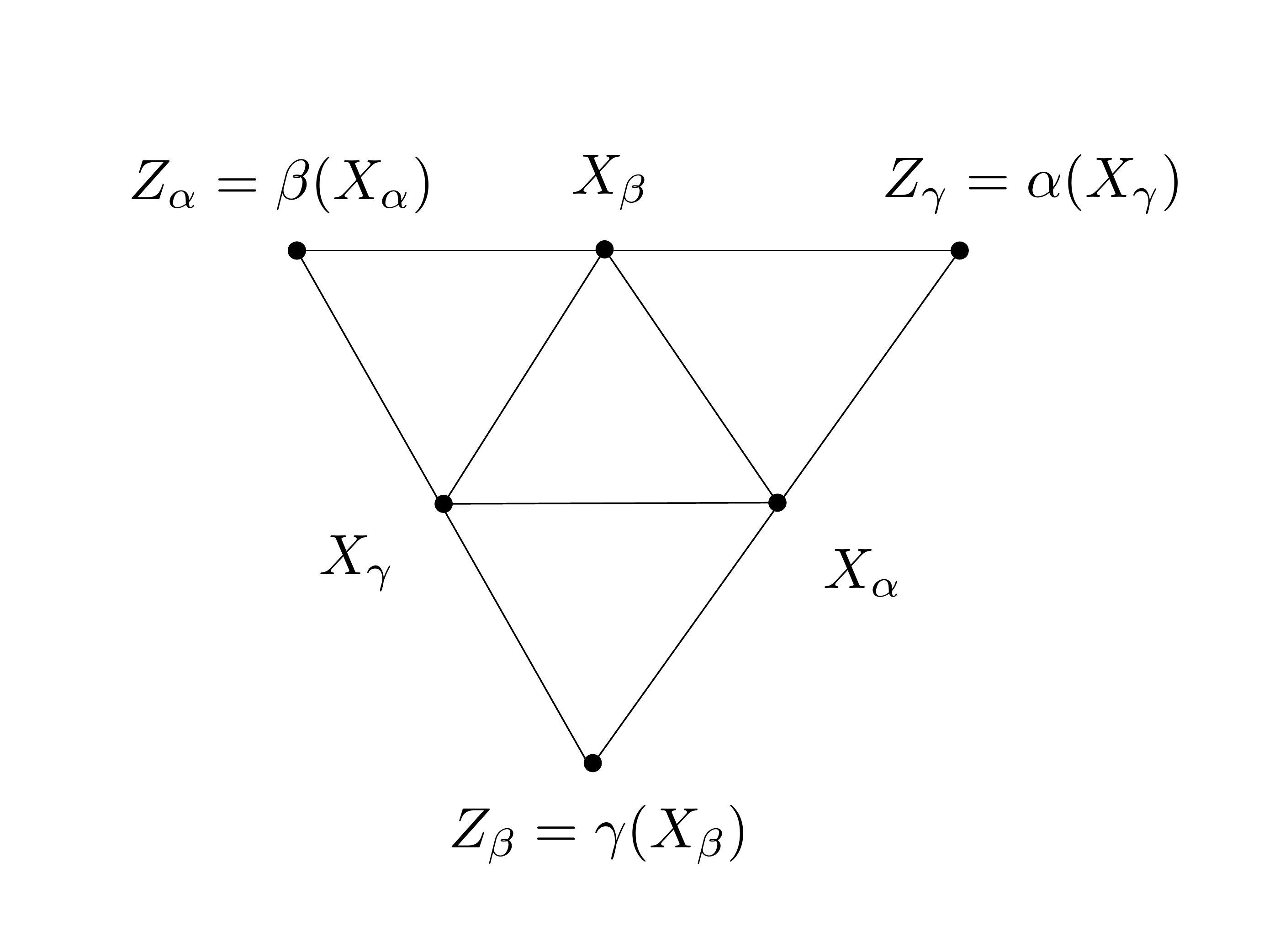} 
   \caption{\em The configuration of  flags
   \protect\label{fig:xz}}
\end{figure}

We describe the configuration of flags  ${\mathcal S}^{u,v,w}_\rho$  in  Figure \ref{fig:xz} analoguously to  Figure \ref{fig:actionpi},   where we identify  the fundamental group with its image under $\rho$,  and use the following notation
\begin{eqnarray*}
Z_\alpha&=&\beta(X_\alpha)=\gamma^{-1}(X_\alpha),\\
Z_\beta&=&\gamma(X_\beta)=\alpha^{-1}(X_\beta),\\
Z_\gamma&=&\alpha(X_\gamma)=\beta^{-1}(X_\gamma),
\end{eqnarray*}
where  $X_\alpha=A^u$, $X_\beta=B^v$ and $X_\gamma=C^w$ are so that 
$$
\mathcal S^{u,v,w}_\rho=(X_\alpha,Z_\gamma, X_\beta,  Z_\alpha, X_\gamma, Z_\beta),$$
is positive. 

In the case $n=2$, the Fock--Goncharov moduli space is nothing 
but the enhanced Teichmüller space of the pair of pants
 \cite{Bonahon:1996}.
 Recall that  a point of the enhanced  Teichmüller space  of a surface with boundary consists
 of a  point of the usual Teichmüller space plus a choice of orientation for 
 each of the boundary components. The discussion of the previous paragraph is exactly a generalisation of this discussion.

\subsubsection{Coordinates on the moduli space of Hitchin representations}\label{coordhitch}

Therefore using the notations of Paragraph \ref{coordpant} and \ref{lifthitch}, we  obtain  a system of coordinates on Hitchin representations
by  choosing a triple $(u,v,w)$ of elements of the Weyl group -- one for each  boundary component --
and considering the functions
\begin{eqnarray*}
X^{m,l,p}_{u,v,w}:\rho&\mapsto& X^{m,l,p}_{u,v,w}(\rho)=X^{m,l,p}(\mathcal S_\rho^{u,v,w}),\\
Y^{m,l,p}_{u,v,w}:\rho&\mapsto& Y^{m,l,p}_{u,v,w}(\rho)=Y^{m,l,p}(\mathcal S_\rho^{u,v,w}),\\
A^{u,v,w}_k:\rho&\mapsto& A^{u,v,w}_k(\rho)=\Delta^A_k(\mathcal S_\rho^{u,v,w}),\\
B^{u,v,w}_k:\rho&\mapsto& B_k^{u,v,w}(\rho)=\Delta^B_k(\mathcal S_\rho^{u,v,w}),\\
C^{u,v,w}_k:\rho&\mapsto& C^{u,v,w}_k(\rho)=\Delta^C_k(\mathcal S_\rho^{u,v,w}).
\end{eqnarray*}

We call this set of functions the {\em Fock--Goncharov coordinates on the moduli space of Hitchin representations} associated to the triple $(u,v,w)$.
This is an extension of Thurston shear coordinates described in 
Section \ref{shear} for Fuchsian representations when $n=2$. 
Observe that this construction leads to  $(n!)^3$ different  coordinate systems on the space of Hitchin representations. To the knowledge of the authors, there is no nice formula and explicit for these changes of coordinates.
\subsection{Gap functions, periods  and coordinates}
In general, it seems difficult to obtain a closed form for the gap functions 
for each of these coordinate systems.
However, for  particular  choices of coordinates, we have a nice formula for the pant gap function.
 
\begin{theorem}\label{theo:goodgap}
Let  $u$, $v$ and $w$ be elements of the Weyl group of $\sln$,
identified with the group of permutations of $\{1,\ldots,n\}$, one for each boundary component. We assume  that 
$$
u(n)=n, \ u(n-1)=1,\  v(1)=1,\  w(1)=n.
$$
Let $\rho$-be an $n$-Hitchin  representation $\pi_1(P)$ into $\sln$. Let $\xi$ be Frenet curve with values in $\mathbb P(\mathbb R^n)$ associated to $\rho$ and $\xi^*$ its hyperplane osculating curve. Let $\bb=\bb_{\xi^*,\xi}$ be the cross ratio associated to the pair $(\xi^*,\xi)$. Then the pant gap function of the cross ratio $\bb$  has the following expression in the Fock--Goncharov coordinates associated to the triple $(u,v,w)$
\begin{eqnarray*}
G_\bb(P)&=&\log\Big(\frac{1+C_1(\rho)e^{\ell_\bb(\alpha)}}{1+C_1(\rho)}\Big)=\log\Big(\frac{1+B_{n-1}(\rho)}{B_{n-1}(\rho)(1+C_1(\rho))}\Big).
\end{eqnarray*}
(We have used the notation introduced in the previous paragraph, but ignored the superscripts $u,v,w$ in order to obtain a readable formula.)
\end{theorem}
\rmks
\begin{enumerate}
\item In other words, we choose invariant flags $A$, $B$ and $C$ so that
\begin{eqnarray*}
A^{n-2}&=&\xi_{n-1}(\alpha^+)\cap\xi_{n-1}(\alpha^-),\\  A^{n-1}&=&\xi_{n-1}(\alpha^+),\\
B^1&=&\xi_1(\beta^+),\\
C^1&=&\xi_1(\gamma^-).
\end{eqnarray*}
where $\xi$ is the limit curve from $\partial_\infty\tilde{P }$ to  $\RPN$, where $\tilde P$ is the universal cover of $P$ equipped with an admissible hyperbolic metric.
\item We shall explain in Section \ref{Xmap} how to obtain, for ${\rm PSL}(3,\mathbb R)$, and the choice of the identity for the elements of the Weyl group, the formula for the pant gap function using a computer assisted proof.
\item So far, it remains a challenge to obtain a closed formula for gap functions in all coordinates. The same remark holds for the boundary gap function, for which we do not have a similar nice formula.  However, one could prove -- but we shall not do it here -- that the gap functions are logarithms of rational functions.

\end{enumerate}
\proof Let $P$ be a pair of pants represented  by a triple $(\alpha,\beta,\gamma)$ of elements of $\grf$. Let us consider the three flags
$$
X_\gamma=C, \ \ X_\alpha=A, \ \ X_\beta=B.
$$
We complete  this configuration by
$$
Z_\gamma=\alpha(C), \ \ Z_\beta=\gamma(C), \ \ Z_\alpha=\beta(A).
$$
Now let $(u,v,w)$ be three elements of the Weyl group as described in the hypothesis of the theorem. 
Let 
$$
V={\mathbb R^n}/\xi_{n-1}(\alpha^+)\cap\xi_{n-1}(\alpha^-)={\mathbb R^n}/A^{n-2}.
$$
To  simplify the notation, we will write $\alpha$ for $\rho(\alpha)$.
We observe that $\alpha$ acts on  $V$. 
Let ${\bf b}_V$ be the cross ratio on  the projective line ${\mathbb P }(V)$. 
Let $\pi$ be the projection on $V$. By definition, the gap function is 
\begin{eqnarray*}
G_\bb(P)&=&\log(\bb_{\xi^*,\xi}(\alpha^+,\gamma^-,\alpha^-,\beta^+)), \\
\end{eqnarray*}
where
\begin{eqnarray*}
\bb_{\xi^*,\xi}(\alpha^+,\gamma^-,\alpha^-,\beta^+)
&=&\frac{\Omega(\xi_1(\gamma^-)\wedge \xi_{n-1}(\alpha^+))\Omega( \xi_1(\beta^+)\wedge \xi_{n-1}(\alpha^-))}{\Omega( \xi_1(\gamma^-)\wedge \xi_{n-1}(\alpha^-))\Omega( \xi_1(\beta^+)\wedge \xi_{n-1}(\alpha^+))}.
\end{eqnarray*}
Let 
\begin{eqnarray*}
c=\pi(\xi_1(\gamma^-)),& &
b=\pi(\xi_1(\beta^+)),\\
a^+=\pi(\xi_{n-1}(\alpha^+)),& &
a^-=\pi(\xi_{n-1}(\alpha^-)).
\end{eqnarray*}
By Equation (\ref{crosspi}), we have
\begin{eqnarray*}
\bb_{\xi^*,\xi}(\alpha^+,\gamma^-,\alpha^-,\beta^+)
&=&{\bf b}_V(a^+,c,a^-,b)
\end{eqnarray*}
Then, using the  that ${\bf b}_V$ is a cross ratio on a projective line and  following the proof  of Equation (\ref{cr3}) in the ${\rm PSL}(2,\mathbb R)$ case,  we obtain
\begin{eqnarray}
\bb_{\xi^*,\xi}(\alpha^+,\gamma^-,\alpha^-,\beta^+)
&=&\frac{1-{\bf b}_V(b, c,a^+,\alpha(c)) {\bf b}_V(a^+,c,a^-,\alpha(c))}{1-{\bf b}_V( b,c,a^+,\alpha(c)) }.\label{pantgb}
\end{eqnarray}
By Equation (\ref{deltacross}) and since $\alpha(c)=\pi(\alpha(\xi_1(\gamma^-)))$
\begin{eqnarray*}
{\bf b}_V(b, c,a^+,\alpha(c))&=&-\delta_{1}(B,\alpha(C),A,C)\\&=&-\Delta^C_1(\mathcal S_\rho)\\&=&-C_1(\rho).
\end{eqnarray*}
Finally, we observe that the two eigenvalues of $\alpha$ on $V$ are the largest and smallest eigenvalues of $\alpha$ on $\mathbb R^n$. Hence
$$
{\bf b}_V(a^+,c,a^-,\alpha(c))=e^{\ell_\bb(\alpha)}.
$$
This concludes the proof of the first equality. For the second equality,
we observe first using  Equation (\ref{l=bc}) that
$$
{\bf b}_V(a^+,c,a^-,\alpha(c)).{\bf b}_V(a^+,b,c,\alpha^{-1}(b)).{\bf b}_V(b,c,a^+,\alpha(c))=1.
$$
Multiplying the numerator and denominator of the left hand side of Equation (\ref{pantgb}) by ${\bf b}_V(a^+,c,a^-,\alpha(c))$, we get that
\begin{eqnarray}
\bb_{\xi^*,\xi}(\alpha^+,\gamma^-,\alpha^-,\beta^+)
&=&\frac{{\bf b}_V(a^+,c,a^-,\alpha(c))-1}{{\bf b}_V(a^+,c,a^-,\alpha(c))\left(1-{\bf b}_V( b,c,a^+,\alpha(c))\right)}
\end{eqnarray}
Finally, by Equation (\ref{deltacross})
\begin{eqnarray*}
{\bf b}_V(a^+,b,c,\alpha^{-1}(b))
&=&-\delta_{n-1}(A,\alpha^{-1}(B),C,B)\\
&=&-B_{n-1}(\rho).
\end{eqnarray*}
The second formula follows.
\qed

\section{Appendix A: the three dimensional case}\label{Xmap}

We use in this section the notations and results of Paragraph 5.2 of \cite{Fock:2007} which describes the monodromy associated to the  Fock--Goncharov coordinates in the case of ${\rm PSL}(3,\mathbb R)$ as that of a local system on a graph.

Our aim is to explain how one computes explicitly the pant gap function for various choices of elements in the Weyl group for the boundary components. Using the notations of Paragraph \ref{coordhitch} and in order to simplify our notations in this particular case, we set
\begin{eqnarray*}
x=X^{u,v,w}_{u,v,w}(\rho),& &y=Y^{u,v,w}_{u,v,w}(\rho),\\
a=A_1^{u,v,w} (\rho),& &A=A_2^{u,v,w} (\rho),\\
b=B_1^{u,v,w} (\rho),& &B=B_2^{u,v,w} (\rho),\\ 
c=C_1^{u,v,w} (\rho),& &C=C_2^{u,v,w} (\rho).
\end{eqnarray*}
We denote by
$$
M(x,y,a,A,b,B,c,C),
$$
The element of the Fock--Goncharov moduli space ${\mathcal M}_{FG}(P )$ whose coordinates are  $x,y,a,A,b,B,c,C$, and by $$
R(x,y,a,A,b,B,c,C),
$$
its holonomy representation. By construction
$$
{\mathcal S}^{u,v,w} (R(x,y,a,A,b,B,c,C))=M(x,y,a,A,b,B,c,C).
$$

We will actually explicit the gap functions for three different choices of coordinates which give very different formulae: Equations (\ref{comp:1}), (\ref{comp:2}) and (\ref{comp:3}). On the arXiv version of this article, we give the computer  assisted proof  of this result and explain the instructions to obtain the 18 different formulae.   

We consider as usual a pair of pants whose boundary are  $\alpha$, $\beta$, $\gamma$, so that the corresponding elements in $\pi_1(P )$ satisfy
$$
\alpha\gamma\beta=1.
$$
We consider positive representations as holonomies of discrete connections on a graph.
\subsection{The construction}
\subsubsection{Some matrices}

We consider the following matrices 
where we adopt the notation of Fock and Goncharov \cite{Fock:2007}
\begin{eqnarray*}
X=T(x),\ \ Y=T(y),& &
Q_a=E(A,a),\ \ Q_b=E(B,b),\ \ Q_c=E(C,c),
\end{eqnarray*}
where 
$$
T(x)=
\left(
\begin{array}{ccc}
0&0&1\\
0&-1&-1
\\ x&1+x&1
\end{array}
\right)
,\ \ E(w,z)=
\left(
\begin{array}{ccc}
0&0&\frac{1}{z}\\
0&-1&0
\\ w&0&0
\end{array}
\right).
$$
We observe that $T(x)^3=x\cdot{\rm Id}$.
\subsubsection{A discrete connection on a graph}
Following \cite{Fock:2007}, we label the edges of a graph by the above matrices following  Figure \ref{fig:graph}.
\begin{figure}
   \centering

   \includegraphics[width=2in]{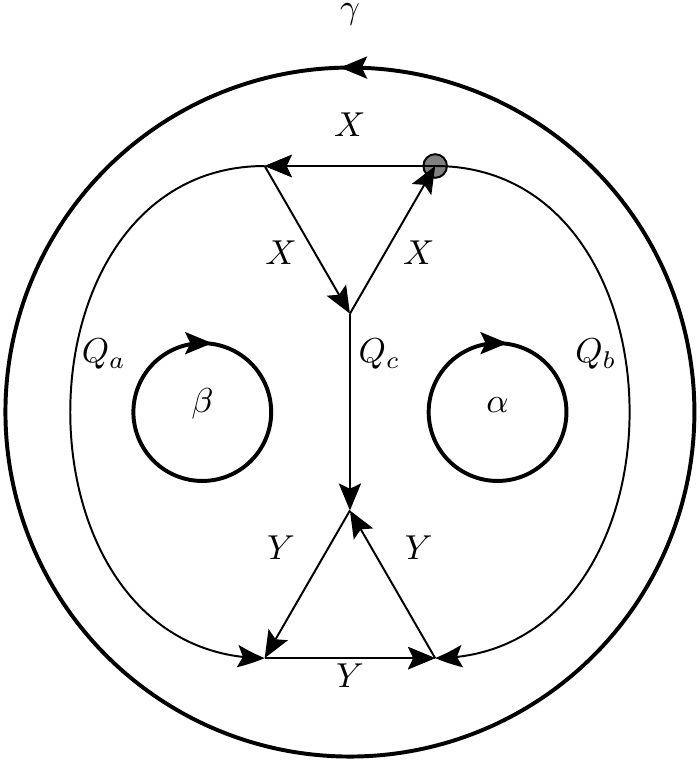} 
 \caption{\em Local system on a graph 
  \protect\label{fig:graph}}
\end{figure}
This way we obtain a discrete connection on this graph and we can compute the holonomies of the boundary components using $(x,y,A,a,B,b,C,c)$. 
Then, the holonomies (in $PGL(3,\mathbb R)$) corresponding to the boundary components $\alpha$, $\beta$, $\gamma$ starting at the point $M$ are respectively
\begin{eqnarray}
\rho(\alpha)=\Hol (A)&=& X\cdot Q_c^{-1}\cdot Y\cdot Q_b,\\
\rho(\beta)=\Hol (B)&=& X\cdot X\cdot Q_a^{-1}\cdot Y\cdot Q_b\cdot X\cdot X,\\
\rho(\gamma)=\Hol (C)&=& Q_b^{-1}\cdot Y\cdot Q_a\cdot X.
\end{eqnarray} 
A computer assisted calculation 
shows that $\Hol (A)$ is a lower triangular matrix, and that
 $\Hol (C)$ is upper triangular.
\subsubsection{Flags}
We now make the link with the description of the coordinates  as in Section \ref{coordpant}. 
First we introduce the following notation: if $(U,V)$ is a pair of independent vectors, 
we denote by $F(U,V)$ the flag $(\mathbb RU, \mathbb RU+\mathbb RV)$.
We now observe that the flags
\begin{eqnarray*}
X_\alpha&=&F((0,0,1),(0,1,0)),\\
X_\gamma&=&F((1,0,0),(0,1,0)).
\end{eqnarray*}
are invariant by $\Hol (A)$ and $\Hol (C)$ respectively, 
which are  respectively lower and upper triangular matrices. 
Furthermore
$$
X_\beta=X.X_\alpha=F((1,-1,1),(0,-1,1+x)),
$$
is invariant by $\Hol (C)$,  which is the conjugate by $X$ of a lower triangular matrix.

It is useful at this stage to calculate the dual invariant  flags. We have
\begin{eqnarray*}
X^*_\alpha&=&F((1,0,0),(0,1,0)),\\
X^*_\beta&=& F((0,0,1),(0,1,0)),\\
X^*_\gamma&=&F((x,1+x,1),(0,1+x,1)).
\end{eqnarray*}

The full configuration of flags as in Section \ref{coordpant} is now completed by 
\begin{eqnarray*}
Z_\alpha&=&\Hol (B)(X_\alpha)=\Hol (C)^{-1}(X_\alpha),\\
Z_\beta&=&\Hol (C)(X_\beta)=\Hol (A)^{-1}(X_\beta),\\
Z_\gamma&=&\Hol (A) (X_\gamma)=\Hol (B)^{-1}(X_\gamma). 
\end{eqnarray*}
Then, by the result of Fock and Goncharov in \cite{Fock:2007}, we have 
\begin{proposition} The configuration of flags $$
\mathcal S=(X_\alpha,Z_\gamma, X_\beta,  Z_\alpha, X_\gamma, Z_\beta),$$
is 
$M(x,y,A,a,B,b,C,c)$.
\end{proposition}
\subsubsection{Eigenvectors and flags}
We use  a computer program to compute the eigenvectors and eigenvalues for $(\Hol (A))^*$, $\Hol (B)$ and $\Hol (C)$. 
\begin{itemize}
\item We denote by $A_1,A_2,A_3$ the eigenvectors of $\Hol (A)^*$ of eigenvalues $\frac{xy}{bC},1, Bc$ respectively.
\item We denote by $B_1,B_2,B_3$ the eigenvectors of $\Hol (B)$ of eigenvalues $\frac{x^2y}{Ac},x, xCa$ respectively.
\item We denote by $C_1,C_2,C_3$ the eigenvectors of $\Hol (C)$ of eigenvalues $bC,1, \frac{xy}{Bc}$ respectively.
\end{itemize}
The results are given in the next section.
We observe that
\begin{eqnarray*}
X_\beta&=&F(B_1,B_2),\\
X_\gamma&=&F(C_3,C_2),\\
X^*_\alpha&=&F(A_3,A_2).
\end{eqnarray*}

\subsubsection{Coordinates for Hitchin representations and inequalities}

We have explained in Paragraph \ref{lifthitch} that various choices of elements $(u,v,w)$ of the Weyl group gives rise to different lifts ${\mathcal S}^{u,v,w}$ from ${\operatorname{Rep}}_H(P)$ to ${\mathcal M}_{FG}(P)$.  

We explain now how to describe the image of $S^{u,v,w}$ in the coordinates $(x,y,a,A,b,B,c,C)$. The answer follows from considering the converse question.
\vskip 0.2 truecm
{\em Given $(x,y,a,A,b,B,c,C)$ so that the holonomy $\rho$ is Hitchin, what is the triple  $(u,v,w)$ of elements of the Weyl group so that $$
S^{u,v,w}_\rho=M(x,y,a,A,b,B,c,C)\ \  ?$$}
\vskip 0.2 truecm
To answer this question, we have to identify the elements of  the Weyl group  which send the attractive flags of the corresponding boundary components to  the invariant flags  $X_\alpha$, $X_\beta$ and $X_\beta$. More precisely, the Weyl group acts on the set of invariant flags by the monodromy of a boundary component -- see the beginning Paragraph \ref{lifthitch}-- and we want to find the elements $u$, $v$ and $w$ so that the $X_\alpha$ is the image by $u$ of the attractive flags of the corresponding monodromy etc ... This depends on the ordering of the eigenvalues  of $\Hol (A)$, $\Hol (B)$ and $\Hol (A)$.  In other words, if we denote by $\alpha^+$ the eigenvector with the largest eigenvalue of $\Hol (A)$, $\alpha^0$ the eigenvector with the middle eigenvalue of $\Hol (A)$ and $\alpha^-$ the eigenvector with the smallest eigenvalue of $\Hol (A)$, we need to  identify who, amongst $A_1,A_2,A_3$,  is $\alpha^+$. We will do that  using identities involving functions  of $x,y,a,A,b,B,c,C$,  in the next paragraph and write down the corresponding gap functions.

\subsection{Pant gap functions in dimension 3}
We have explained that
the ordering of the eigenvalues of $\Hol (A)$, $\Hol (B)$ and $\Hol (C)$ -- described by inequalities depending on  $x,y,a,A,b,B,c,C$ -- corresponds to  different choices of elements of the Weyl group on every boundary component, that is to  the various parametrisations of Hitchin representations. 

In this paragraph, we give the explicit expressions corresponding to various choices of parametrisation in the case of 
$\operatorname{PSL}(3,\mathbb R)$. The formulae were obtained using a computing software and the instructions can be found online on the ArXiv version.

We just want to make two informal observations
\begin{itemize}
\item The formulae are highly sensitive on the choice of elements of the Weyl group and it seems difficult to hope for a simple expression of the gap function in all possible cases.
\item  From the cluster algebra point of view, the gap functions are positive functions. Therefore the identities are relations in a completion of the cluster algebra. Does these identities have a pure cluster meaning ?  
\end{itemize}

\vskip 0.5truecm
\noindent{\sc The identity on the boundary}: 
The case $(u,,v,w)=(\operatorname{Id},\operatorname{Id},\operatorname{Id})$ is obtained whenever $(x,y,a,A,b,B,c,C)$  satisfies the following inequalities :
\begin{eqnarray*}
Bc>1>\frac{xy}{bC},\ \
\frac{xy}{Ac}>1>Ca,\ \ 
\frac{xy}{Ba}>1>bA.
\end{eqnarray*}
Indeed, in this case $X_\alpha=F(\alpha^+,\alpha^0)$, $X_\alpha^*=F(\alpha^+,\alpha^0)$,  $X_\gamma=F(\gamma^+,\gamma^0)$ and  $X_\beta=F(\beta^+,\beta^0)$.
Then 
$$
B_1=\beta^+,\ \ C_1=\gamma^-,\ A_3=\alpha^+, \ A_1=\alpha^-.
$$
The pant gap function is then the logarithm of
$$
\bb(B_1,A_1,C_1,C_3).
$$
Using a computing software whose  instructions can be found online on the ArXiv version, we obtained that the exponential of the pant gap function is equal to the following rational fraction
\begin{eqnarray}\label{comp:1}
e^{G(P)}= \Big((ayA + yA + xyA + xy + A^{2}a
b + abA + aybA + ybA) \cr
(xby + xy + xb + xbBC + xCb + bBx + b^{2}
BC + bBC)xy\Big)/\Big(bx^{3}y^{3} + bx^{3}y^{2} \cr
 + bx^{3}Cy^{2} + x^{2}y^{2}ACb + b^{2}B
x^{2}y^{2}C + bBx^{2}y^{2}C + bBx^{3}Cy^{2}
 + bBx^{3}y^{2} \cr
 + x^{3}y^{3} + 2b^{2}Bx^{2}y^{2}AC + b^{2}
Bx^{2}y^{2}AcC + b^{3}BxayAcC^{2} + x^{
2}y^{2}b^{2}ABaC \cr
 + xb^{3}AB^{2}acC^{2} + CbABax^{2}
y^{2} + x^{3}y^{2}ACb + xb^{3}A^{2}B^{2}acC
^{2} \cr
 + b^{4}A^{2}B^{2}acC^{2} + b^{2}AB^{2}a
cxyC + b^{2}ABaxyC + b^{2}Ax^{2}y^{2}C
 \cr
 + x^{2}y^{2}AcBCb + b^{3}A^{2}B^{2}ac
C^{2} + x^{3}y^{2}bA + x^{3}y^{2}bAB + 2xb^{2
}Bay^{2}AcC \cr
 + xb^{3}BA^{2}acCy + xy^{2}bBAC
 + xb^{3}By^{2}AcC + 2xb^{3}BaAcCy
 \cr
 + xb^{3}BaAcy^{2}C + xb^{3}AB^{2}a
cyC + xy^{2}b^{3}ABaC + 2xb^{2}By^{2}
AcC \cr
 + 2xb^{2}By^{2}AC + xb^{3}By^{2}AC
 + xb^{3}BA^{2}acC + xb^{3}A^{2}B^{2}acC
 \cr
 + xb^{3}AB^{2}acC + xb^{3}BaAcC^{
2} + xb^{3}BaAcC + xy^{2}bBaAcC \cr
 + 2xyb^{2}BaAcC + xyb^{2}BaAc
C^{2} + xyb^{3}ABaC + xyb^{3}BAcC \cr
 + xyb^{2}BA^{2}acC + xyb^{2}BAc
C + x^{3}y^{2}bBAC + 2x^{2}y^{2}bBAC \cr
 + x^{2}y^{2}b^{2}BA + x^{2}y^{2}bAB + x
y^{2}CbABa + xy^{2}bBAcC + xb^{3}BA^{
2}acC^{2} \cr
 + 2xy^{2}b^{2}ABaC\Big) .
\end{eqnarray}

\vskip 0.5truecm

\noindent{\sc The case of Theorem \ref{theo:goodgap}}, we have
$$
B_1=\beta^+,\ \ C_3=\gamma^-\ \ A_2=\alpha^+,\ \ A_3=\alpha^-.
$$
which gives the expected  reasonable formula 
\begin{eqnarray}\label{comp:2}
G(P)=\log\left(
\frac{1+B}{B(1+c)}
\right).
\end{eqnarray}

\vskip 0.5truecm

\noindent{\sc An intermediate choice} A somewhat intermediate choice is 
$$
B_1=\beta+,\ \ C_3=\gamma^-\ \ A_1=\alpha^+,\ \ A_3=\alpha^-.
$$
which gives
\begin{eqnarray}\label{comp:3}
G(P)=\log\left(\frac {(xby + xy + x
b + xbBC + xCb + bBx + b^{2}BC + bBC)y}{b
B(bcC^{2} + cCb + bcyC + yCb + cyC + y
C + Cxy + xy)}\right).
\end{eqnarray}

\section{Appendix B: computer instructions}

 \def\emptyline{\vspace{12pt}}
\DefineParaStyle{Heading 1}
\DefineParaStyle{Maple Output}
\DefineParaStyle{Warning}
\DefineCharStyle{2D Math}
\DefineCharStyle{2D Output}
\pagestyle{empty}
\begin{maplegroup}
\begin{mapleinput}
\mapleinline{active}{1d}{restart;}{%
}
\end{mapleinput}

\begin{mapleinput}
\mapleinline{active}{1d}{with(linalg):}{%
}
\end{mapleinput}

\mapleresult
\end{maplegroup}

\subsection{Matrices and holonomies}\label{comp:hol}

\begin{maplegroup}
\begin{mapleinput}
\mapleinline{active}{1d}{Id:=linalg[matrix](3,3,[1,0,0,0,1,0,0,0,1]):}{%
}
\end{mapleinput}

\begin{mapleinput}
\mapleinline{active}{1d}{Xf:=x->linalg[matrix](3,3,[0,0,1,0,-1,-1,x,1+x,1]):}{%
}
\end{mapleinput}

\begin{mapleinput}
\mapleinline{active}{1d}{Af:=(a,b)->linalg[matrix](3,3,[0,0,1/b,0,-1,0,a,0,0]):}{%
}
\end{mapleinput}

\begin{mapleinput}
\mapleinline{active}{1d}{X:=Xf(x):}{%
}
\end{mapleinput}

\begin{mapleinput}
\mapleinline{active}{1d}{Y:=Xf(y):}{%
}
\end{mapleinput}

\begin{mapleinput}
\mapleinline{active}{1d}{Qb:=Af(B,b):}{%
}
\end{mapleinput}

\begin{mapleinput}
\mapleinline{active}{1d}{Qa:=Af(A,a):}{%
}
\end{mapleinput}

\begin{mapleinput}
\mapleinline{active}{1d}{Qc:=Af(C,c):}{%
}
\end{mapleinput}

\begin{mapleinput}
\mapleinline{active}{1d}{HolA:=evalm(X&*inverse(Qc)&*Y&*Qb);}{%
}
\end{mapleinput}

\begin{mapleinput}
\mapleinline{active}{1d}{HolB:=evalm(X&*X&*inverse(Qa)&*Y&*Qc&*X&*X):}{%
}
\end{mapleinput}

\begin{mapleinput}
\mapleinline{active}{1d}{HolC:=evalm(inverse(Qb)&*Y&*Qa&*X);}{%
}
\end{mapleinput}

\mapleresult
\begin{maplelatex}
\mapleinline{inert}{2d}{HolA := matrix([[c*B, 0, 0], [(-c-1)*B, 1, 0], [(c+1+x+x/C)*B,
-1-x-x/C*(1+y), x/C*y/b]]);}{%
\[
\mathit{HolA} :=  \left[ 
{\begin{array}{ccc}
c\,B & 0 & 0 \\
( - c - 1)\,B & 1 & 0 \\
(c + 1 + x + {\displaystyle \frac {x}{C}} )\,B &  - 1 - x - 
{\displaystyle \frac {x\,(1 + y)}{C}}  & {\displaystyle \frac {x
\,y}{C\,b}} 
\end{array}}
 \right] 
\]
}
\end{maplelatex}

\begin{maplelatex}
\mapleinline{inert}{2d}{HolC := matrix([[1/B*y/a*x, 1/B*(1+y)+1/B*y/a*(1+x),
1/B*A+1/B*(1+y)+1/B*y/a], [0, 1, A+1], [0, 0, b*A]]);}{%
\[
\mathit{HolC} :=  \left[ 
{\begin{array}{ccc}
{\displaystyle \frac {y\,x}{B\,a}}  & {\displaystyle \frac {1 + y
}{B}}  + {\displaystyle \frac {y\,(1 + x)}{B\,a}}  & 
{\displaystyle \frac {A}{B}}  + {\displaystyle \frac {1 + y}{B}} 
 + {\displaystyle \frac {y}{B\,a}}  \\ [2ex]
0 & 1 & A + 1 \\
0 & 0 & b\,A
\end{array}}
 \right] 
\]
}
\end{maplelatex}

\end{maplegroup}

\subsection{Eigenvectors and eigenvalues}

\begin{maplegroup}
\begin{mapleinput}
\mapleinline{active}{1d}{IA:=inverse(transpose(HolA)):}{%
}
\end{mapleinput}

\begin{mapleinput}
\mapleinline{active}{1d}{egA:=eigenvalues(IA):}{%
}
\end{mapleinput}

\begin{mapleinput}
\mapleinline{active}{1d}{A3:=kernel(IA-scalarmul(Id,1/(c*B)))[1];}{%
}
\end{mapleinput}

\begin{mapleinput}
\mapleinline{active}{1d}{A2:=kernel(IA-scalarmul(Id,1))[1];}{%
}
\end{mapleinput}

\begin{mapleinput}
\mapleinline{active}{1d}{A1:=kernel(IA-scalarmul(Id,(C*b)/(x*y)))[1];}{%
}
\end{mapleinput}

\mapleresult
\begin{maplelatex}
\mapleinline{inert}{2d}{A3 := vector([1, 0, 0]);}{%
\[
\mathit{A3} := [1, \,0, \,0]
\]
}
\end{maplelatex}

\begin{maplelatex}
\mapleinline{inert}{2d}{A2 := vector([(c+1)*B/(c*B-1), 1, 0]);}{%
\[
\mathit{A2} :=  \left[  \! {\displaystyle \frac {(c + 1)\,B}{c\,B
 - 1}} , \,1, \,0 \!  \right] 
\]
}
\end{maplelatex}

\begin{maplelatex}
\mapleinline{inert}{2d}{A1 :=
vector([b*B*x*(b*c*C^2+c*C*b+b*c*y*C+y*C*b+c*y*C+y*C+C*x*y+x*y)/(c*B*C
^2*b^2-c*B*x*y*C*b-x*y*C*b+x^2*y^2), -b*(C+x*C+x+x*y)/(-C*b+x*y),
1]);}{%
\maplemultiline{
\mathit{A1} :=  \left[ {\vrule height0.80em width0em depth0.80em}
 \right. \!  \! {\displaystyle \frac {b\,B\,x\,(b\,c\,C^{2} + c\,
C\,b + b\,c\,y\,C + y\,C\,b + c\,y\,C + y\,C + C\,x\,y + x\,y)}{c
\,B\,C^{2}\,b^{2} - c\,B\,x\,y\,C\,b - x\,y\,C\,b + x^{2}\,y^{2}}
} ,  \\
 - {\displaystyle \frac {b\,(C + x\,C + x + x\,y)}{ - C\,b + x\,y
}} , \,1 \! \! \left. {\vrule height0.80em width0em depth0.80em}
 \right]  }
}
\end{maplelatex}

\end{maplegroup}
\begin{maplegroup}
\begin{mapleinput}
\mapleinline{active}{1d}{IB:=HolB:}{%
}
\end{mapleinput}

\begin{mapleinput}
\mapleinline{active}{1d}{B3:=kernel(IB-scalarmul(Id,x*C*a))[1];}{%
}
\end{mapleinput}

\begin{mapleinput}
\mapleinline{active}{1d}{B2:=kernel(IB-scalarmul(Id,x))[1];}{%
}
\end{mapleinput}

\begin{mapleinput}
\mapleinline{active}{1d}{B1:=kernel(IB-scalarmul(Id,(x*x*y)/(A*c)))[1];}{%
}
\end{mapleinput}

\begin{mapleinput}
\mapleinline{active}{1d}{egB:=eigenvalues(IB):}{%
}
\end{mapleinput}

\mapleresult
\begin{maplelatex}
\mapleinline{inert}{2d}{B3 :=
vector([-c*(a*A*x+x*y*a+x*C*a*A+a*x+x*y+x*C*a+C*a^2*A+C*a*A)/x/(a*c+y+
a*A*c+C*a*A*c+a*c*y+C*a*c+a*y+c*y), 1,
-(C^2*a*c+C*a*c+C*a*c*y+C*a*y+c*y*C+y*C+C*x*y+x*y)/(a*c+y+a*A*c+C*a*A*
c+a*c*y+C*a*c+a*y+c*y)/C]);}{%
\maplemultiline{
\mathit{B3} :=  \left[ {\vrule height0.80em width0em depth0.80em}
 \right. \!  \!  - {\displaystyle \frac {c\,(a\,A\,x + x\,y\,a + 
x\,C\,a\,A + a\,x + x\,y + x\,C\,a + C\,a^{2}\,A + C\,a\,A)}{x\,(
a\,c + y + a\,A\,c + C\,a\,A\,c + a\,c\,y + C\,a\,c + a\,y + c\,y
)}} , \,1,  \\
 - {\displaystyle \frac {C^{2}\,a\,c + C\,a\,c + C\,a\,c\,y + C\,
a\,y + c\,y\,C + y\,C + C\,x\,y + x\,y}{(a\,c + y + a\,A\,c + C\,
a\,A\,c + a\,c\,y + C\,a\,c + a\,y + c\,y)\,C}}  \! \! \left. 
{\vrule height0.80em width0em depth0.80em} \right]  }
}
\end{maplelatex}

\begin{maplelatex}
\mapleinline{inert}{2d}{B2 := vector([-(A+x*A+x+x*y)*c/x/(c+y+c*y+A*c), 1,
-(x*y+y+c*y+c)/(c+y+c*y+A*c)]);}{%
\[
\mathit{B2} :=  \left[  \!  - {\displaystyle \frac {(A + x\,A + x
 + x\,y)\,c}{x\,(c + y + c\,y + A\,c)}} , \,1, \, - 
{\displaystyle \frac {x\,y + y + c\,y + c}{c + y + c\,y + A\,c}} 
 \!  \right] 
\]
}
\end{maplelatex}

\begin{maplelatex}
\mapleinline{inert}{2d}{B1 := vector([-1, 1, -1]);}{%
\[
\mathit{B1} := [-1, \,1, \,-1]
\]
}
\end{maplelatex}

\end{maplegroup}
\begin{maplegroup}
\begin{mapleinput}
\mapleinline{active}{1d}{IC:=HolC:}{%
}
\end{mapleinput}

\begin{mapleinput}
\mapleinline{active}{1d}{C3:=kernel(IC-scalarmul(Id,(x*y)/(B*a)))[1];}{%
}
\end{mapleinput}

\begin{mapleinput}
\mapleinline{active}{1d}{C2:=kernel(IC-scalarmul(Id,1))[1];}{%
}
\end{mapleinput}

\begin{mapleinput}
\mapleinline{active}{1d}{C1:=kernel(IC-scalarmul(Id,(b*A)))[1];}{%
}
\end{mapleinput}

\begin{mapleinput}
\mapleinline{active}{1d}{egC:=eigenvalues(IC):}{%
}
\end{mapleinput}

\mapleresult
\begin{maplelatex}
\mapleinline{inert}{2d}{C3 := vector([1, 0, 0]);}{%
\[
\mathit{C3} := [1, \,0, \,0]
\]
}
\end{maplelatex}

\begin{maplelatex}
\mapleinline{inert}{2d}{C2 := vector([(a+a*y+y+x*y)/(B*a-x*y), 1, 0]);}{%
\[
\mathit{C2} :=  \left[  \! {\displaystyle \frac {a + a\,y + y + x
\,y}{B\,a - x\,y}} , \,1, \,0 \!  \right] 
\]
}
\end{maplelatex}

\begin{maplelatex}
\mapleinline{inert}{2d}{C1 :=
vector([(a*y*A+y*A+x*y*A+x*y+A^2*a*b+a*b*A+a*y*b*A+y*b*A)/(-x*y*A-x*y+
b*A^2*B*a+b*A*B*a), 1, (-1+b*A)/(A+1)]);}{%
\[
\mathit{C1} :=  \left[  \! {\displaystyle \frac {a\,y\,A + y\,A
 + x\,y\,A + x\,y + A^{2}\,a\,b + a\,b\,A + a\,y\,b\,A + y\,b\,A
}{ - x\,y\,A - x\,y + b\,A^{2}\,B\,a + b\,A\,B\,a}} , \,1, \,
{\displaystyle \frac { - 1 + b\,A}{A + 1}}  \!  \right] 
\]
}
\end{maplelatex}

\end{maplegroup}

\subsection{Pant Gap Function}\label{par:bir}

\begin{maplegroup}
\begin{mapleinput}
\mapleinline{active}{1d}{S:=(u,v)->multiply(transpose(convert(u,vector)),convert(v,vector)):}{
}
\end{mapleinput}

\begin{mapleinput}
\mapleinline{active}{1d}{bir:=(u,v,w,z)->eval(S(u,v)*S(w,z)/(S(u,z)*S(w,v))):}{%
}
\end{mapleinput}

\begin{mapleinput}
\mapleinline{active}{1d}{pantgap:=factor(simplify(evalm(bir(B1,A2,C3,A3))));}{%
}
\end{mapleinput}

\mapleresult
\begin{maplelatex}
\mapleinline{inert}{2d}{pantgap := (B+1)/(c+1)/B;}{%
\[
\mathit{birapport} := {\displaystyle \frac {B + 1}{(c + 1)\,B}} 
\]
}
\end{maplelatex}

\end{maplegroup}
\begin{maplegroup}
\begin{mapleinput}
\mapleinline{active}{1d}{pantgap:=factor(simplify(evalm(bir(B1,A1,C1,A3))));}{%
}
\end{mapleinput}

\mapleresult
\begin{maplelatex}
\mapleinline{inert}{2d}{pantgap :=
(a*y*A+y*A+x*y*A+x*y+A^2*a*b+a*b*A+a*y*b*A+y*b*A)*(x*b*y+x*y+x*b+x*b*B
*C+x*C*b+b*B*x+b^2*B*C+b*B*C)*x*y/(b*x^3*y^3+b*x^3*y^2+b*x^3*C*y^2+x^2
*y^2*A*C*b+b^2*B*x^2*y^2*C+b*B*x^2*y^2*C+b*B*x^3*C*y^2+b*B*x^3*y^2+x^3
*y^3+2*b^2*B*x^2*y^2*A*C+b^2*B*x^2*y^2*A*c*C+b^3*B*x*a*y*A*c*C^2+x^2*y
^2*b^2*A*B*a*C+x*b^3*A*B^2*a*c*C^2+C*b*A*B*a*x^2*y^2+x^3*y^2*A*C*b+x*b
^3*A^2*B^2*a*c*C^2+b^4*A^2*B^2*a*c*C^2+b^2*A*B^2*a*c*x*y*C+b^2*A*B*a*x
*y*C+b^2*A*x^2*y^2*C+x^2*y^2*A*c*B*C*b+b^3*A^2*B^2*a*c*C^2+x^3*y^2*b*A
+x^3*y^2*b*A*B+2*x*b^2*B*a*y^2*A*c*C+x*b^3*B*A^2*a*c*C*y+x*y^2*b*B*A*C
+x*b^3*B*y^2*A*c*C+2*x*b^3*B*a*A*c*C*y+x*b^3*B*a*A*c*y^2*C+x*b^3*A*B^2
*a*c*y*C+x*y^2*b^3*A*B*a*C+2*x*b^2*B*y^2*A*c*C+2*x*b^2*B*y^2*A*C+x*b^3
*B*y^2*A*C+x*b^3*B*A^2*a*c*C+x*b^3*A^2*B^2*a*c*C+x*b^3*A*B^2*a*c*C+x*b
^3*B*a*A*c*C^2+x*b^3*B*a*A*c*C+x*y^2*b*B*a*A*c*C+2*x*y*b^2*B*a*A*c*C+x
*y*b^2*B*a*A*c*C^2+x*y*b^3*A*B*a*C+x*y*b^3*B*A*c*C+x*y*b^2*B*A^2*a*c*C
+x*y*b^2*B*A*c*C+x^3*y^2*b*B*A*C+2*x^2*y^2*b*B*A*C+x^2*y^2*b^2*B*A+x^2
*y^2*b*A*B+x*y^2*C*b*A*B*a+x*y^2*b*B*A*c*C+x*b^3*B*A^2*a*c*C^2+2*x*y^2
*b^2*A*B*a*C);}{%
\maplemultiline{
\mathit{pantgap} := (a\,y\,A + y\,A + x\,y\,A + x\,y + A^{2}\,a
\,b + a\,b\,A + a\,y\,b\,A + y\,b\,A) \\
(x\,b\,y + x\,y + x\,b + x\,b\,B\,C + x\,C\,b + b\,B\,x + b^{2}\,
B\,C + b\,B\,C)\,x\,y/(b\,x^{3}\,y^{3} + b\,x^{3}\,y^{2} \\
\mbox{} + b\,x^{3}\,C\,y^{2} + x^{2}\,y^{2}\,A\,C\,b + b^{2}\,B\,
x^{2}\,y^{2}\,C + b\,B\,x^{2}\,y^{2}\,C + b\,B\,x^{3}\,C\,y^{2}
 + b\,B\,x^{3}\,y^{2} \\
\mbox{} + x^{3}\,y^{3} + 2\,b^{2}\,B\,x^{2}\,y^{2}\,A\,C + b^{2}
\,B\,x^{2}\,y^{2}\,A\,c\,C + b^{3}\,B\,x\,a\,y\,A\,c\,C^{2} + x^{
2}\,y^{2}\,b^{2}\,A\,B\,a\,C \\
\mbox{} + x\,b^{3}\,A\,B^{2}\,a\,c\,C^{2} + C\,b\,A\,B\,a\,x^{2}
\,y^{2} + x^{3}\,y^{2}\,A\,C\,b + x\,b^{3}\,A^{2}\,B^{2}\,a\,c\,C
^{2} \\
\mbox{} + b^{4}\,A^{2}\,B^{2}\,a\,c\,C^{2} + b^{2}\,A\,B^{2}\,a\,
c\,x\,y\,C + b^{2}\,A\,B\,a\,x\,y\,C + b^{2}\,A\,x^{2}\,y^{2}\,C
 \\
\mbox{} + x^{2}\,y^{2}\,A\,c\,B\,C\,b + b^{3}\,A^{2}\,B^{2}\,a\,c
\,C^{2} + x^{3}\,y^{2}\,b\,A + x^{3}\,y^{2}\,b\,A\,B + 2\,x\,b^{2
}\,B\,a\,y^{2}\,A\,c\,C \\
\mbox{} + x\,b^{3}\,B\,A^{2}\,a\,c\,C\,y + x\,y^{2}\,b\,B\,A\,C
 + x\,b^{3}\,B\,y^{2}\,A\,c\,C + 2\,x\,b^{3}\,B\,a\,A\,c\,C\,y
 \\
\mbox{} + x\,b^{3}\,B\,a\,A\,c\,y^{2}\,C + x\,b^{3}\,A\,B^{2}\,a
\,c\,y\,C + x\,y^{2}\,b^{3}\,A\,B\,a\,C + 2\,x\,b^{2}\,B\,y^{2}\,
A\,c\,C \\
\mbox{} + 2\,x\,b^{2}\,B\,y^{2}\,A\,C + x\,b^{3}\,B\,y^{2}\,A\,C
 + x\,b^{3}\,B\,A^{2}\,a\,c\,C + x\,b^{3}\,A^{2}\,B^{2}\,a\,c\,C
 \\
\mbox{} + x\,b^{3}\,A\,B^{2}\,a\,c\,C + x\,b^{3}\,B\,a\,A\,c\,C^{
2} + x\,b^{3}\,B\,a\,A\,c\,C + x\,y^{2}\,b\,B\,a\,A\,c\,C \\
\mbox{} + 2\,x\,y\,b^{2}\,B\,a\,A\,c\,C + x\,y\,b^{2}\,B\,a\,A\,c
\,C^{2} + x\,y\,b^{3}\,A\,B\,a\,C + x\,y\,b^{3}\,B\,A\,c\,C \\
\mbox{} + x\,y\,b^{2}\,B\,A^{2}\,a\,c\,C + x\,y\,b^{2}\,B\,A\,c\,
C + x^{3}\,y^{2}\,b\,B\,A\,C + 2\,x^{2}\,y^{2}\,b\,B\,A\,C \\
\mbox{} + x^{2}\,y^{2}\,b^{2}\,B\,A + x^{2}\,y^{2}\,b\,A\,B + x\,
y^{2}\,C\,b\,A\,B\,a + x\,y^{2}\,b\,B\,A\,c\,C + x\,b^{3}\,B\,A^{
2}\,a\,c\,C^{2} \\
\mbox{} + 2\,x\,y^{2}\,b^{2}\,A\,B\,a\,C) }
}
\end{maplelatex}

\end{maplegroup}
\begin{maplegroup}
\begin{mapleinput}
\mapleinline{active}{1d}{birapport:=factor(simplify(evalm(bir(B1,A1,C3,A3))));}{%
}
\end{mapleinput}

\mapleresult
\begin{maplelatex}
\mapleinline{inert}{2d}{pantgap :=
(x*b*y+x*y+x*b+x*b*B*C+x*C*b+b*B*x+b^2*B*C+b*B*C)*y/b/B/(b*c*C^2+c*C*b
+b*c*y*C+y*C*b+c*y*C+y*C+C*x*y+x*y);}{%
\[
\mathit{pantgap} := {\displaystyle \frac {(x\,b\,y + x\,y + x\,
b + x\,b\,B\,C + x\,C\,b + b\,B\,x + b^{2}\,B\,C + b\,B\,C)\,y}{b
\,B\,(b\,c\,C^{2} + c\,C\,b + b\,c\,y\,C + y\,C\,b + c\,y\,C + y
\,C + C\,x\,y + x\,y)}} 
\]
}
\end{maplelatex}

\end{maplegroup}
\vskip 0.5 truecm
\section{Appendix C: positivity of Frenet curves}

\subsection{Positivity of triples} We first prove.

\begin{proposition} \label{disjoint intervals}
Let $\xi^1$ be a Frenet curve. Let $I_1$, $I_2$ and $I_3$ be  three disjoint subintervals of $\torus$.
For  $i=1,2,3$, let  $X_i=(x^i_1,\ldots,x^i_{n-1})$  be an $(n-1)$-tuple of distinct points  of $I_i$. 
 Let $F_i$ be the flag given by
$$
F_i^k=\sum_{j=1}^{k}\xi^1 (x^i_j).
$$
Then the triple $(F_1,F_2,F_3)$ is positive.
\end{proposition}

\proof  Let $t$ be a point in $\torus\setminus (I_1\cup I_2\cup I_3)$. 
We choose  a continuous map $\widehat\xi$ 
$$
\mapping{\torus\setminus\{ t\}}{\mathbb R^n\setminus\{0\},}{u}{\widehat\xi(u)\in\xi(u).}
$$ 
We choose an orientation $\Omega$ on $\mathbb R^n$ 
such that for any positively oriented 
 $n$-tuple of points $\{y_1,\ldots,y_n\}$ in $\torus\setminus\{t\}$ we have 
\begin{eqnarray}
\Omega(\widehat\xi(y_1) \wedge \ldots \wedge \widehat\xi(y_n))>0.\label{omega} 
\end{eqnarray}
We write
$$
X_i^p=\widehat{\xi}^1(x^i_1)\wedge\ldots\wedge\widehat\xi^1(x^i_p). 
$$
For   each  $i, 1 \leq i \leq 3$ and each $p, 1 \leq p \leq n$
let   $\sigma_{i,p}$  be a permutation such that  
$
(x^i_{\sigma_{i,p}(1)},\ldots,x^i_{\sigma_{i,p}(p)}),
$
is positively oriented. 

We also define $\alpha(I_1,I_2,I_3)=+1$, if $(a_1,a_2,a_3)$ is positively oriented for any $a_i$ in $I_i$ and 
 $\alpha(I_1,I_2,I_3)=-1$ in the opposite case.
Thus the sign of 
$$
\Omega( X_1^m\wedge X_2^l\wedge X_3^p),
$$
coincides with the sign of 
$$
\varepsilon(\sigma_{1,m}).\varepsilon(\sigma_{2,l}).\varepsilon(\sigma_{3,p}))\alpha(I_1,I_2,I_3)^{pl},
$$
where $\varepsilon(\sigma)$ is the signature of the permutation $\sigma$.

When we report this information in the ratio given in Definition \ref{triratio}, we obtain that each of the signatures above appear twice, and thus the sign of $T^{m,l,p}(F_1,F_2,F_3)$
coincides with the sign of 
$$
\alpha(I_1,I_2,I_3)^{(p-1)l + p(l+1) + (p+1)(l+1)-p(l-1)-(p-1)(l+1)-(p+1)l}=\alpha(I_1,I_2,I_3)^{2p+2}.$$

Thus  the triple $(F_1,F_2,F_3)$ is positive.\qed

The following is essentially a corollary
\begin{proposition}\label{postrip2}
Let  $\xi^1$ be  a Frenet curve from $\torus$ to $\RPN$  with osculating flag curve 
 $\xi=(\xi^{1},\xi^{2},\ldots,\xi^{n-1})$.
 Then  $(\xi(y_1),\xi(y_2),\xi(y_3))$ is a positive triple of flags 
 whenever  $(y_1,y_2,y_3)$ is a  triple of distinct points of $\torus$.

More generally, if
$
(y_1^+,y_1^-,y_2^+,y_2^-,y_3^+,y_3^-)
$
is a positively 
oriented sextuplet  of points in $S^1$ and $Y_i,\,i= 1,2,3$  is a flag compatible with 
$(\xi(y_i^+),\xi(y_i^-))$ then $(Y_1,Y_2,Y_3)$ is a positive flag.
\end{proposition}
\proof For  the first part, let $I_1$, $I_2$ and $I_3$ be three disjoint subintervals of $\torus$  
and   $X_i=(x^i_1,\ldots,x^i_{n-1})$ for $i=1,2,3$ be  three $(n-1)$-tuples of distinct points  $X_i\subset I_i$. 
As before set 
$$
F_i^k=\sum_{j=1}^{k}\xi^1 (x^i_j).
$$
It follows that the triple of flags $(F_1,F_2,F_3)$ satisfy the hypothesis of  Proposition \ref{disjoint intervals}
above. Hence, this triple is positive. 

Now we let $x^i_j$ tend to $y_j$.
Recall that $\xi$ is a Frenet curve and so satisfies Conditions (\ref{fre2})  and (\ref{fre3}).
Firstly, since $\xi$   satisfies  Condition (\ref{fre2})   by Definition \ref{triratio}, 
all the triple ratios associated to  $(\xi(y_1),\xi(y_2),\xi(y_3))$ are nonzero.
Secondly,  since  $\xi$   satisfies   Condition (\ref{fre3}), $\xi(y_i)$ is a limit of $F_i$, 
hence all the triple ratios of the 
$(\xi(y_1),\xi(y_2),\xi(y_3))$ are nonnegative. As a conclusion, $(\xi(y_1),\xi(y_2),\xi(y_3))$ is a positive triple.

The proof of the second part is a natural extension of this argument. 
One  merely observes that a  flag compatible with $(\xi(x_+),\xi(x_-))$ is a limit of flags  constructed from  direct sums  
of $\xi(x^i_+)$ and $\xi(x^j_-)$ for points $x^i_+$ close to $x_+$ and $x^j_-$ close to $x_-$. \qed

\subsection{Positivity of quadruples} 
The same argument as in the previous paragraph yields
\begin{proposition}\label{posquad2}
Let $\xi^1$ be a Frenet curve. Let $I_1$, $I_2$, $I_3$ and $I_4$
be  four disjoint subintervals of $\torus$ 
such that some (and so any) quadruple of points $t_{k}\in I_{k}$ is cyclically ordered. 
For $i=1,2,3,4$  let $X_i=(x^i_1,\ldots,x^i_{n-1})$  be an  $(n-1)$-tuple 
of distinct points  of  $ I_i$.
 Let $F_i$ be the flag given by
$$
F^k_i=\sum_{j=1}^{k}\xi^1(x_j^{i}). 
$$
Then the quadruple $(F_1,F_2,F_3,F_4)$ is positive.
\end{proposition}
\proof   By construction, the sign of $\delta_i(F_1,F_2,F_3,F_4)$ depends continuously on $(X_1,X_2,X_3,X_4)$. Moreover, when $x_1^1$ converges to $x^4_1$, $\delta_i(F_1,F_2,F_3,F_4)$ 
 converges to 1 by Equation (\ref{deltabir}). 
 The proposition follows from this remark and Proposition  \ref{disjoint intervals}. \qed
\begin{proposition}\label{posquadflag}
Let $\xi^1$ be a Frenet curve from $S^1$ to $\mathbb P(\mathbb R^{n})$  with osculating flag curve
 $\xi=(\xi^{1},\xi^{2},\ldots,\xi^{n-1})$. 
 Then for every positively oriented quadruple of distinct points $(y_1,y_2,y_3,y_4)$ in $S^1$,
  $(\xi(y_1),\xi(y_2),\xi(y_3),\xi(y_4))$ is a positive quadruple of flags.
More generally, let 
$$
(y_1^+,y_1^-,y_2^+,y_2^-,y_3^+,y_3^-,y_4^+,y_4^-),
$$
be a cyclically ordered octuplet of points in $S^1$. Let $Y_i$  be a flag compatible with 
$(\xi(y_i^+),\xi(y_i^-))$, then $(Y_1,Y_2,Y_3,Y_4)$ is a positive quadruple.
\end{proposition}

\proof We use a similar argument as in Proposition \ref{postrip2} in order to obtain the flags $Y_i$ as limits of flags of the type described in Proposition \ref{posquad2}.  Then, the edge functions for the flags $Y_i$ are nonnegative by a limiting argument, and nonzero by the Frenet property. The result follows. \qed
%

\begin{thebibliography}{10}

\bibitem{Akiyoshi:2006}
Hirotaka Akiyoshi, Hideki Miyachi, and Makoto Sakuma, \emph{Variations of
  {M}c{S}hane's identity for punctured surface groups}, Spaces of Kleinian
  groups, London Math. Soc. Lecture Note Ser., vol. 329, Cambridge Univ. Press,
  Cambridge, 2006, pp.~151--185.

\bibitem{Birman:1985}
Joan~S. Birman and Caroline Series, \emph{Geodesics with bounded intersection
  number on surfaces are sparsely distributed}, Topology \textbf{24} (1985),
  no.~2, 217--225.

\bibitem{Bonahon:1996}
Francis Bonahon, \emph{Shearing hyperbolic surfaces, bending pleated surfaces
  and {T}hurston's symplectic form}, Ann. Fac. Sci. Toulouse Math. (6)
  \textbf{5} (1996), no.~2, 233--297.

\bibitem{Bonahon:2001}
Francis Bonahon and Ya{\c{s}}ar S{\"o}zen, \emph{The {W}eil-{P}etersson and
  {T}hurston symplectic forms}, Duke Math. J. \textbf{108} (2001), no.~3,
  581--597.

\bibitem{Bourdon:1996}
Marc Bourdon, \emph{Sur le birapport au bord des {${\rm CAT}(-1)$}-espaces},
  Inst. Hautes \'Etudes Sci. Publ. Math. (1996), no.~83, 95--104.

\bibitem{Bowditch:1996}
Brian~H. Bowditch, \emph{A proof of {M}c{S}hane's identity via {M}arkoff
  triples}, Bull. London Math. Soc. \textbf{28} (1996), no.~1, 73--78.

\bibitem{Fock:2006a}
Vladimir~V. Fock and Alexander~B. Goncharov, \emph{Moduli spaces of local
  systems and higher {T}eichm{\"u}ller theory}, Inst. Hautes \'Etudes Sci.
  Publ. Math. (2006), no.~103, 1--211.

\bibitem{Fock:2007}
\bysame, \emph{Moduli spaces of convex projective structures and surfaces},
  Adv. in Math. \textbf{208} (2007), 249--273.

\bibitem{Hamenstadt:1997}
Ursula Hamenst{\"a}dt, \emph{Cocycles, {H}ausdorff measures and cross ratios},
  Ergodic Theory Dynam. Systems \textbf{17} (1997), no.~5, 1061--1081.

\bibitem{Labourie:2005a}
Fran{\c{c}}ois Labourie, \emph{Cross ratios, {A}nosov representations and the
  energy functional on {T}eichm{\"u}ller space.}, Ann. Sci. \'Ecole Norm. Sup.
  (4) (To appear).

\bibitem{Labourie:2006}
\bysame, \emph{Anosov flows, surface groups and curves in projective space},
  Invent. Math. \textbf{165} (2006), no.~1, 51--114.

\bibitem{Labourie:2005}
\bysame, \emph{Cross ratios, surface groups, {${\rm PSL}(n,{\bf R})$} and
  diffeomorphisms of the circle}, Publ. Math. Inst. Hautes \'Etudes Sci.
  (2007), no.~106, 139--213.

\bibitem{Ledrappier:1995}
Fran{\c{c}}ois Ledrappier, \emph{Structure au bord des vari\'et\'es \`a
  courbure n\'egative}, S\'eminaire de Th\'eorie Spectrale et G\'eom\'etrie,
  No. 13, Ann\'ee 1994--1995, S\'emin. Th\'eor. Spectr. G\'eom., vol.~13, Univ.
  Grenoble I, Saint, 1995, pp.~97--122.

\bibitem{Lusztig:1994}
George Lusztig, \emph{Total positivity in reductive groups}, Lie theory and
  geometry, Progr. Math., vol. 123, Birkh\"auser Boston, Boston, MA, 1994,
  pp.~531--568.

\bibitem{McShane:1998}
Gregory McShane, \emph{Simple geodesics and a series constant over
  {T}eichm{{\"u}}ller space}, Invent. Math. \textbf{132} (1998), no.~3,
  607--632.

\bibitem{Mirzakhani:2007a}
Maryam Mirzakhani, \emph{Simple geodesics and {W}eil-{P}etersson volumes of
  moduli spaces of bordered {R}iemann surfaces}, Invent. Math. \textbf{167}
  (2007), no.~1, 179--222.

\bibitem{Otal:1990}
Jean-Pierre Otal, \emph{Le spectre marqu\'e des longueurs des surfaces \`a
  courbure n\'egative}, Ann. of Math. (2) \textbf{131} (1990), no.~1, 151--162.

\bibitem{Otal:1992}
\bysame, \emph{Sur la g\'eometrie symplectique de l'espace des g\'eod\'esiques
  d'une vari\'et\'e \`a courbure n\'egative}, Rev. Mat. Iberoamericana
  \textbf{8} (1992), no.~3, 441--456.

\bibitem{Tan:2004}
Ser~Peow Tan, Yan~Loi Wong, and Ying Zhang, \emph{Mc{S}hane's identities for
  classical {S}chottky groups}, ArXiv: math.GT/0411628, 2004.

\bibitem{Thurston:1984}
William~P. Thurston, \emph{Minimal stretch maps between hyperbolic surfaces},
  ArXiv: math.GT/9801039, 1984.

\end{thebibliography}
\def\cprime{$'$}
\providecommand{\bysame}{\leavevmode\hbox to3em{\hrulefill}\thinspace}
\providecommand{\MR}{\relax\ifhmode\unskip\space\fi MR }
\providecommand{\MRhref}[2]{%
  \href{http://www.ams.org/mathscinet-getitem?mr=#1}{#2}
}
\providecommand{\href}[2]{#2}

\end{document}